\documentclass[12pt]{article}
\usepackage{amssymb}
\usepackage{mathrsfs}
\usepackage[top=72pt,bottom=96pt,left=72pt,right=72pt]{geometry}
\def\ad{\mathrm{ad}}
\def\alg{\mathrm{alg}}
\def\C{\mathbb{C}}
\def\cl{\mathrm{cl}}
\def\crit{\mathrm{crit}}
\def\D{\mathscr{L}}
\def\Der{\mathrm{Der}}
\def\dim{\mathrm{dim}}
\def\dist{\mathrm{dist}}

\def\ext{\mathrm{ext}}
\def\End{\mathrm{End}\,}
\def\FS{\mathrm{FS}}
\def\free{\mathrm{free}}
\def\g{\mathfrak{g}}
\def\Gr{\mathrm{Gr}}
\def\H{\mathscr{H}}
\def\Hess{\mathrm{Hess}\,}
\def\Hom{\mathrm{Hom}\,}
\def\id{\mathrm{id}}
\def\Im{\mathrm{Im}}
\def\im{\mathrm{im}}
\def\ins{\lrcorner}
\def\ker{\mathrm{ker}}
\def\Killing{\mathrm{Killing}}
\def\knc{\mathrm{knc}}
\def\Lam{\Lambda}
\def\LGr{\mathrm{LGr}}
\def\loc{\mathrm{loc}}
\def\N{\mathbb{N}}
\def\ort{\mathrm{or}}
\def\pfill{\par\vskip10pt plus3pt minus3pt\noindent}
\def\pr{\mathrm{pr}}
\def\proof{\pfill\textbf{Proof:}\quad}
\def\PT{\mathbf{PT}}
\def\Q{\mathbb{Q}}
\def\qed{\ensuremath{\hfill\Box}}
\def\R{\mathbb{R}}
\def\Re{\mathrm{Re}}
\def\rest{\mathrm{rest}}
\def\S{\mathrm{Sym}}
\def\skew{\mathrm{skew}}
\def\stab{\mathfrak{stab}}
\def\T{\mathfrak{T}}
\def\tr{\mathrm{tr}}
\def\X{\mathfrak{X}}
\def\Y{\mathfrak{Y}}
\newtheorem{Lemma}{Lemma}[section]
\newtheorem{Remark}[Lemma]{Remark}
\newtheorem{Theorem}[Lemma]{Theorem}

\newtheorem{Corollary}[Lemma]{Corollary}
\newtheorem{Definition}[Lemma]{Definition}

\newenvironment{Reference}[1]{\pfill\textbf{#1} \textit\bgroup}{\egroup\par}
\begin{document}
\title{Riemannian and K\"ahlerian Normal Coordinates}
\author{Tillmann Jentsch%
 \footnote{Lehrstuhl f\"ur Geometrie, Institut f\"ur Geometrie und Topologie,
 Fachbereich Mathematik, Universit\"at Stuttgart, Pfaffenwaldring 57,
 70569 Stuttgart, ALLEMAGNE;\ \texttt{tilljentsch@web.de}}
 \quad\&\quad Gregor Weingart%
 \footnote{Unidad Cuernavaca del Instituto de Matem\'aticas, Universidad
 Nacional Aut\'onoma de M\'exico, Avenida Universidad s/n, Lomas de
 Chamilpa, 62210 Cuernavaca, Morelos, MEXIQUE;\ \texttt{gw@matcuer.unam.mx}}}
\maketitle
\begin{center}
 \textbf{Abstract}
 \\[11pt]
 \parbox{400pt}{%
  In every point of a K\"ahler manifold there exist special holomorphic
  coordinates well adapted to the underlying geometry. Comparing these
  K\"ahler normal coordinates with the Riemannian normal coordinates defined
  via the exponential map we prove that their difference is a universal power
  series in the curvature tensor and its iterated covariant derivatives and
  devise an algorithm to calculate this power series to arbitrary order. As
  a byproduct we generalize K\"ahler normal coordinates to the class of
  complex affine manifolds with $(1,1)$--curvature tensor. Moreover we
  describe the Spencer connection on the infinite order Taylor series of
  the K\"ahler normal potential and obtain explicit formulas for the Taylor
  series of all relevant geometric objects on symmetric spaces.}
 \\[11pt]
\end{center}
\section{Introduction}
 Situated at the crossroads of differential and algebraic geometry the
 geometry of K\"ahler manifolds is a very attractive topic of research,
 in particular both analytical and algebraic tools and ideas have been
 brought to bear on the topic. The main motivation of the article at
 hand is to understand the interrelationship between the Taylor series
 of infinite order of several objects relevant for K\"ahler geometry.
 Philosophically our study of K\"ahler manifolds is thus based in the
 jet calculus of differential geometry, a calculus which in contrast to
 tensor or exterior calculus tries to avoid taking actual derivatives
 at all cost. Instead of taking derivatives jet calculus focusses on the
 algebraic constraints satisfied by the jets of all relevant objects, ideally
 then these algebraic constraints are sufficient to determine all the jets
 by multiplication, which at higher orders is much simpler than
 differentiation.

 \pfill
 Not to the least the success of this strategy depends on our ability to
 isolate the objects to which it can be applied in the first place. In
 affine geometry for example the relevant object turns out to be the backward
 parallel transport $\Phi^{-1}$, which encodes the differential of the
 exponential map and describes the Taylor series of all covariantly parallel
 tensors. The parallel transport equation proved in Lemma \ref{pte} of this
 article is the algebraic constraint
 $$
  N\,(\,N\,+\,1\,)\,\Phi^{-1}
  \;\;=\;\;
  \mathscr{R}\,\Phi^{-1}\;+\;N\,(\,\mathscr{T}\,\Phi^{-1}\,)
 $$
 on the Taylor series of the backward parallel transport, whose solution
 $\Phi^{-1}$ is easily found by multiplication. In passing we remark that
 in the torsion free case the Taylor series of $\Phi^{-1}$ was calculated
 by several authors before using either implicitly or explicitly properties
 of iterated covariant derivatives. The novelity in our argument is that
 the parallel transport equation appears directly as an integrability
 condition for the Jacobi equation.

 In order to find the relevant objects for the study of K\"ahler geometry
 we recall an observation originally due to Bochner \cite{bo}: Riemannian
 normal coordinates on K\"ahler manifolds are not complex coordinates unless
 the manifold is flat. Bochner found remedy to this nuisance in singling out
 complex coordinates centered in an arbitrary point $p$ of a K\"ahler manifold
 $M$, which are unique up to unitary transformations. Their lack of uniqueness
 is easily fixed by thinking of K\"ahler normal coordinates as anchored
 coordinates, local diffeomorphisms
 $$
  \knc_p:\;\;T_pM\;\longrightarrow\;M,\qquad X\;\longmapsto\;\knc_pX
 $$
 mapping the origin to $p$ with differential $(\knc_p)_{*,0}\,=\,\id_{T_pM}$
 under the usual identification $T_0(T_pM)\,\cong\,T_pM$. The difference
 element $K\,:=\,\exp_p^{-1}\circ\knc_p$ and its inverse $K^{-1}$ are
 ``hidden'' relevant objects in K\"ahler geometry, because their Taylor
 series relate to commonly considered objects like the K\"ahler potential.
 Inevitably then a pivotal role in this article is played by the proof
 of the following statement about the Taylor series of $K$ and $K^{-1}$:

 \begin{Reference}{Theorem \ref{univ}\ 
  (Universality of K\"ahler Normal Coordinates)}\hfill\break
  Every term in the Taylor series of the difference element $K\,=\,
  \exp_p^{-1}\circ\knc_p$ and its inverse $K^{-1}$ in the origin of K\"ahler
  normal coordinates is a universal polynomial, independent of the manifold
  and its dimension, in the complex structure $I_p$, the curvature tensor
  $R_p$ and its iterated covariant derivatives $(\nabla R)_p,\,(\nabla^2R)_p,
  \,\ldots$ evaluated at $p\,\in\,M$.
 \end{Reference}

 \pfill
 Of course the theorem in itself can hardly be called surprising, its
 validity for example is implicitly taken for granted without even a
 fleeting comment in \cite{hin}. The rigorous proof in Section \ref{element}
 however allows us to draw serveral conclusions about the nature of K\"ahler
 normal coordinates, for example K\"ahler normal coordinates are inherited
 by totally geodesic complex submanifolds according to Corollary \ref{sub}.
 Moreover the proof of Theorem \ref{univ} provides us with a simple recursion
 formula for the Taylor series of the difference element $K^{-1}$.

 Perhaps the most striking insight however to be taken from the proof of
 Theorem \ref{univ} is that K\"ahler normal coordinates, very much like
 their Riemannian counterpart, do only depend on the complex affine
 geometry and {\em not} on the metric structure, despite the fact that
 their characterization apparently involves the metric in form of the
 potential. In fact K\"ahler normal coordinates generalize to the much
 more general class of balanced complex affine manifolds: Complex manifolds
 $M$ endowed with a torsion free connection $\nabla$ on their tangent
 bundle $TM$ such that the associated almost complex structure is parallel
 $\nabla I\,=0$ and the curvature tensor is a $(1,1)$--form in the
 sense $R_{IX,\,IY}\,=\,R_{X,\,Y}$ for all $X,\,Y$.

 \pfill
 Using the recursion formula for $K^{-1}$ formulated in Remark \ref{recs}
 in analogy to the parallel transport equation we calculate the initial terms
 of the Taylor series of $K^{-1}$, of the K\"ahler normal potential and the
 Riemannian distance to the origin. Actually the recursion formula is simple
 enough to be readily programmed in a computer algebra system, the problem
 is to standardize the numerous terms cropping up in the process. An
 interesting conclusion from these calculations is in any case that the
 total holomorphic sectional curvature
 $$
  S_{\mathrm{total}}(\;X\;)
  \;\;:=\;\;
  \sum_{k\,\geq\,4}\frac1{(k-4)!}\;
  g(\;(\nabla^{k-4}_{X,\,\ldots,\,X}R)_{X,\,IX}IX,\;X\;)
  \;\;\in\;\;
  \Gamma(\;\overline\S^{\geq(2,2)}T^*M\;)
 $$
 written as a sum of its bihomogeneous components $S_{\kappa,\,\overline
 \kappa}$ of degrees $\kappa,\,\overline\kappa\,\geq\,2$ is congruent
 \begin{eqnarray}
  \theta(\,X\,)\qquad
  &\equiv&
  g(\,X,\,X\,)\,-\,\sum_{\kappa,\,\overline\kappa\,\geq\,2}
  \frac1{2\,\kappa\,(\kappa-1)\,\overline\kappa\,(\overline\kappa-1)}\;
  S_{\kappa,\,\overline\kappa}(\,X\,)
  \label{tcong}
  \\[3pt]
  \dist^2_g(p,\knc_pX)
  &\equiv&
  g(\,X,\,X\,)\,-\,\sum_{\kappa,\,\overline\kappa\,\geq\,2}
  \frac1{2\,(\kappa+\overline\kappa-1)\,(\kappa-1)\,(\overline\kappa-1)}\;
  S_{\kappa,\,\overline\kappa}(\;X\;)
  \label{dcong}
 \end{eqnarray}
 to the K\"ahler normal potential $\theta$ and the Riemannian distance
 to the origin modulo terms at least {\em quadratic} in the curvature tensor
 and all its iterated covariant derivatives. A simple induction based on
 these congruences implies that the total holomorphic sectional curvature,
 the K\"ahler normal potential and the Riemannian distance to the origin
 all parametrize the underlying K\"ahler geometry uniquely up to covering.

 For the time being the K\"ahler normal potential is certainly the most
 convenient of these three parameters, because the covariant derivative
 of the K\"ahler normal potential considered as a section $\theta\,\in\,
 \Gamma(\,\overline\S\,T^*M\,)$ can be calculated using the concept of
 holomorphically extended vector fields. More precisely we find in Lemma
 \ref{spcon} an explicit formula of the form
 $$
  \nabla_Z\theta
  \;\;=\;\;
  \pr_{\geq(2,2)}\,(\;Z\;\ins\;\theta^\free\;)
  \;-\;\pr_{[\,1\,]}(\;Z\;\ins\;\theta^\crit\;)\;\bullet\;\theta^\free
 $$
 where $\bullet$ is a simple bilinear operation, while $\theta^\free$ and
 $\theta^\crit$ denote the sums of all bihomogeneous components of $\theta$
 at least or exactly quadratic respectively in the holomorphic or
 antiholomorphic coordinates. Note that this result is best thought
 of as describing the so called Spencer connection on the parameter
 vector bundle $\overline\S^{\geq(2,2)}T^*M$.

 \pfill
 Specializing from general K\"ahler manifolds to concrete examples we
 study the locally symmetric K\"ahler manifolds, usually called hermitean
 locally symmetric spaces, in the second part of this article. In general
 the Lie theoretic setup \cite{s} \cite{br} makes many calculations like
 the determination of the Taylor series of the difference element $K$ and
 its inverse $K^{-1}$ feasible for symmetric spaces, making good use this
 simplification we obtain the formula:

 \begin{Reference}{Theorem \ref{main}\ 
  (Difference Elements of Hermitean Symmetric Spaces)}\hfill\break
  For every hermitean locally symmetric space the difference element
  $K\,:=\,\exp_p^{-1}\,\circ\,\knc_p$ measuring the deviation between the
  Riemannian and K\"ahlerian normal coordinates reads:
  $$
   KX
   \;\;=\;\;
   \frac{\mathrm{artanh}\,(\,\frac12\,\ad\,IX\,)}{\frac12\,\ad\,IX}\;X
   \qquad\qquad
   K^{-1}X
   \;\;=\;\;
   \frac{\tanh(\,\frac12\,\ad\,IX\,)}{\frac12\,\ad\,IX}\;X
  $$
 \end{Reference}

 \pfill
 Although the arguments and statements in this second part are formulated
 for hermitean locally symmetric spaces, all results except of course
 Corollary \ref{kps} hold true for the locally symmetric among the
 balanced complex affine spaces, for which K\"ahler normal coordinates
 are defined. Unfortunately complex symmetric spaces are never balanced
 unless they are flat, nevertheless there certainly exist balanced complex
 affine manifolds, which are symmetric, but not even pseudo--hermitean
 symmetric spaces.

 In Section \ref{kaehler} we briefly recall several important definitions
 for K\"ahler manifolds with a view towards the parametrization of K\"ahler
 geometry by a single power series in the parameter bundle $\S^{\geq(2,2)}
 T^*M$. In Section \ref{affine} we prove the parallel transport equation
 and discuss the notion of exponentially extended vector fields. Difference
 elements and their relation with other objects relevant for K\"ahler geometry
 are the topic of Section \ref{element}, in which we prove Theorem \ref{univ}
 and describe the Spencer connection for the potential in Lemma \ref{spcon}.
 In the very rigid context of locally symmetric spaces the characteristic
 power series for K\"ahler normal coordinates are calculated in the final
 Section \ref{spaces}. Appendix \ref{examples} details explicit calculations
 for the four families of compact hermitean symmetric undertaken by the
 authors to vindicate the formula of Theorem \ref{main} and verify Corollaries
 \ref{sub} and \ref{symevf}.

 \pfill
 The first author would like to thank the Institute of Mathematics at
 Cuernavaca of the National Autonomous University of Mexico for its
 hospitality during two prolonged stays. The second author is similarly
 indebted to U.~Semmelmann for many fruitful mathematical discussions
 about K\"ahler and quaternionic K\"ahler manifolds as well as for his
 benevolence and generocity in numerous visits to the University of
 Stuttgart. 
\section{K\"ahler Manifolds and K\"ahler Normal Coordinates}
\label{kaehler}
 K\"ahler geometry is a classical topic of Differential Geometry and
 blends the metric structure characteristic for Riemannian geometry with
 the complex structure making complex analytical tools available. Every
 K\"ahler manifold is real analytic, hence the infinite order Taylor series
 of the metric and the complex structure in an arbitrary point determine
 the manifold completely up to coverings. In this section we briefly
 recall the definition of K\"ahler manifolds and then establish a
 parametrization of these two Taylor series by a single power series
 in the parameter vector bundle $\overline\S^{\geq(2,2)}T^*M$ using
 two independent arguments. Excellent introductory texts to K\"ahler
 geometry are the recent textbooks \cite{ba} and \cite{mor}.

 \pfill
 A K\"ahler manifold is a smooth manifold $M$ endowed with a Riemannian
 metric $g$ and an orthogonal almost complex structure $I\,\in\,\Gamma
 (\,\End\,TM\,)$ satisfying the rather strong integrability condition
 that $I$ is parallel with respect to the Levi--Civita connection
 $\nabla$ for $g$:
 \begin{equation}\label{pari}
  \nabla_XI\;\;=\;\;0
 \end{equation}
 Pseudo--K\"ahler manifolds generalize K\"ahler manifolds by weakening
 the positive--definiteness of the Riemannian metric $g$ in the definition
 to non--degeneracy. Every K\"ahler manifold $M$ is actually a complex
 manifold, because the integrability condition implies the vanishing of
 the Nijenhuis tensor and so the Theorem of Newlander--Nierenberg applies.
 A powerful tool to study the topological properties of (pseudo) K\"ahler
 manifolds is the Hodge decomposition
 $$
  \Lam\;T^*M
  \;\;=\;\;
  \bigoplus_{k\,\geq\,0}\Lam^kT^*M\,\otimes_\R\C
  \;\;=\;\;
  \bigoplus_{\kappa,\,\overline\kappa\,\geq\,0}
  \Lam^{\kappa,\,\overline\kappa}T^*M
 $$
 of complex--valued differential forms on $M$ into the eigenspaces of the
 derivation
 \begin{equation}\label{der}
  (\,\Der_I\eta\,)(\,X_1,\,X_2,\,\ldots,\,X_k\,)
  \;\;:=\;\;
  \eta(\,IX_1,\,X_2,\,\ldots,\,X_k\,)
  \;+\;\cdots\;+\;
  \eta(\,X_1,\,X_2,\,\ldots,\,IX_k\,)
 \end{equation}
 extending the complex structure $I$, more precisely
 $$
  \Lam^{\kappa,\,\overline\kappa}T^*
  \;\;:=\;\;
  \{\;\;\eta\,\in\,\Lam^{\kappa+\overline\kappa}T^*\otimes_\R\C\;\;|\;\;\;
  \Der_I\eta\,=\,i\,(\kappa-\overline\kappa)\,\eta\;\;\;\}
 $$
 is the eigenspace of $\Der_I$ for the eigenvalue $i(\kappa-\overline\kappa)$.
 In passing we remark that $\Der_I\,=\,-I\star$ agrees up to sign with the
 representation of the Lie algebra bundle $\End\,TM$ on forms. Replacing
 alternating by symmetric forms on $T$ we obtain the analoguous decomposition
 \begin{equation}\label{sdec}
  \S\,T^*M\otimes_\R\C
  \;\;=\;\;
  \bigoplus_{k\,\geq\,0}\S^kT^*M\otimes_\R\C
  \;\;=\;\;
  \bigoplus_{\kappa,\,\overline\kappa\,\geq\,0}
  \S^{\kappa,\,\overline\kappa}T^*M
 \end{equation}
 into the eigenspaces of $\Der_I$ for the eigenvalues
 $i\,(\kappa-\overline\kappa)\,\in\,\C$ on $\S^{\kappa+\overline\kappa}
 T^*M\otimes_\R\C$. The conjugation of the value of a complex valued
 symmetric multilinear form clearly commutes with $\Der_I$ and thus
 induces isomorphisms $\S^{\kappa,\,\overline\kappa}T^*M\longrightarrow
 \S^{\overline\kappa,\,\kappa}T^*M,\,\eta\longmapsto\overline\eta,$
 for all $\kappa,\,\overline\kappa\,\geq\,0$. In particular a real valued
 symmetric form $\eta\,=\,\overline\eta$ has conjugated bihomogeneous
 components in $\S^{\kappa,\,\overline\kappa}T^*M$ and $\S^{\overline\kappa,
 \,\kappa}T^*M$ respectively. Bihomogeneous real valued symmetric forms
 $\eta\,=\,\overline\eta$ exist only for $\kappa\,=\,\overline\kappa$
 and are characterized by $\Der_I\eta\,=\,0$.
 
 \pfill
 With the complex structure $I$ being parallel on a K\"ahler manifold $M$
 the curvature tensor $R$ of the Levi--Civita connection $\nabla$ associated
 to the Riemannian metric $g$ commutes with $I$
 $$
  R_{X,\,Y}IZ\;\;=\;\;I\,R_{X,\,Y}Z
 $$
 in addition to the standard symmetries of a Riemannian curvature tensor,
 namely
 $$
  g(\,R_{X,\,Y}Z,\,W\,)
  \;\;=\;\;
  -\,g(\,R_{Y,\,X}Z,\,W\,)
  \;\;=\;\;
  -\,g(\,R_{X,\,Y}W,\,Z\,)
  \;\;\stackrel!=\;\;
  +\,g(\,R_{Z,\,W}X,\,Y\,)
 $$
 and of course the first Bianchi identity $R_{X,\,Y}Z\,+\,R_{Y,\,Z}X\,+\,
 R_{Z,\,X}Y\,=\,0$. Combining these standard identities with the
 characteristic commutativity of K\"ahler geometry we obtain
 \begin{equation}\label{sym}
  g(\,R_{X,\,Y}IU,\,IV\,)
  \;\;=\;\;
  g(\,R_{X,\,Y}U,\,V\,)
  \;\;=\;\;
  g(\,R_{IX,\,IY}U,\,V\,)
 \end{equation}
 and conclude that $R_{X,\,IY}\,=\,R_{Y,\,IX}$ is symmetric in $X,\,Y$.
 Via polarization the curvature tensor $R$ of a K\"ahler manifold $M$ is
 thus completely determined by the biquadratic polynomial $T_pM\times T_pM
 \longrightarrow\R,\,(X,Y)\longmapsto-g(\,R_{X,IX}Y,\,IY\,),$ called the
 biholomorphic sectional curvature of $M$ in order to distinguish it from
 the holomorphic sectional curvature:
 
 \begin{Definition}[Holomorphic Sectional Curvature Tensor]
 \hfill\label{hsc}\break
  The holomorphic sectional curvature tensor of a K\"ahler manifold
  $M$ with Riemannian metric $g$ and orthogonal complex structure $I$
  is the section $S\,\in\,\Gamma(\,\S^4T^*M\,)$ defined by:
  $$
   S(\;X,\,Y,\,U,\,V\;)
   \;\;:=\;\;
   8\,\Big(\;g(\;R_{X,\,IY}IU,\,V\;)\;+\;g(\;R_{X,\,IU}IV,\,Y\;)
   \;+\;g(\;R_{X,\,IV}IY,\,U\;)\;\Big)
  $$
  The holomorphic sectional curvature is the associated quartic polynomial
  on $TM$ defined by:
  $$
   S:\quad TM\;\longrightarrow\;\R,\qquad
   X\;\longmapsto\;\frac1{4!}\,S(\;X,\,X,\,X,\,X\;)
   \;\;=\;\;
   g(\;R_{X,\,IX}IX,\,X\;)
  $$
 \end{Definition}

 \noindent
 In passing we observe that the symmetries (\ref{sym}) of the curvature tensor
 $R$ of a K\"ahler manifold $M$ imply that $S$ is actually symmetric in all its
 four arguments, in particular it is completely determined by the associated
 quartic polynomial on $TM$, this is the holomorphic sectional curvature. In
 a similar vein we may deduce $\Der_IS\,=\,0$ via polarization from the
 identity:
 $$
  (\,\Der_IS\,)(\,X,X,X,X\,)
  \;\;=\;\;
  4\,S(\,IX,X,X,X\,)
  \;\;=\;\;
  96\,g(\,R_{IX,\,IX}IX,\,X\,)
  \;\;=\;\;
  0
 $$
 Geometrically $\Der_IS\,=\,0$ is equivalent to the statement that the
 holomorphic sectional curvature $S$ is constant along the fibers of the
 Hopf fibration from the unit sphere $S(\,T_pM\,)$ to the set $\mathbb{P}
 (\,T_pM\,)$ of complex lines in $T_pM$ with respect to the complex structure
 $I_p$.

 \begin{Lemma}[Description of Curvature Tensor]
 \hfill\label{dct}\break
  The curvature tensor $R$ of a K\"ahler manifold $M$ is determined by the
  holomorphic sectional curvature tensor $S\,\in\,\Gamma(\,\S^{2,2}T^*M\,)$
  and can be reconstructed from $S$ by means of:
  $$
   g(\;R_{X,\,Y}U,\,V\;)
   \;\;=\;\;
   \frac1{32}\,\Big(\;S(\,X,IY,IU,V\,)\;-\;S(\,X,IY,U,IV\,)\;\Big)
  $$
 \end{Lemma}

 \noindent
 The proof of this lemma is completely straightforward: Expanding the
 definition of $S$ 
 \begin{eqnarray*}
  +\,S(\,X,IY,IU,V\,)
  &=&
  \hbox to90pt{\hfil$8\,g(\,R_{X,Y}U,\,V\,)$\hfil}
  \;-\;\hbox to90pt{\hfil$8\,g(\,R_{X,U}IV,\,IY\,)$\hfil}
  \;-\;\hbox to90pt{\hfil$8\,g(\,R_{X,IV}Y,\,IU\,)$\hfil}
  \\
  -\,S(\,X,IY,U,IV\,)
  &=&
  \hbox to90pt{\hfil$8\,g(\,R_{X,Y}IU,\,IV\,)$\hfil}
  \;+\;\hbox to90pt{\hfil$8\,g(\,R_{X,IU}V,\,IY\,)$\hfil}
  \;-\;\hbox to90pt{\hfil$8\,g(\,R_{X,V}Y,\,U\,)$\hfil}
 \end{eqnarray*}
 adding and using the first Bianchi identity twice we obtain $4\cdot 8\,
 g(R_{X,Y}U,V)$ as claimed. In particular we can calculate the sectional
 curvature of a K\"ahler manifold $M$ for an arbitrary plane $\mathrm{span}
 \{\,X,Y\,\}\,\subset\,T_pM$ from the holomorphic sectional curvatures:
 \begin{eqnarray*}
  g(\;R_{X,\,Y}Y,\,X\;)
  &=&
  \frac1{32}\,\Big(\;S(\,X,X,IY,IY\,)\;-\;S(\,IX,X,IY,Y\,)\;\Big)
  \\
  &=&
  \frac1{64}\,\Big(\;3\,S(\,X,X,IY,IY\,)\;-\;S(\,X,X,Y,Y\,)\;\Big)
 \end{eqnarray*}
 In the second equality we have replaced $S(\,IX,X,IY,Y\,)$ using the
 argument
 $$
  -4\,S(\,X,X,Y,Y\,)\,+\,4\,S(\,X,X,IY,IY\,)\,+\,8\,S(\,IX,X,IY,Y\,)
  \;\;=\;\;
  0
 $$
 which follows directly from $\Der_IS\,=\,0$ by expanding $(\Der^2_IS)
 (\,X,X,Y,Y\,)\,=\,0$. Generalizing the holomorphic curvature tensor
 to higher orders to incorporate information about the iterated covariant
 derivatives $\nabla R,\,\nabla^2R,\,\ldots$ of the curvature tensor $R$
 as well we arrive at:

 \begin{Definition}[Higher Holomorphic Sectional Curvature Tensors]
 \hfill\label{hhc}\break
  For all $k\,\geq\,4$ the higher holomorphic sectional curvature
  $S_k\,\in\,\Gamma(\,\S^kT^*M\,)$ is defined by:
  $$
   S_k(\;X\;)
   \;\;\widehat=\;\; 
   \frac1{k!}\,S_k(\;X,\,\ldots,\,X\;)
   \;\;:=\;\;
   \frac1{(k-4)!}\,
   g(\;(\,\nabla^{k-4}_{X,\,\ldots,\,X}R\,)^{}_{X,\,IX}IX,\;X\;)
  $$
 \end{Definition}

 \noindent
 Calculating $\Der_IS_k$ we observe that the last four summands expected
 from the definition of $\Der_I$ all vanish by the symmetries (\ref{sym})
 of curvature tensors of K\"ahler type. In consequence the bihomogeneous
 components $S_{\kappa,\,\overline\kappa}\,\in\,\Gamma(\,\S^{\kappa,\,
 \overline\kappa}T^*M\otimes_\R\C\,)$ of the higher order holomorphic
 sectional curvature tensor $S_k$ vanish unless $-k+4\,\leq\,\kappa-
 \overline\kappa\,\leq\,k-4$, this is to say:
 \begin{equation}\label{stot}
  S_{\mathrm{total}}
  \;\;:=\;\;
  \bigoplus_{k\,\geq\,4}S_k
  \;\;=\;\;
  \bigoplus_{\kappa,\,\overline\kappa\,\geq\,2}S_{\kappa,\,\overline\kappa}
  \;\;\in\;\;
  \Gamma(\;\overline\S^{\geq(2,2)}T^*M\;)
 \end{equation}
 Of course the bihomogeneous components $S_{\kappa,\,\overline\kappa}$
 and $S_{\overline\kappa,\,\kappa}$ are conjugated, because the higher
 holomorphic sectional curvature tensor $S_k\,\in\,\Gamma(\,\S^kT^*M\,)$
 is a real valued polynomial by definition. In light of Lemma \ref{dct}
 we may replace the curvature tensor $R$ with $S\,=\,S_4$, similarly we
 may identify $S_5$ with the covariant derivative $\nabla R$ in the
 following way:
 
 \begin{Corollary}[Covariant Derivative of Curvature]
 \hfill\label{cdr}\break
  The first higher holomorphic sectional curvature $S_5\,\in\,
  \Gamma(\,\S^5T^*M\,)$ of a K\"ahler manifold $M$ describes the covariant
  derivative $\nabla R$ of the curvature tensor by means of the formula: 
  \begin{eqnarray*}
   g(\;(\nabla_XR)_{Y,\,Z}U,\,V\;)
   &=&
   \frac1{192}\,\Big(\,
    +\;S_5(\,X,\,Y,\,IZ,\,IU,\,V\,)\;-\;S_5(\,X,\,Y,\,IZ,\,U,\,IV\,)
   \\[-2pt]
   &&
   \qquad\quad
   +\;S_5(\,X,\,IY,\,Z,\,U,\,IV\,)\;-\;S_5(\,X,\,IY,\,Z,\,IU,\,V\,)\;\Big)
  \end{eqnarray*}
 \end{Corollary}

 \noindent
 This formula for $\nabla R$ is a direct corollary of the description
 of $R$ in Lemma \ref{dct}, because
 $$
  \frac1{4!}\,S_5(\,Y,\,X,\,X,\,X,\,X\,)
  \;\;=\;\;
  g(\,(\nabla_YR)_{X,\,IX}IX,\,X\,)
  \;+\;4\,g(\,(\nabla_XR)_{Y,\,IX}IX,\,X\,)
 $$
 implies via the second Bianchi identity:
 \begin{eqnarray*}
  \lefteqn{\frac1{4!}\,S_5(\,Y,\,X,\,X,\,X,\,X\,)
   \;+\;\frac1{4!}\,S_5(\,Y,\,IX,\,IX,\,IX,\,IX\,)}
  \qquad
  &&
  \\[3pt]
  &=&
  2\,g(\,(\nabla_YR)_{X,\,IX}IX,\,X\,)
  \;+\;4\,g(\,(\nabla_XR)_{Y,\,IX}IX,\,X\,)
  \;+\;4\,g(\,(\nabla_{IX}R)_{X,\,Y}IX,\,X\,)
  \\
  &=&
  6\,g(\,(\nabla_YR)_{X,\,IX}IX,\,X\,)
  \;\;=\;\;
  \frac6{4!}\,(\,\nabla_YS_4\,)(\,X,\,X,\,X,\,X\,)
 \end{eqnarray*}
 In consequence of Lemma \ref{dct} and Corollary \ref{cdr} every K\"ahler
 geometry can be constructed up to covering from the total holomorphic
 sectional curvature $S_{\mathrm{total}}$ defined in (\ref{stot}),
 in particular all iterated covariant derivatives $\nabla R,\,\nabla^2R,\,
 \ldots$ and thus the infinite order Taylor series of both $g$ and $I$ in
 exponential coordinates are determined by $S_{\mathrm{total}}$. For the
 time being however this reconstruction is a theoretical possibility based
 on counting the free parameters in the sequence $R,\,\nabla R,\,\ldots$
 using the formal theory of partial differential equations \cite{bc3}.
 En nuce the problem addressed in this article is that we are lacking
 the explicit formulas needed for this reconstruction to work in practise.
 Moreover $S_{\mathrm{total}}$ is not the only power series with the correct
 number of parameters, another interesting candidate parametrizing K\"ahler
 geometry by a section of the parameter bundle $\overline\S^{\geq(2,2)}T^*M$
 arises from a suitable potential:

 \begin{Definition}[Local Potential Functions]
 \hfill\label{pf}\break
  A local potential function for a K\"ahler manifold $M$ with Riemannian
  metric $g$ and complex structure $I$ is a smooth function $\theta^\loc:
  \,U\longrightarrow\R$ defined on an open subset $U\,\subset\,M$, which
  is a preimage of the K\"ahler form $\omega(X,Y)\,:=\,g(IX,Y)$ in the
  sense of the $\partial\overline\partial$--Lemma
  $$
   \omega
   \;\;=\;\;
   \frac i2\,\partial\overline\partial\,\theta^\loc
  $$
  where $\partial$ and $\overline\partial$ are the $(1,0)$ and
  $(0,1)$--components of the exterior derivative $d\,:=\,\partial\,+\,
  \overline\partial$. Note that $\omega\,\in\,\Gamma(\,\Lam^2T^*M\,)$ is
  parallel under the Levi--Civita connection and so is closed.
 \end{Definition}

 \noindent
 A local potential function is never unique as it can always be modified
 $\theta^\loc\,\rightsquigarrow\,\theta^\loc+\Re\,f$ by adding the real
 part of a holomorphic function $f:\,U\longrightarrow\C$. Instead of the
 K\"ahler form $\omega$ we may employ the hermitean form $h\,:=\,g\,+\,i\,
 \omega$ associated to a K\"ahler geometry on a manifold $M$ in the equation
 characterizing a local potential function $\theta^\loc$
 \begin{equation}\label{cgih}
  h(\;X,\;Y\;)
  \;\;:=\;\;
  g(\;X,\;Y\;)\;+\;i\,\omega(\;X,\;Y\;)
  \;\;=\;\;
  (\,\Hess\,\theta^\loc\,)(\;\pr^{0,1}X,\;\pr^{1,0}Y\;)
 \end{equation}
 where $\pr^{1,0}\,:=\,\frac12(\id-iI)$ and $\pr^{0,1}\,:=\,\frac12(\id+iI)$
 are the $I$--eigenprojections, or alternatively:
 \begin{equation}\label{cgi}
  g(\;X,\,Y\;)
  \;\;=\;\;
  \frac14\,\Big(\;(\,\Hess\,\theta^\loc\,)(\;X,\,Y\;)
  \;+\;(\,\Hess\;\theta^\loc\,)(\;IX,\,IY\;)\;\Big)
 \end{equation}
 In both equations the Hessian $D\,\circ\,d$ can actually be taken with
 respect to an arbitrary torsion free connection $D$ on the tangent bundle
 making $I$ parallel, it is not necessary to use the Levi--Civita connection
 $\nabla$ for this purpose. In fact the difference $A_XY\,:=\,D_XY\,-\,
 \hat D_XY$ of two such connections is symmetric $A_XY\,=\,A_YX$ with
 $A_XIY\,=\,IA_XY$ and thus automatically satisfies $A_XY\,+\,A_{IX}IY\,=\,
 0\,=\,A_{\pr^{0,1}X}\pr^{1,0}Y$. The liberty granted by this observation
 is very useful in verifying equations (\ref{cgih}) and (\ref{cgi}) directly
 in local holomorphic coordinates $(\,z^1,\,\ldots,\,z^n\,)$, because we may
 simply choose $D$ to be the trivial connection arising from the local
 trivialization of the complexified tangent bundle by the Wirtinger vector
 fields
 $$
  \frac\partial{\partial z^\mu}
  \;\;:=\;\;
  \frac12\,\Big(\,\frac\partial{\partial x^\mu}\,-\,i
  \frac\partial{\partial y^\mu}\;\Big)
  \qquad\qquad
  \frac\partial{\partial\overline z^\mu}
  \;\;:=\;\;
  \frac12\,\Big(\,\frac\partial{\partial x^\mu}\,+\,i
  \frac\partial{\partial y^\mu}\,\Big)
 $$
 where $(\,x^1,\,y^1,\,\ldots,\,x^n,\,y^n\,)$ are the real and imaginary
 parts of $(\,z^1,\,\ldots,\,z^n\,)$. The Riemannian metric $g$ is completely
 determined by its mixed components $g_{\mu\,\overline\nu}\,:=\,g(\frac\partial
 {\partial z^\mu},\frac\partial{\partial\overline z^\nu})$ in holomorphic
 coordinates, because the eigensubbundles $T^{1,0}M$ and $T^{0,1}M$ of
 the complexified tangent bundle $TM\otimes_\R\C$ under the skew symmetric
 endomorphism $I$ are isotropic, so that all pure components $g_{\mu\,\nu}
 \,=\,0\,=\,g_{\overline\mu\,\overline\nu}$ vanish. In holomorphic
 coordinates we thus find the following local expansions of the metric,
 the K\"ahler and the hermitean form:
 $$
  g
  \;\;=\;\;
  \sum_{\mu\,\nu}g_{\mu\,\overline\nu}\,dz^\mu\,\cdot\,d\overline z^\nu
  \qquad
  \omega
  \;\;=\;\;
  i\,\sum_{\mu\,\nu}g_{\mu\,\overline\nu}\,dz^\mu\,\wedge\,d\overline z^\nu
  \qquad
  h
  \;\;=\;\;
  2\,\sum_{\mu\,\nu}g_{\mu\,\overline\nu}\,d\overline z^\nu\,\otimes\,dz^\mu
 $$
 Comparing the expansion for $\omega$ with the analoguous local expansion
 of $\frac i2\,\partial\overline\partial\,\theta^\loc$ we find
 $$
  \partial\,\overline\partial\,\theta^\loc
  \;\;=\;\;
  \sum_{\mu\,\nu}
  \frac{\partial^2\;\theta^\loc}{\partial z^\mu\partial\overline z^\nu}
  \,dz^\mu\,\wedge\,d\overline z^\nu
  \qquad\Longrightarrow\qquad
  g_{\mu\,\overline\nu}
  \;\;=\;\;
  \frac12\,\frac{\partial^2}{\partial z^\mu\partial\overline z^\nu}
  \,\theta^\loc
 $$
 and reinserting this relation between $\theta^\loc$ and the mixed components
 of $g$ into the expansions of the hermitean form $h$ we verify equation
 (\ref{cgih}) and in turn equation (\ref{cgi}) for the Hessian $D\,\circ\,d$
 associated to the local trivial connection $D\,\frac\partial {\partial z^\mu}
 \,=\,0\,=\,D\frac\partial{\partial\overline z^\nu}$.

 The concept of local potentials for K\"ahler geometries allows us to single
 out special holomorphic coordinates charts $(\,z^1,\,\ldots,\,z^n\,)$ on a
 K\"ahler manifold by requiring that the local potential, after adding
 the  real part of a holomorphic function, takes as simple a form as possible.
 Uniqueness of these holomorphic coordinates is usually stipulated modulo
 unitary transformations only, in order to remove this ambiguity we borrow
 from the Riemannian exponential map the idea of anchored coordinates centered
 in a point $p\,\in\,M$. An anchored local coordinate chart is a smooth map
 $\varphi:\,T_pM\longrightarrow M$, which is a diffeomorphism of some
 neighborhood of $0\,\in\,T_pM$ to a neighborhood of $\varphi(0)\,=\,p$
 such that the differential
 $$
  \varphi_{*,\,0}:\quad
  T_pM\;\;\cong\;\;T_0(\,T_pM\,)\;\;\longrightarrow\;T_pM,\qquad
  Z\;\longmapsto\;Z
 $$
 equals $\id_{T_pM}$ under the natural identification $T_0(T_pM)\,\cong\,T_pM$.
 Actually we do not insist on $T_pM$ to be the domain of $\varphi$, hence
 it may not be defined outside a neighborhood of $0$.

 \begin{Theorem}[K\"ahler Normal Coordinates and Normal Potentials \cite{bo}]
 \hfill\label{knc}\break
  In every point $p\,\in\,M$ of a K\"ahler manifold $M$ there exist unique
  anchored local holomorphic coordinates $\knc_p:\,T_pM\longrightarrow M$
  centered in $p$ and a unique local potential function $\theta^\loc_p$
  such that the infinite order Taylor series $\theta_p$ of the pull back
  of $\theta^\loc_p$ along $\knc_p$ is congruent
  $$
   \theta^\loc_p(\,\knc_pX\,)
   \;\;\raise-6pt\hbox{$\stackrel\sim{\scriptstyle X\,\to\,0}$}\;\;
   \theta_p(\;X\;)
   \;\;\equiv\;\;
   g_p(\;X,\,X\;)
   \qquad\textrm{mod}\qquad
   \overline\S^{\geq(2,2)}T^*_pM
  $$
  to the norm square function on $T_pM$ modulo a power series, which is at
  least quadratic in both holomorphic and antiholomorphic coordinates. Due
  to their uniqueness the holomorphic coordinates $\knc_p$ are called K\"ahler
  normal coordinates with normal potential $\theta^\loc_p\,\widehat=\,
  \theta_p$.
 \end{Theorem}

 \noindent
 The double uniqueness statement in this theorem implies in particular
 that for given anchored local holomorphic coordinates $\varphi:\,T_pM
 \longrightarrow M$ the existence of a power series potential $\theta_p
 \,\in\,\overline\S\,T^*_pM$ satisfying the normalization congruence is
 sufficient to conclude that $\theta_p$ is the Taylor series of the normal
 potential $\theta_p^\loc$ and that $\varphi\,=\,\knc_p$ are the K\"ahler
 normal coordinates, this argument will be used in explicit examples in
 Appendix \ref{examples}. The original account \cite{bo} of Bochner
 contains a very readable, elementary proof of Theorem \ref{knc}, for
 this reason we will not discuss Theorem \ref{knc} independently of
 Theorem \ref{univ} below.
\section{Affine Exponential Coordinates}
\label{affine}
 In the study of the geometry of Riemannian or more generally affine
 manifolds the exponential map provides an indispensable tool in both
 explicit calculations and theoretic considerations. Classically Jacobi
 vector fields are used to describe the differential of the exponential
 map, alternatively the forward and backward parallel transport defined
 in this section can be used for this purpose. In particular we will use
 Jacobi vector fields to give an apriori proof of the parallel transport
 equation, which describe the Taylor series of parallel tensors, without
 calculating its solution first. Moreover we will introduce the concept
 of exponentially extended vector fields for general affine manifolds,
 which generalize the left invariant vector fields on Lie groups and
 the transvection Killing vector fields on symmetric spaces.
 
 \pfill
 Affine linear spaces are characterized by the presence of a distinguished
 family of curves, the family of straight lines. The choice of a connection
 $\nabla$ on the tangent bundle $TM$ of a manifold $M$ replaces this
 distinguished family by the family of geodesics, the solutions $\gamma:\,
 \R\longrightarrow M$ to the geodesic equation associated to the choice of
 connection. According to the Theorem of Picard--Lindel\"of a unique solution
 $\gamma$ to the geodesic equation $\frac\nabla{dt}\dot\gamma\,=\,0$ exists
 for all initial values so that we may define the exponential map centered
 in a point
 $$
  \exp_p:\;\;T_pM\;\longrightarrow\;M,\qquad X\;\longmapsto\;\gamma^X(\,1\,)
 $$
 by sending a tangent vector $X$ to the value $\gamma^X(1)$ of the
 unique geodesic $\gamma^X$ with initial values $\gamma^X(0)\,=\,p$ and
 $\dot\gamma^X(0)\,=\,X$. Constant reparametrizations of geodesics are
 geodesics $\gamma^{\lambda X}(t)\,=\,\gamma^X(\lambda t)$, hence the image
 of the ray $t\longmapsto tX$ under exponential coordinates is the geodesic
 $t\longmapsto\gamma^{tX}(1)\,=\,\gamma^X(t)$ we started with. The parallel
 transport along this ray geodesic $\PT^\nabla(X)\,:=\,\PT^\nabla_{\gamma^X}
 (1)$ is by construction a linear isomorphism of tangent spaces:
 $$
  \PT^\nabla(\,X\,):\quad
  T_pM\;\stackrel\cong\longrightarrow\;T_{\exp_pX}M
 $$
 This isomorphism allows us to describe the differential of exponential
 coordinates $\exp_p$ centered in $p\,\in\,M$, either by using the forward
 parallel transport defined as the composition
 $$
  \Phi(\,X\,):\qquad
  T_pM\;\stackrel{\PT^\nabla(X)}\longrightarrow\;T_{\exp_pX}M
  \;\stackrel{(\exp_p)^{-1}_{*,X}}\longrightarrow\;T_XT_pM
  \;\stackrel\cong\longrightarrow\;T_pM
 $$
 wherever the exponential map $\exp_p$ is a local diffeomorphism, or the
 backward transport
 $$
  \Phi^{-1}(\,X\,):\qquad
  T_pM\;\stackrel\cong\longrightarrow\;T_XTpM
  \;\stackrel{(\exp_p)_{*,X}}\longrightarrow\;T_{\exp_pX}M
  \;\stackrel{\PT^\nabla(X)^{-1}}\longrightarrow\;T_pM
 $$
 defined for all $X\,\in\,T_pM$ under the rather mild assumption that
 $M$ is a complete affine manifold. The infinite order Taylor series of
 $\Phi^{-1}$ in the origin $0\,\in\,T_pM$ is an explicitly known power
 series in the curvature tensor $R$ of $M$ and its iterated covariant
 derivatives $\nabla R,\,\nabla^2R,\,\ldots$ evaluated at $p$. More
 precisely the Taylor series of $\Phi^{-1}$ in the origin equals the unique
 power series solution $\Phi^{-1}\,\in\,\overline\S\;T^*_pM\otimes\End\,T_pM$
 to a formal differential equation with initial value $\Phi(0)\,=\,\id$
 involving the number operator $N$ on power series:

 \begin{Lemma}[Parallel Transport Equation]
 \hfill\label{pte}\break
  The infinite order Taylor series of the backward parallel transport
  $\Phi^{-1}:\,T_pM\longrightarrow\End\,T_pM$ in the center $0\,\in\,T_pM$
  of exponential coordinates of an affine manifold $M$ is characterized
  as a formal power series $\Phi^{-1}\,\in\,\overline\S\,T^*_pM\otimes
  \End\,T_pM$ by the formal differential equation
  $$
   N\,(\,N\,+\,1\,)\,\Phi^{-1}
   \;\;=\;\;
   \mathscr{R}\,\Phi^{-1}\;+\;N(\,\mathscr{T}\,\Phi^{-1}\,)
  $$
  in which $N$ denotes the number operator on power series and the power
  series $\mathscr{R}$ and $\mathscr{T}$ reflect the infinite  order Taylor
  series of the curvature and the torsion respectively:
  $$
   \mathscr{R}(\,X\,)\,Y
   \;\;:=\;\;
   \sum_{k\geq0}\frac1{k!}\,(\,\nabla^k_{X,\ldots,X}R\,)^{}_{X,\,Y}X
   \qquad
   \mathscr{T}(\,X\,)\,Y
   \;\;:=\;\;
   \sum_{k\geq0}\frac1{k!}\,(\,\nabla^k_{X,\ldots,X}T\,)(\,X,\,Y\,)
  $$
 \end{Lemma}

 \proof
 In order to prove the lemma we want to study the standard Jacobi
 equation for a vector field $J\,\in\,\Gamma(\gamma^*TM)$ along a
 geodesic $\gamma$ in a manifold $M$ with respect to an affine, not
 necessarily torsion free connection $\nabla$ on the tangent bundle
 with curvature $R$ and torsion $T$:
 $$
  \mathrm{Jac}_\gamma J
  \;\;:=\;\;
  \frac{\nabla^2}{dt^2}\,J
  \;+\;\frac\nabla{dt}\,T(\,J,\,\dot\gamma\,)
  \;+\;R_{J,\,\dot\gamma}\dot\gamma
  \;\;\stackrel?=\;\;
  0
 $$
 Multiplying by $t^2$ we obtain for every solution $J$ to this equation
 the equality:
 \begin{equation}\label{jacx}
  t\,\frac\nabla{dt}\,\Big(\,t\,\frac\nabla{dt}\,-\,1\,\Big)\,J
  \;\;=\;\;
  R_{t\,\dot\gamma,\,J}(\,t\,\dot\gamma\,)\;+\;
  \Big(\,t\,\frac\nabla{dt}\,-\,1\,\Big)\,T(\;t\,\dot\gamma,\;J\;)
 \end{equation}
 In particular we are interested in the family of solutions to the Jacobi
 equation, which arise naturally in exponential coordinates by varying the
 geodesic rays. More precisely we consider for a point $p\,\in\,M$ and an
 endomorphism $A\,\in\,\End\,T_pM$ of its tangent space the geodesic variation
 $(s,t)\longmapsto\exp_p(\,t\,e^{sA}X\,)$ for some tangent vector $X\,\in\,
 T_pM$. This geodesic variation induces a Jacobi field $J_A(t)$ along the
 geodesic ray $\gamma:\,t\longmapsto\exp_p(tX)$ determined by $X$:
 $$
  J_A(\;t\;)
  \;\;:=\;\;
  \left.\frac\partial{\partial s}\right|_0\exp_p(\;t\,e^{sA}\,X\;)
  \;\;=\;\;
  (\,\exp_p\,)_{*,\,tX}(\,t\,AX\,)
  \;\;=\;\;
  t\,\PT^\nabla(\,tX\,)\;\Phi^{-1}(\,tX\,)\,AX
 $$
 Multiplying the Jacobi equation (\ref{jacx}) for this Jacobi field by
 $\PT^\nabla(tX)^{-1}$ and using the cha\-racteristic property of the
 parallel transport $\PT^\nabla(tX)^{-1}\circ\frac\nabla{dt}\,=\,\frac d{dt}
 \circ\PT^\nabla(tX)^{-1}$ we obtain
 \begin{eqnarray}
  t\frac d{dt}\,(\,t\frac d{dt}-1\,)
  \Big(\,t\,\Phi^{-1}(tX)\,AX\,\Big)
  &=&
  \Big(\,\PT^\nabla(tX)^{-1}R\,\Big)_{tX,\,t\,\Phi^{-1}(tX)\,AX}(\,tX\,)
  \label{jacz}
  \\
  &&
  \;+\;\Big(\,t\,\frac d{dt}\,-\,1\,\Big)
  \;\Big(\,\PT^\nabla(tX)^{-1}T\,\Big)(\,tX,\,t\,\Phi^{-1}(tX)\,AX\,)
  \nonumber
 \end{eqnarray}
 where $\PT^\nabla(tX)^{-1}R$ for example denotes the parallel
 transport for the curvature tensor
 $$
  \Big(\;\PT^\nabla(\,tX\,)^{-1}R\;\Big)_{U,\,V}W
  \;\;:=\;\;
  \PT^\nabla(\,tX\,)^{-1}\;\Big(\;R_{\PT^\nabla(\,tX\,)\,U,
  \,\PT^\nabla(\,tX\,)\,V}\PT^\nabla(\,tX\,)\,W\;\Big)
 $$
 the analoguous definition of the parallel transport
 $\PT^\nabla(\,tX\,)^{-1}T$ for the torsion tensor is omitted. In general
 the derivative of a geodesic like $t\longmapsto\exp_p(tX)$ is given by
 the parallel transport of the initial tangent vector $\dot\gamma(t)\,=\,
 (\exp_p)_{*,\,tX}X\,=\,\PT^\nabla(\,tX\,)\,X$, incidentally this argument
 directly implies the so called Gau\ss\ Lemma valid for all $X\,\in\,T_pM$
 and $t\,\in\,\R$:
 \begin{equation}\label{gss}
  \Phi^{-1}(\,tX\,)\,X\;\;=\;\;X\;\;=\;\;\Phi(\,tX\,)\,X
 \end{equation}
 Of course the point in defining the parallel transport for the
 curvature and torsion tensor is that $\PT^\nabla(\,tX\,)^{-1}R$ and
 $\PT^\nabla(\,tX\,)^{-1}T$ are trilinear and bilinear maps on $T_pM$
 respectively with values in $T_pM$ independent of the argument $X\,\in\,T_pM$.
 Hence it makes sense to talk about their infinite order Taylor series as
 $X$ approaches $0\,\in\,T_pM$, which are given by
 \begin{equation}\label{aex}
  \PT^\nabla(\,X\,)^{-1}R
  \;\;\raise-6pt\hbox{$\stackrel\sim{\scriptstyle X\,\to\,0}$}\;\;
  \sum_{k\,\geq\,0}\frac1{k!}\nabla^k_{X,\,\ldots,\,X}R
  \qquad\qquad
  \PT^\nabla(\,X\,)^{-1}T
  \;\;\raise-6pt\hbox{$\stackrel\sim{\scriptstyle X\,\to\,0}$}\;\;
  \sum_{k\,\geq\,0}\frac1{k!}\nabla^k_{X,\,\ldots,\,X}T
 \end{equation}
 due to the definition of iterated covariant derivatives, a more detailed
 derivation of these asymptotic expansions can be found for example in
 \cite{w2}. Replacing $\PT^\nabla(\,tX\,)^{-1}R$ and $\PT^\nabla(\,tX\,)^{-1}T$
 as well as the backward parallel transport $\Phi^{-1}(\,tX\,)$ by their
 respective infinite order Taylor series in equation (\ref{jacz}) we obtain
 the following identity of power series
 \begin{eqnarray*}
  \lefteqn{t\frac d{dt}\,
  (\,t\frac d{dt}-1\,)\,\Big(\;\Phi^{-1}(\,tX\,)\,A(\,tX\,)\;\Big)}\qquad
  &&
  \\
  &=&
  \mathscr{R}(\,tX\,)\,\Phi^{-1}(\,tX\,)\,A(\,tX\,)\;+\;(\,t\frac d{dt}-1\,)
  \,\Big(\;\mathscr{T}(\,tX\,)\,\Phi^{-1}(\,tX\,)\,A(\,tX\,)\;\Big)
 \end{eqnarray*}
 in which every occurrence of the argument $X\,\in\,T_pM$ comes along with
 a factor $t$ and vice versa. In turn we may replace the differential operator
 $t\frac d{dt}$ by the number operator $N$ on power series in $X$ and evaluate
 at $t\,=\,1$ to reduce this power series identity to:
 $$
  N\,(\,N\,-\,1\,)\,(\,\Phi^{-1}(X)\,AX\,)
  \;\;=\;\;
  \mathscr{R}(\,X\,)\,\Phi^{-1}(\,X\,)\,AX
  \;+\;(\,N\,-\,1\,)\,(\,\mathscr{T}(\,X\,)\,\Phi^{-1}(\,X\,)\,AX\,)
 $$
 Recall now that the endomorphism $A\,\in\,\End\,T_pM$ of the tangent space
 $T_pM$ can be chosen arbitrarily in this identity. A fortiori the infinite
 order Taylor series $\Phi^{-1}$ of the backward parallel transport satisfies
 the differential equation $(N+1)N\Phi^{-1}\,=\,\mathscr{R}\,\Phi^{-1}\,+\,
 N(\,\mathscr{T}\,\Phi^{-1}\,)$ in light of the observation that the following
 homogeneous linear map of degree $+1$
 $$
  \iota:\;\;\S^\bullet T^*\otimes\End\,T\;\longrightarrow\;
  \Hom(\,\End\,T,\,\S^{\bullet+1}T^*\otimes T\,),\qquad
  \Psi\;\longmapsto\;\Big(\,A\,\longmapsto\,\Psi(\,\cdot\,)
  \,A(\,\cdot\,)\,\Big)
 $$
 is injective for every finite dimensional vector space $T$; the degree
 $+1$ homogeneity of $\iota$ evidently accounts for the shift $N\,
 \rightsquigarrow\,N+1$ in the differential equation for $\Phi^{-1}$.
 To justify our observation about the injectivity of $\iota$ we specify
 a linear map $\iota^*$ in the opposite direction
 $$
  \iota^*:\;\;\Hom(\,\End\,T,\,\S^{\bullet+1}T^*\otimes T\,)\;
  \longrightarrow\;\S^\bullet T^*\otimes\End\,T,\qquad F\;\longmapsto\;\iota^*F
 $$
 by summing over a pair $\{\,E_\mu\,\}$ and $\{\,dE_\mu\,\}$ of dual bases
 for $T$ and $T^*$ respectively:
 $$
  [\;\iota^*F\;](\,X\,)\,Y
  \;\;:=\;\;
  \sum_\mu(\,E_\mu\,\ins\,F(\,dE_\mu\otimes Y\,)\,)(\,X\,)
 $$
 The injectivity of $\iota$ is then a direct consequence of the following
 identity for $\iota^*\iota$:
 $$
  [\,\iota^*\iota\Psi\,](\,X\,)\,Y
  \;\;=\;\;
  \sum_\mu
  \Big(E_\mu\ins(\,dE_\mu(\,\cdot\,)\,\Psi(\,\cdot\,)\,Y\,)\Big)(X)
  \;\;=\;\;
  [\,(\,N+\dim\,T\,)\,\Psi\,](\,X\,)\,Y
 $$
 \vskip-25pt\qed

 \pfill
 Perhaps it is a good idea to remind the reader that the parallel transport
 equation formulated in Lemma \ref{pte} is only a formal differential equation
 for the infinite order Taylor series $\Phi^{-1}$ of the backward parallel
 transport. Its formal solution is easily expanded to arbitrary order in $X$,
 but it is much easier to recall the formal differential equation than its
 explicit solution. Expanding the solution to order $4$ in $X$ for example
 we obtain in the torsion free case:
 \begin{eqnarray*}
  \Phi^{-1}(\,X\,)\,Y
  &=&
  Y\;+\;\frac16\,R_{X,\,Y}X\;+\;\frac1{12}\,(\nabla_XR)_{X,\,Y}X
  \\
  &&
  \qquad+\;\frac1{40}\,(\nabla^2_{X,\,X}R)^{}_{X,\,Y}X
  \;+\;\frac1{120}\,R_{X,\,R_{X,\,Y}X}X\;+\;O(\,X^5\,)
 \end{eqnarray*}
 Similarly the parallel transport equation can be solved explicitly for Lie
 groups. Consider for example the flat connection $\nabla^L$ on the tangent
 bundle $TG$ of a Lie group $G$, which makes all left invariant vector fields
 parallel. By construction this connection is flat $R^L\,=\,0$ and its
 torsion $T^L(X,Y)\,=\,-[\,X,\,Y\,]_{\mathrm{algebraic}}$ is parallel
 with respect to $\nabla^L$ so that $\mathscr{T}^L(X)\,=\,-\ad\,X$. Under
 the ansatz $\Phi^{-1}(\,X\,)\,=\,\varphi^{-1}(\,\ad\,X\,)$ the parallel
 transport equation becomes
 $$
  (\,x\frac d{dx}\,+\,1\,)\,x\frac d{dx}\,\varphi^{-1}(\,x\,)
  \;\;=\;\;
  x\frac d{dx}\,(\,-\,x\,\varphi^{-1}(\,x\,)\,)
  \qquad\qquad
  \varphi^{-1}(\,0\,)\;\;=\;\;1
 $$
 with unique solution $\varphi^{-1}(x)\,=\,\frac{e^{-x}-1}{-x}$, in turn
 we find the well--known explicit solution:
 $$
  \Phi^{-1}(\;X\;)\;\;=\;\;\frac{e^{-\ad\,X}\,-\,\id}{-\ad\,X}
  \qquad\qquad
  \Phi(\;X\;)\;\;=\;\;\frac{-\ad\,X}{e^{-\ad\,X}\,-\,\id}
 $$
 The parallel transport equation for symmetric spaces can be solved using
 a similar ansatz.

 \begin{Corollary}[Taylor Series of Parallel Tensors]
 \hfill\label{taypar}\break
  The unique power series solution $\Phi^{-1}\,\in\,\overline{\S}\,T^*_pM
  \,\otimes\,\End\,T_pM$ to the parallel transport equation $N\,(N+1)\,
  \Phi^{-1}\,=\,\mathscr{R}\,\Phi^{-1}\,+\,N\,(\,\mathscr{T}\,\Phi^{-1}\,)$
  determines the Taylor series of every tensor parallel with respect to
  $\nabla$ in exponential coordinates. For a purely covariant parallel
  tensor $\eta$ say the resulting Taylor series reads for $X\,\in\,T_pM$
  and $A_1,\,\ldots,\,A_r\,\in\,T_X(T_pM)\,\cong\,T_pM$:
  $$
   (\,\mathrm{taylor}_0(\,\exp_p^*\eta\,)\,)_X(\,A_1,\,\ldots,\,A_r\,)
   \;\;=\;\;
   \eta_p(\;\Phi^{-1}(\,X\,)\,A_1,\;\ldots,\;\Phi^{-1}(\,X\,)\,A_r\;)
  $$
 \end{Corollary}

 \noindent
 Of course the restriction of a tensor $\eta$ parallel with respect
 to $\nabla$ on all of $M$ is parallel along every radial geodesic
 $t\longmapsto\exp_p(\,tX\,)$, with this in mind we find the relation
 \begin{eqnarray*}
  (\,\exp^*_p\eta\,)_X(\,A_1,\,\ldots,\,A_r\,)
  &=&
  \eta_{\exp_pX}(\,(\,\exp_p\,)_{*,\,X}A_1,\,\ldots,\,
  (\,\exp_p\,)_{*,\,X}A_r\,)
  \\
  &=&
  \eta_p(\,\PT^\nabla(X)^{-1}(\,\exp_p\,)_{*,\,X}A_1,\,\ldots,\,
  \PT^\nabla(X)^{-1}(\,\exp_p\,)_{*,\,X}A_r\,)
 \end{eqnarray*}
 and so the definition of $\Phi^{-1}$ implies the Corollary. For mixed
 co-- and contravariant tensors the formulas become slightly more
 complicated, but the argument in itself remains valid.

 \pfill
 A notion closely related to the construction of the forward and backward
 parallel transport is the notion of exponentially extended vector fields,
 which in a sense describe the covariant derivative of the exponential map.
 The exponentially extended vector field associated to a tangent vector
 $Z\,\in\,T_pM$ is the smooth vector field defined on the domain of the
 exponential map $\exp_p$ in $T_pM$ as the derivative of the following
 family of diffeomorphisms of $T_pM$
 $$
  Z^{\exp}\;\;:=\;\;\left.\frac d{dt}\right|_0
  \Big(\;T_pM\;\stackrel{\PT^\nabla_\gamma(t)}\longrightarrow\;
  T_{\gamma(t)}M\;\stackrel{\exp_{\gamma(t)}}\longrightarrow\;
  M\;\stackrel{(\exp_p)^{-1}}\longrightarrow\;T_pM\;\Big)
  \;\;\in\;\;\Gamma(\;T\,(T_pM)\;)
 $$
 where $\PT^\nabla_\gamma(\,t\,)$ denotes the parallel transport along a
 curve $\gamma$ representing $Z\,=\,\left.\frac d{dt}\right|_0\gamma$.
 The exponentially extended vector field $Z^{\exp}$ does not depend on
 the curve $\gamma$ used in its definition, moreover its value $Z^{\exp}
 (\,0\,)\,=\,Z$ in the origin equals the tangent vector we started with
 under the identification $T_0(T_pM)\,\cong\,T_pM$. Somewhat more general the
 infinite order Taylor series of $Z^{\exp}$ in $0\,\in\,T_pM$ is determined
 by a power series $\Theta$ on $T_pM$ with values in $\End\,T_pM$ via
 $$
  Z^{\exp}(\;X\;)
  \;\;\raise-6pt\hbox{$\stackrel\sim{\scriptstyle X\,\to\,0}$}\;\;
  \Theta(\,X\,)\,Z
  \qquad\qquad
  \Theta\;\;\in\;\;\overline{\S}\,T^*_pM\,\otimes\,\End\,T_pM
 $$
 which has a universal expansion in terms of the curvature and the torsion
 of the connection $\nabla$ together with all their covariant derivatives.
 In the torsion free case $T\,=\,0$ for example
 \begin{eqnarray*}
  \lefteqn{Z^{\mathrm{exp}}(\,X\,)}
  &&
  \\
  &=&
  Z\;+\;\frac13\,R_{X,Z}X\;+\;\frac1{12}\,(\nabla_XR)_{X,Z}X
  \;+\;\Big(\,\frac1{60}\,(\nabla^2_{X,X}R)_{X,Z}X\,-\,\frac1{45}\,
  R_{X,R_{X,Z}X}X\,\Big)
  \\
  &&
  +\;\Big(\,\frac1{360}\,
  (\nabla^3_{X,X,X}R)_{X,Z}X\,-\,\frac1{120}\,R_{X,(\nabla_XR)_{X,Z}X}X
  \,-\,\frac1{120}\,(\nabla_XR)_{X,R_{X,Z}X}X\,\Big)\;+\;O(\,X^6\,)
 \end{eqnarray*}
 reflects the terms of the power series $\Theta$ up to order $5$. A simple
 algorithm to calculate the asymptotic expansion of $\Theta$ to arbitrary
 order is based on the asymptotic expansion
 \begin{equation}\label{aef}
  f(\;\exp_pX\;)
  \;\;\raise-6pt\hbox{$\stackrel\sim{\scriptstyle X\,\to\,0}$}\;\;
  \sum_{k\,\geq\,0}\frac1{k!}\,(\,\nabla^k_{X,\,\ldots,\,X}f\,)(\;p\;)
 \end{equation}
 of the pull back of functions via the exponential map, which is the
 analogue for sections $f\,\in\,C^\infty M$ of the trivial line bundle
 of the expansions (\ref{aex}) used in the proof of Lemma \ref{pte}.
 According to the definition of the exponentially extended vector field
 $Z^{\exp}$ the equality
 $$
  \left.\frac d{dt}\right|_0
  \exp_{\gamma(t)}\Big(\;\PT^\nabla(\,t\,)\,X\;\Big)
  \;\;=\;\;
  \left.\frac d{dt}\right|_0
  \exp_p\Big(\;X\;+\;t\,Z^{\exp}(\,X\,)\;\Big)
  \;\;\in\;\;
  T_{\exp_pX}M
 $$
 of tangent vectors holds true for all $X\,\in\,T_pM$ and all curves
 $\gamma$ representing $Z\,=\,\left.\frac d{dt}\right|_0\gamma$, in turn
 the asymptotic expansion (\ref{aef}) implies an equality of asymptotic
 expansions of the form
 \begin{eqnarray*}
  \left.\frac d{dt}\right|_0
  f(\;\exp_{\gamma(t)}\PT^\nabla(\,t\,)X\;)
  &\raise-6pt\hbox{$\stackrel\sim{\scriptstyle X\,\to\,0}$}&
  \left.\frac d{dt}\right|_0\sum_{k\,\geq\,0}\frac1{k!}
  (\,\nabla^k_{\PT^\nabla(t)X,\,\ldots,\,\PT^\nabla(t)X}f\,)
  (\;\gamma(\,t\,)\;)
  \\
  &=&
  \sum_{k\,\geq\,0}\frac1{k!}(\,\nabla^{k+1}_{Z,\,X,\,\ldots,\,X}f\,)(\;p\;)
  \\
  &\raise-6pt\hbox{$\stackrel\sim{\scriptstyle X\,\to\,0}$}&
  \left.\frac d{dt}\right|_0\sum_{k\,\geq\,0}\frac1{k!}
  (\,\nabla^k_{X+tZ^{\exp}(X),\,\ldots,\,X+tZ^{\exp}(X)}f\,)(\;p\;)
 \end{eqnarray*}
 The algorithm to calculate the Taylor series of the exponentially extended
 vector field $Z^{\exp}$ in $0\,\in\,T_pM$ is thus based on calculating the
 unique power series $\Theta^{-1}$ satisfying the identity
 $$
  \left.\frac d{dt}\right|_0
  \sum_{k\,\geq\,0}\frac1{k!}\;\nabla^k_{X\,+\,tZ,\,\ldots,\,X+tZ}f
  \;\;=\;\;
  \sum_{k\,\geq\,0}\frac1{k!}\;
  \nabla^{k+1}_{\Theta^{-1}(\,X\,)Z,\,X,\,\ldots,\,X}f
 $$
 of power series in $X\,\in\,T_pM$ for every function $f\,\in\,C^\infty M$,
 after formally inverting this power series we obtain the power series
 $\Theta$ describing the Taylor series of exponentially extended vector
 fields via $Z^{\exp}(X)\,\sim\,\Theta(X)\,Z$. The details of this algorithm
 can be worked out in analogy to the calculation of the backward parallel
 transport $\Phi^{-1}$ in \cite{w1}, somewhat hidden in equation (4.8) of
 the same reference the reader may find the relatively explicit formula
 \begin{equation}\label{eext}
  N\,(\,N\,-\,1\,)\,\Theta
  \;\;=\;\;
  N\,\Big(\;(\,\id\,-\,\Phi\,)\,\Phi^*\;\Big)\;+\;(\,N\,\Phi^{-1}\,)\,
  \Phi\,\Phi^*
 \end{equation}
 valid in the torsion free case. The power series $\Phi^*$ occuring in this
 formula is obtained by expanding the power series $\Phi$ in terms of the
 homogeneous components of the power series $\mathscr{R}$ considered as
 indeterminates and reverse the order of the factors in all monomials.
 Alternatively $\Phi^*$ is the inverse of the unique power series solution
 $\Phi^{-*}$ of the wrong sided parallel transport equation $N(N+1)\,\Phi^{-*}
 \,=\,\Phi^{-*}\,\mathscr{R}$ with initial value $\Phi^{-*}(\,0\,)\,=\,\id$.
\section{The Taylor Series of the Difference Element}
\label{element}
 In jet calculus the key idea is to find the appropiate difference elements
 and study the algebraic constraints they satisfy. Following this strategy
 we devise in this section an algorithm to calculate the infinite order
 Taylor series of the difference element $K\,:=\,\exp^{-1}\circ\knc_p$
 and its inverse in terms of the curvature tensor and its iterated covariant
 derivatives proving universality of the resulting expression on the way.
 Using this algorithm we calculate the lowest order terms of the Taylor
 series of the K\"ahler normal potential and the Riemannian distance function.
 Last but not least we define holomorphically extended vector fields and use
 them to determine the Spencer connection for the K\"ahler normal potential.

 \begin{Definition}[Difference Element for K\"ahler Normal Coordinates]
 \hfill\label{ekn}\break
  In order to compare the unique K\"ahler normal coordinates $\knc^{}_p:
  \,T_pM\longrightarrow M$ centered in point $p\,\in\,M$ of a K\"ahler
  manifold $M$ with the exponential map $\exp_p:\,T_pM\longrightarrow M$
  associated to the Levi--Civita connection $\nabla$ we consider their
  difference as a smooth map:
  $$
   K:\quad T_pM\;\longrightarrow\;T_pM,
   \qquad X\;\longmapsto\;\exp^{-1}_p\Big(\;\knc^{}_pX\;\Big)
  $$
  The analogue of the backward parallel transport $\Phi^{-1}$ in Riemannian
  normal coordinates is the K\"ahler backward parallel transport $\Psi^{-1}
  \,\in\,C^\infty(\,T_pM,\,\End\,T_pM\,)$ defined by means of:
  $$
   \Psi^{-1}(\,X\,)\,Y
   \;\;:=\;\;
   \Phi^{-1}(\,KX\,)\;DK(\,X\,)\;Y
   \;\;:=\;\;
   \Phi^{-1}(\,KX\,)\;\left.\frac d{dt}\right|_0K(\,X\,+\,t\,Y\,)
  $$
 \end{Definition}

 \noindent
 Recall that the main role of the backward parallel transport $\Phi^{-1}$
 in affine exponential coordinates is to describe the differential of the
 exponential map in covariant terms. In complete analogy the K\"ahler
 backward parallel transport $\Psi^{-1}$ describes the differential of
 K\"ahler normal coordinates $\knc_p:\,T_pM\longrightarrow M$ as a power
 series of covariant tensors. More precisely the implicit form $\knc_p\,=\,
 \exp_p\,\circ\,K$ of the definiton of the difference element $K$ implies
 $(\knc_p)_{*,\,X}\,=\,(\exp_p)_{*,\,KX}\circ K_{*,\,X}$ on differentials
 leading to the commutative diagram
 \begin{center}
  \begin{picture}(400,106)(0,0)
   \put(162, 93){\vector(+1, 0){70}}\put(174, 98){$(\knc_p)_{*,\,X}$}
   \put(148, 85){\vector(+1,-1){27}}\put(132, 66){$K_{*,\,X}$}
   \put(213, 58){\vector(+1,+1){27}}\put(233, 66){$(\exp_p)_{*,\,KX}$}
   \put( 27, 93){\vector(+1, 0){75}}\put( 60, 98){$\cong$}
   \put(284, 93){\vector(+1, 0){85}}\put(295, 98){$\PT^\nabla(KX)^{-1}$}
   \put(194, 39){\vector( 0,-1){28}}\put(198, 22){$\cong$}
   \put(194, 11){\vector( 0,+1){28}}
   \put( 27, 85){\vector(+2,-1){153}}\put( 54, 38){$DK(X)$}
   \put(212,  9){\vector(+2,+1){153}}\put(298, 38){$\Phi^{-1}(KX)$}
   \put(  0, 90){$T_pM$}
   \put(106, 90){$T_X(\,T_pM\,)$}
   \put(237, 90){$T_{\knc_pX}M$}
   \put(373, 90){$T_pM$}
   \put(165, 45){$T_{KX}(\,T_pM\,)$}
   \put(184,  0){$T_pM$}
  \end{picture}
 \end{center}
 which allows us to identify the power series $\Psi^{-1}$ from Definition
 \ref{ekn} with the composition:
 \begin{equation}\label{kpti}
  \Psi^{-1}(\,X\,):\;\;
  T_pM\;\cong\;T_X(\,T_pM\,)
  \;\stackrel{(\knc_p)_{*,X}}\longrightarrow\;
  T_{\knc_pX}M\;\stackrel{\PT^\nabla(KX)^{-1}}\longrightarrow\;T_pM
 \end{equation}
 In particular the pull back of the Riemannian metric $g$ to K\"ahler
 normal coordinates reads:
 \begin{equation}\label{pull}
  (\knc^*_pg)_X(\,A,B\,)
  \;\;=\;\;
  (\exp_p^*g)_{KX}(\,K_{*,X}A,K_{*,X}B\,)
  \;\;=\;\;
  g_p(\,\Psi^{-1}(X)A,\,\Psi^{-1}(X)B\,)
 \end{equation}
 Combined with equation (\ref{cgi}) this description of the Riemannian
 metric $g$ in K\"ahler normal coordinates can be used to calculate the
 Taylor series $\theta_p$ of the pull back $\knc_p^*\theta^\loc_p$ of the
 K\"ahler normal potential $\theta_p^\loc$ to the tangent space $T_pM$
 via $\knc_p$. According to our discussion of equation (\ref{cgi}) we
 are free to choose the connection $D$ among the torsion free connections
 making $I$ parallel in order to evaluate $\Hess\,=\,D\,\circ\,d$. With
 $\knc_p:\,T_pM\longrightarrow M$ being holomorphic with respect to the
 constant complex structure $I_p$ on $T_pM$ we may simply use the trivial
 connection $D$ for this purpose and so we are lead to consider the
 consequence
 \begin{eqnarray*}
  4\;(\,\knc_p^*g\,)_X(\;X,\;X\;)
  &=&
  4\;g_p(\;\Psi^{-1}(X)\,X,\;\Psi^{-1}(X)\,X\;)
  \\
  &=&
  (\,\Hess\,\theta_p\,)_X(\;X,\;X\;)\;+\;(\,\Hess\,\theta_p\,)_X(\;IX,\;IX\;)
 \end{eqnarray*}
 of equations (\ref{cgi}) and (\ref{pull}), which implies that the formal
 differential equation
 \begin{equation}\label{pott}
  [\;(\,N^2\,+\,\Der_I^2\,)\,\theta_p\;](\,X\,)
  \;\;=\;\;
  4\,g(\;\Psi^{-1}(X)\,X,\;\Psi^{-1}(X)\,X\;)
 \end{equation}
 is obeyed by the Taylor series $\theta_p$ of the pull back $\knc_p^*
 \theta_p^\loc$ of the normal potential to the tangent space $T_pM$; in
 fact the Hessian with respect to the trivial connection $D$ satisfies
 \begin{eqnarray*}
  (\,\Hess\,\psi\,)_X(\;\,X\,,\;\;X\;)
  &=&
  [\;N\,(\,N\;-\;1\,)\;\psi\;](\;X\;)
  \\
  (\,\Hess\,\psi\,)_X(\,IX,\,IX\,)
  &=&
  [\;(\,\Der_I^2\,+\,N\,)\,\psi\;](\;X\;)
 \end{eqnarray*}
 for every polynomial $\psi$ on $T_pM$. Besides determining the pull backs
 of the Riemannian metric $g$ and the K\"ahler normal potential $\theta^\loc$
 respectively to K\"ahler normal coordinates via equations (\ref{pull}) and
 (\ref{pott}) the difference element $K$ offers a very direct description
 of the square of the Riemannian distance function in K\"ahler normal
 coordinates as well:

 \begin{Remark}[Distance Function in K\"ahler Normal Coordinates]
 \hfill\label{dfk}\break
  In Riemannian geometry the exponential map is a radial isometry in the
  sense that the Riemannian distance between $p\,\in\,M$ and $\exp_pX$
  equals $\dist^2_g(\,p,\,\exp_pX\,)\,=\,g_p(\,X,\,X\,)$ for $X$ sufficiently
  small. In K\"ahler normal coordinates this property of the distance
  becomes:
  $$
   \dist^2_g(\;p,\;\knc^{}_pX\;)
   \;\;=\;\;
   \dist^2_g(\;p,\;\exp_p(\,KX\,)\;)
   \;\;=\;\;
   g_p(\,KX,\,KX\,)
  $$
 \end{Remark}

 \noindent
 Let us now come to the more difficult task of finding an efficient method
 to calculate the Taylor series of the difference elements $K$ and $K^{-1}$.
 For that purpose we recall the fact that a vector field $X\,\in\,\Gamma
 (\,TM\,)$ on a manifold $M$ with integrable almost complex structure
 $I$ is the real part of a holomorphic vector field in the sense that
 the map $X-iIX:\,M\longrightarrow T^{1,0}M$ between complex manifolds
 is actually holomorphic, if and only if
 $$
  \mathfrak{Lie}_XI\;\;=\;\;0
 $$
 compare for example \cite{ba} and \cite{mor}. In order to construct a
 resolution for the symbolic differential operator $\D$ corresponding
 to this complex Killing equation in the formal theory of partial
 differential equations \cite{bc3}, \cite{xd} we need the concept
 of special alternating forms:

 \begin{Definition}[Special Alternating Forms]
 \hfill\label{saf}\break
  Consider the graded vector space $\Lam^\circ T^*\otimes T$ of alternating
  forms on a real vector space $T$ with values in $T$. The choice of a
  complex structure $I\,\in\,\End\,T$ on $T$ singles out the graded subspace
  $\Sigma^\circ\,\subset\,\Lam^\circ T^*\otimes T$ of special alternating
  forms defined in all degrees $k\,\in\,\N_0$ by:
  $$
   \Sigma^k
   \;\;:=\;\;
   \{\;\;F\,\in\,\Lam^kT^*\otimes T\;\;|\;\;(\,\Der_I\,\otimes\,I\,)
   \,F\;=\;k\,F\;\;\}
  $$
  Alternatively a special alternating $k$--form $F:\,T\times\ldots\times T
  \longrightarrow T$ is characterized by
  $$
   I\,F(\;X_1,\,\ldots,\,X_k\;)
   \;\;=\;\;
   -\,F(\;IX_1,\,X_2,\,\ldots,\,X_k\;)
  $$
  in particular $\Sigma^1\,\subset\,T^*\otimes T$ is simply the subspace of
  endomorphisms anticommuting with $I$:
  $$
   \pr_{\Sigma^1}:\quad T^*\,\otimes\,T\;\longrightarrow\;\Sigma^1,
   \qquad F\;\longmapsto\;{\textstyle\frac12}(\,F\,+\,I\,F\,I\,)
  $$ 
 \end{Definition}

 \noindent
 The link between the complex Killing equation and special alternating
 forms is established directly by rewriting the definition of the Lie
 derivative $\mathfrak{Lie}_XI$ in terms of a torson free connection
 $\nabla$ on the tangent bundle making $I$ parallel $\nabla I\,=\,0$.
 On a K\"ahler manifold $M$ we take for example the  Levi--Civita connection
 and find for the Lie derivative of the orthogonal complex structure $I$
 in the direction of a vector field $X\,\in\,\Gamma(TM)$:
 \begin{eqnarray*}
  (\,\mathfrak{Lie}_XI\,)\,Y
  &=&
  \Big(\;\nabla_X(\,IY\,)\;-\;\nabla_{IY}X\;\Big)\;-\;I\,
  \Big(\;\nabla_XY\;-\;\nabla_YX\;\Big)
  \\
  &=&
  I\,\Big(\;\nabla_YX\;+\;I\,\nabla_{IY}X\;\Big)
  \;\;=\;\;
  2\,I\,\pr_{\Sigma^1}(\,\nabla X\,)\;Y
 \end{eqnarray*}
 Hence the real parts of holomorphic vector fields are precisely the sections
 $X\,\in\,\Gamma(\,TM\,)$ in the kernel of the differential operator
 $\D_{\mathrm{diff}}$ defined as the composition of the covariant
 derivative $\nabla$ with the projection $\pr_{\Sigma^1}$ to the
 subbundle $\Sigma^1M\,\subset\,\End\,TM$
 \begin{equation}\label{dd}
  \D_{\mathrm{diff}}\,:\qquad
  \Gamma(\,TM\,)
  \;\stackrel\nabla\longrightarrow\:
  \Gamma(\,T^*M\otimes TM\,)
  \;\stackrel{\pr_{\Sigma^1}}\longrightarrow\;
  \Gamma(\,\Sigma^1M\,),
  \qquad X\;\longmapsto\;-\,\frac12\,I\,\mathfrak{Lie}_XI
 \end{equation}
 of endomorphisms anticommuting with $I$. From the definition of
 $\D_{\mathrm{diff}}$ as the composition $\pr_{\Sigma^1}\circ\nabla$
 we read off its principal symbol and in turn the associated symbol
 comodule \cite{xd}:
 $$
  \H^\bullet
  \;\;:=\;\;
  \ker\Big(\,\D:\;\;\S^\bullet T^*\otimes\Sigma^0
  \;\stackrel\Delta\longrightarrow\;
  \S^{\bullet-1}T^*\otimes(T^*\otimes T)
  \;\stackrel{\id\otimes\pr_{\Sigma^1}}\longrightarrow\;
  \S^{\bullet-1}T^*\otimes\Sigma^1\,\Big)
 $$
 The symbolic differential operator $\D$ defining the comodule $\H$ extends
 to a complete resolution of the comodule $\H$ by free comodules using an
 adjoint pair of boundary operators:
 
 \begin{Definition}[Symbolic Differential Operators]
 \hfill\label{bo}\break
  Consider a real vector space $T$ endowed with a complex structure
  $I\,\in\,\End\,T$. On the bigraded vector space $\S^\bullet T^*\otimes
  \Sigma^\circ\,\subset\,\S^\bullet T^*\otimes\Lam^\circ T^*\otimes T$ of
  special alternating forms on $T$ with polynomial coefficients we define
  two bigraded boundary operators $\D^*$ and $\D$ by:
  \begin{eqnarray*}
   [\hbox to27pt{\hfill$\D^*F$\hfill}]_X(Z_2,\ldots,Z_r)
   &:=&
   \hphantom{\frac12}F_X(\;X,\,Z_2,\,\ldots,\,Z_r\;)
   \\[2pt]
   [\hbox to27pt{\hfill$\D\,F$\hfill}]_X(Z_0,\ldots,Z_r)
   &:=&
   \frac12\sum_{\mu\,=\,0}^r(-1)^\mu\Big(\,(\,Z_\mu\,
   \ins\,F\,)_X(Z_0,\ldots,Z_r)\,+\,I\,(IZ_\mu\ins F)_X(Z_0,\ldots,Z_r)
   \,\Big)
  \end{eqnarray*}
 \end{Definition}

 \noindent
 In order to clarify this definition we recall the convention adopted in this
 article concerning the identification of a symmetric $k$--multilinear form
 $F\,=\,\S^kT^*$ on $T$ with a homogeneous polynomial on $T$ of degree
 $k\,\in\,\N_0$: The polynomial corresponding to $F$ is defined by $F(X)
 \,:=\,\frac1{k!}F(X,\ldots,X)$ so that the operation $Z\,\ins\,$ of
 inserting the first argument agrees with the directional derivative
 $\frac\partial{\partial Z}$ of the polynomial in direction $Z$. With
 this convention in place we may alternatively define the boundary operators
 $\D^*$ and $\D$ in terms of the tensor product decomposition
 $\S^\bullet T^*\otimes\Sigma^\circ\,\subset\,\S^\bullet T^*\otimes
 \Lam^\circ T^*\otimes T$ as sums
 \begin{eqnarray*}
  \D^*
  &:=&
  \hphantom{\frac12}
  \sum_\mu dE_\mu\,\cdot\,\otimes\,E_\mu\,\ins\,\otimes\,\id
  \\[-1pt]
  \D&:=&
  \frac12\sum_\mu\Big(\,E_\mu\,\ins\,\otimes\,dE_\mu\,\wedge\,\otimes\,\id
  \;+\;IE_\mu\,\ins\,\otimes\,dE_\mu\,\wedge\,\otimes\,I\;\Big)
 \end{eqnarray*}
 over a dual pair of bases $\{\,E_\mu\,\}$ and $\{\,dE_\mu\,\}$ of $T$ and
 $T^*$. Alternatively we could complexify the domain $\S^\bullet T^*\otimes
 \Lam^\circ T^*\otimes T$ of the boundary operator $\D$ and choose a complex
 basis $\{\,F_\alpha\,\}$ of the subspace $T^{1,0}\,\subset\,T\otimes_\R\C$
 with dual basis $\{\,dF_\alpha\,\}$ of $T^{1,0*}\,\subset\,T^*\otimes_\R\C$
 to rewrite
 $$
  \D
  \;\;:=\;\;
  \sum_\alpha\Big(\;
  F_\alpha\,\ins\,\otimes\,dF_\alpha\,\wedge\,\otimes\,\pr^{0,1}
  \;+\;\overline{F}_\alpha\,\ins\,\otimes\,d\overline{F}_\alpha\,
  \wedge\,\otimes\,\pr^{0,1}\;\Big)
 $$
 with the eigenprojections $\pr^{1,0}\,:=\,\frac12(\id-iI)$ and
 $\pr^{0,1}\,:=\,\frac12(\id+iI)$ to $T^{1,0}$ and $T^{0,1}$, in this
 alternative formulation the operator $\D$ evidently preserves the subspace
 of special alternating forms with polynomial coefficients. The reader is
 invited to use this description of the operator $\D$ in order to provide
 a more enlightening proof of Corollary \ref{frhv} below. Staying in the
 real domain we calculate the anticommutator $\{\,\D,\,\D^*\,\}$ of the
 boundary operators $\D$ and $\D^*$ using the canonical commutation and
 anticommutation relations
 \begin{eqnarray*}
  \{\,\D,\,\D^*\,\}
  &=&
  \frac12\sum_{\mu\,\nu}\Big(\,[\,E_\mu\ins,\,dE_\nu\cdot\,]\otimes
  dE_\mu\wedge E_\nu\ins\otimes\id\,+\,dE_\nu\cdot E_\mu\ins\otimes
  \{\,dE_\mu\wedge,\,E_\nu\ins\,\}\otimes\id
  \\[-8pt]
  &&
  \qquad+\;[\,IE_\mu\ins,\,dE_\nu\cdot\,]\otimes
  dE_\mu\wedge E_\nu\ins\otimes I\,+\,dE_\nu\cdot IE_\mu\ins\otimes
  \{\,dE_\mu\wedge,\,E_\nu\ins\,\}\otimes I\,\Big)
  \\
  &=&
  \frac12\,(\;\id\,\otimes\,N\,\otimes\,\id\;+\;N\,\otimes\,\id\,\otimes\,\id
  \;+\;\id\,\otimes\,\Der_I\,\otimes\,I\;+\;\Der_I\,\otimes\,\id\,\otimes\,I\;)
 \end{eqnarray*}
 where $N$ denotes the number operator either on polynomials or on forms
 depending on context. The equality of linear maps $\Der_I\otimes I\,=\,
 N\otimes\id$ defining the subspace of special alternating forms $\Sigma^\circ
 \,\subset\,\Lam^\circ T^*\otimes T$ allows us to write this identity in the
 simpler form:
 \begin{equation}\label{lap}
  \Delta
  \;\;:=\;\;
  \{\;\D,\;\D^*\;\}
  \;\;=\;\;
  \id\,\otimes\,N\,\otimes\,\id\;+\;\frac12\,
  (\;N\,\otimes\,\id\,\otimes\,\id\;+\;\Der_I\,\otimes\,\id\,\otimes\,I\;)
 \end{equation}
 A short inspection reveals that the formal Laplace operator $\Delta$ acts
 diagonalizable on the space $\S^\bullet T^*\otimes\Sigma^\circ$ of special
 alternating forms with polynomial coefficients with eigenspaces
 $$
  (\,\S^{\kappa,\,\overline\kappa}T^*\otimes\Lam^dT^{0,1*}\otimes T^{1,0}\,)
  \;\oplus\;
  (\,\S^{\overline\kappa,\,\kappa}T^*\otimes\Lam^dT^{1,0*}\otimes T^{0,1}\,)
  \;\;\subset\;\;
  (\,\S^\bullet T^*\,\otimes\,\Sigma^d\,)\,\otimes_\R\C
 $$
 and eigenvalues $d\,+\,\overline\kappa$ parametrized by $\kappa,\,
 \overline\kappa\,\geq\,0$ and $d\,\geq\,0$ with $\bullet\,=\,
 \kappa+\overline\kappa$. With $\Delta$ being diagonalizable we
 conclude from Hodge theory that the homology of the $\D$--complex
 equals the homology of the eigensubcomplex $\ker\,\Delta$, because
 $\frac1\lambda\D^*$ is a zero homotopy for the eigensubcomplexes
 $\ker\,(\,\Delta-\lambda\,\id\,)$ for all eigenvalues $\lambda\,\neq\,0$.
 On the other hand $\ker\,\Delta$ equals
 $$
  (\,\S^{\bullet,\,0}T^*\otimes T^{1,0}\,)
  \;\oplus\;
  (\,\S^{0,\,\bullet}T^*\otimes T^{0,1}\,)
  \;\;=\;\;
  \H^\bullet\,\otimes_\R\C
  \;\;\subset\;\;
  (\,\S^\bullet T^*\otimes T\,)\,\otimes_\R\C
 $$
 and so $\D\,=\,0$ vanishes on $\ker\,\Delta$ for trivial reasons.
 Summarizing this argument we conclude:

 \begin{Corollary}[Free Resolution of Holomorphic Vector Fields]
 \hfill\label{frhv}\break
  The total prolongation comodule $\H^\bullet$ associated to the principal
  symbol $\pr_{\Sigma^1}:\,T^*\otimes T\longrightarrow\Sigma^1$ of the
  differential operator $\D_{\mathrm{diff}}$ characterizing the real parts
  of holomorphic vector fields in complex dimension $n\,:=\,\frac12\,\dim\,T$
  has a resolution by free comodules of the form:
  $$
   0\;\longrightarrow\;\H^\bullet
   \;\stackrel\subset\longrightarrow\;\S^\bullet T^*\otimes\Sigma^0
   \;\stackrel\D\longrightarrow\;\S^{\bullet-1}T^*\otimes\Sigma^1
   \;\stackrel\D\longrightarrow\;\ldots\;
   \;\stackrel\D\longrightarrow\;\S^{\bullet-n}T^*\otimes\Sigma^n
   \;\longrightarrow\;0
  $$
 \end{Corollary}

 \noindent
 An important property of the real parts of holomorphic vector fields
 needed below is that they can be reconstructed up to a linear term
 from their closure, their image under the map
 \begin{equation}\label{cl}
  \cl:\quad\S^\bullet T^*\,\otimes\,T\;\longrightarrow\;
  \S^{\bullet+1}T^*,\qquad Z\;\longmapsto\;\cl\,Z
 \end{equation}
 from vector fields to polynomials defined by $(\cl\,Z)(X)\,:=\,g(Z(X),X)$.
 Since the eigenspaces $T^{1,0}$ and $T^{0,1}$ of $I$ are isotropic subspaces
 of $T\otimes_\R\C$, the closure of the real part of a holomorphic vector
 field lies in the vector space $\S^\bullet_{[\,1\,]}T^*\,:=\,\ker(\,(N-2)^2
 \,+\,\Der_I^2\,)$ of polynomials linear in the holomorphic or antiholomorphic
 coordinates with complexification:
 \begin{equation}\label{1s}
  \S^\bullet_{[\,1\,]}T^*\,\otimes_\R\C
  \;\;=\;\;
  \bigoplus_{{\scriptstyle\kappa,\,\overline\kappa\,\geq\,0}\atop
  {\scriptstyle\kappa\,=\,1\;\mathrm{or}\;\overline\kappa\,=\,1}}
  \S^{\kappa,\,\overline\kappa}T^*
 \end{equation}
 Restricted to the real parts of holomorphic vector fields $\cl$ induces
 in fact an exact sequence
 \begin{equation}\label{ses}
  0\;\longrightarrow\;
  \delta_{\bullet\,=\,1}\;\mathfrak{u}(\,T,\,g,\,I\,)
  \;\stackrel\subset\longrightarrow\;
  \H^\bullet
  \;\stackrel{\mathrm{cl}}\longrightarrow\;
  \S^{\bullet+1}_{[\,1\,]}T^*
  \;\longrightarrow\;0
 \end{equation}
 in which $\mathfrak{u}(\,T,\,g,\,I\,)$ denotes the vector space of linear
 vector fields on $T$ corresponding to infinitesimal unitary transformations
 with respect to the hermitean form $h\,=\,g+i\omega$.

 \pfill
 Coming back to our original aim to derive recursion formulas for the
 infinite order Taylor series of the difference elements $K$ and $K^{-1}$
 we recall that multilinear maps on a vector space with values in this
 vector space can be composed to produce new multilinear forms
 $$
  (\,A\,\circ_\mu B\,)(\,X_1,\ldots,X_{a+b+1}\,)
  \;\;:=\;\;
  A(\,X_1,\ldots,,X_{\mu-1},\,B(\,X_\mu,\ldots,X_{\mu+b}\,),\,X_{\mu+b+1},
  \ldots,X_{a+b+1}\,)
 $$
 for some position $\mu\,=\,1,\ldots,a$. We will call a multilinear form
 arising from a fixed set of basic multilinear forms by iterated compositions
 in arbitrary positions a composition polynomial in the selected set of basic
 multilinear forms. In particular we will consider composition polynomials in
 the complex structure $I_p$, the curvature tensor $R_p$ and its iterated
 covariant derivatives $(\nabla R)_p,\,(\nabla^2R)_p,\,\ldots$ on the tangent
 space $T_pM$ of a K\"ahler manifold $M$ in a point $p\,\in\,M$. In general
 multilinear forms on $T_pM$ are bigraded by degree and weight: The degree of
 a multilinear form is the number of arguments it takes minus $1$ to make the
 degree additive under composition, while its weight is the eigenvalue for
 the weight operator:
 \begin{equation}\label{wo}
  (\,\delta A\,)(\,X_1,\ldots,X_k\,)
  \;\;:=\;\;
  I\,A(\,X_1,\,\ldots,\,X_k\,)
  \;-\;(\,\Der_IA\,)(\,X_1,\,\ldots,\,X_k\,)
 \end{equation}
 The weight qualifies as a grading, because the weight operator is a
 derivation for composition
 $$
  \delta\,(\,A\,\circ_\mu B\,)
  \;\;=\;\;
  (\,\delta A\,)\,\circ_\mu B\;+\;A\,\circ_\mu(\,\delta B\,)
 $$
 as the term $A(\,X_1,\ldots,X_{\mu-1},\,IB(X_\mu,\ldots,X_{\mu+b}),
 \,X_{\mu+b+1},\ldots,X_{a+b+1}\,)$ appears twice with different sign in
 the expansion of the right hand side. The curvature tensor $R$ of the
 K\"ahler manifold $M$ for example is a degree $2$ multilinear form on
 $T_pM$ of weight $0$
 $$
  (\,\delta R\,)_{X,\,Y}Z
  \;\;=\;\;
  R_{X,\,Y}IZ\;-\;R_{I^2X,\,IY}Z\;-\;R_{X,\,IY}Z\;-\;R_{X,\,Y}IZ
  \;\;=\;\;
  0
 $$
 in consequence its iterated covariant derivatives $\nabla^rR$ decompose
 after complexification into a sum of multilinear forms of degree $r+2$ and
 weight $-ri,\,\ldots,\,+ri$. As both degree and weight are additive under
 composition a quadratic composition polynomial in the curvature tensor and
 its iterated covariant derivatives needs to have a least degree $r+4$ to
 decompose similarly into multilinear forms of weight $-ri,\,\ldots,\,+ri$
 and so on. The possible combinations of degree and weight of multilinear
 forms in general and composition polynomials in the curvature tensor and
 its iterated covariant derivatives can be read off from the diagram:
 \begin{equation}\label{diag}
  \vcenter
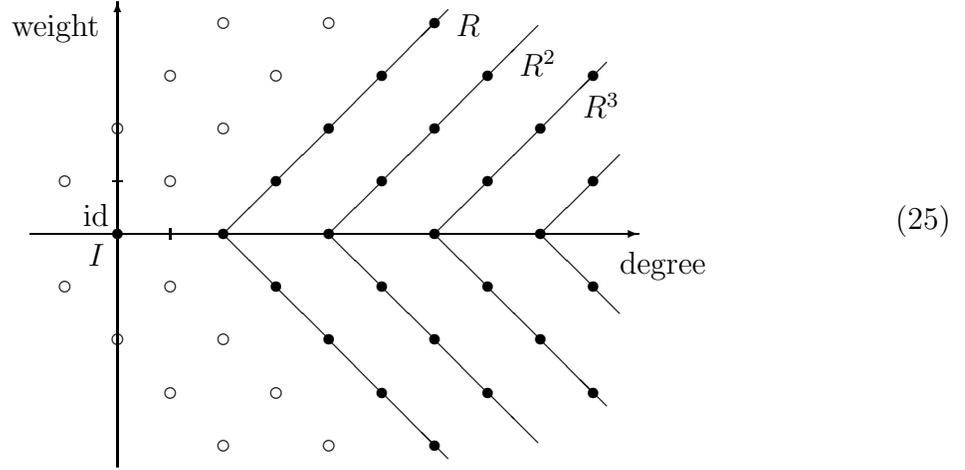
{\hbox{\begin{picture}(240,170)(0,0)
   \put( 40, -8){\vector( 0,+1){176}}\put(  0,155){$\mathrm{weight}$}
   \put(  7, 80){\vector(+1, 0){230}}\put(230, 65){$\mathrm{degree}$}
   \put( 60, 78){\line( 0,+1){4}}\put( 38,100){\line(+1, 0){4}}
   \put( 27, 83){$\id$}\put( 29, 68){$I$}
   \put(168,155){$R$}\put(192,140){$R^2$}\put(216,124){$R^3$}
   \put( 80, 80){\line(+1,+1){85}}\put( 80, 80){\line(+1,-1){85}}
   \put(120, 80){\line(+1,+1){79}}\put(120, 80){\line(+1,-1){79}}
   \put(160, 80){\line(+1,+1){65}}\put(160, 80){\line(+1,-1){65}}
   \put(200, 80){\line(+1,+1){30}}\put(200, 80){\line(+1,-1){30}}
   \put( 80,  0){\circle4} \put(120,  0){\circle4} \put(160,  0){\circle*4}
   \put( 60, 20){\circle4} \put(100, 20){\circle4} \put(140, 20){\circle*4}
   \put(180, 20){\circle*4}\put(220, 20){\circle*4}\put( 40, 40){\circle4}
   \put( 80, 40){\circle4} \put(120, 40){\circle*4}\put(160, 40){\circle*4}
   \put(200, 40){\circle*4}\put( 20, 60){\circle4} \put( 60, 60){\circle4}
   \put(100, 60){\circle*4}\put(140, 60){\circle*4}\put(180, 60){\circle*4}
   \put(220, 60){\circle*4}\put( 40, 80){\circle*4}\put( 80, 80){\circle*4}
   \put(120, 80){\circle*4}\put(160, 80){\circle*4}\put(200, 80){\circle*4}
   \put( 20,100){\circle4} \put( 60,100){\circle4} \put(100,100){\circle*4}
   \put(140,100){\circle*4}\put(180,100){\circle*4}\put(220,100){\circle*4}
   \put( 40,120){\circle4} \put( 80,120){\circle4} \put(120,120){\circle*4}
   \put(160,120){\circle*4}\put(200,120){\circle*4}\put( 60,140){\circle4}
   \put(100,140){\circle4} \put(140,140){\circle*4}\put(180,140){\circle*4}
   \put(220,140){\circle*4}\put( 80,160){\circle4} \put(120,160){\circle4}
   \put(160,160){\circle*4}
  \end{picture}}}
 \end{equation}
 \pfill
 These considerations play a crucial role in the proof of the following
 theorem:

 \begin{Theorem}[Universality of K\"ahler Normal Coordinates]
 \hfill\label{univ}\break
  Every term in the Taylor series of the difference elements $K\,=\,
  \exp_p^{-1}\circ\knc_p$ and $K^{-1}$ in the origin $0\,\in\,T_pM$
  is a universal composition polynomial in the complex structure $I$, the
  curvature tensor $R$ and its iterated covariant derivatives $\nabla R,\,
  \nabla^2R,\ldots$ evaluated in $p\,\in\,M$.
 \end{Theorem}

 \noindent
 In this context universality refers to the statement that these
 composition polynomials can be chosen without taking the K\"ahler
 manifold $M$ or its dimension into account. Universality does not
 imply uniqueness of course, and actually uniqueness wouldn't make
 sense anyhow:

 \begin{Remark}[Non--Uniqueness of Composition Polynomials]
 \hfill\label{scc}\break
  Composition polynomials in the complex structure, the curvature tensor
  and its covariant derivatives $I,\,R,\,\nabla R,\,\ldots$ are in general
  not unique due to additional identities arising from Ricci type constraints.
  The following identity for example is valid on K\"ahler manifolds $M$
  \begin{eqnarray*}
   (\,\nabla^2_{X,\,IX}R\,)_{X,\,IX}
   \;-\;(\,\nabla^2_{IX,\,X}R\,)_{X,\,IX}
   &=&
   [\,R_{X,\,IX},\,R_{X,\,IX}\,]
   \;-\;R_{R_{X,\,IX}X,\,IX}\;-\;R_{X,\,R_{X,\,IX}IX}
   \\
   &=&
   -\,2\,R_{R_{X,\,IX}X,\,IX}
  \end{eqnarray*}
  for all $X\,\in\,TM$ and implies $R_{X,\,R_{X,\,IX}IX}\,=\,0$ on hermitean
  locally symmetric spaces.
 \end{Remark}

 \noindent\textbf{Proof of Theorem \ref{univ}:\quad}
 Both $\exp_p$ and $\knc_p$ are anchored coordinates for a K\"ahler manifold
 $M$, hence the differential of $K$ in $0\,\in\,T_pM$ equals the identity.
 In turn the decomposition of the infinite order Taylor series of $K$ into
 homogeneous components reads
 $$
  KX
  \;\;=\;\;
  X\,+\,K_2X\,+\,K_3X\,+\,\ldots
  \qquad\quad
  K_n
  \;\;\in\;\;
  \S^nT^*_pM\otimes T_pM
  \;\;=\;\;
  \H^n\,\oplus\,\H^{n\perp}
 $$
 where the complement $\H^{n\perp}$ to $\H^n\,=\,\ker\,\Delta$ is simply
 the direct sum of all eigenspaces of $2\Delta\,=\,(\,N\otimes\id\,+\,
 \Der_I\otimes I\,)$ with non--vanishing eigenvalue and thus the proper
 domain of:
 \begin{equation}\label{pi} 
  (\;N\,\otimes\,\id\;+\;\Der_I\,\otimes\,I\;)^{-1}:
  \quad \H^{n\perp}\;\longrightarrow\;\H^{n\perp}
 \end{equation}
 Essentially the proof consists in verifying by induction on the degree $n$
 that the characteristic normalization constraint imposed on the K\"ahler
 normal potential $\theta$ is actually equivalent to the normalization
 constraint $K_n\,\in\,\H^{n\perp}$ for all $n\,\geq\,2$.

 In a first step we want to analyse the implications of the holomorphicity
 of K\"ahler normal coordinates $\knc_p$ on the difference element $K\,=\,
 \exp_p^{-1}\,\circ\,\knc_p$. On a K\"ahler manifold the complex structure
 $I$ is a parallel tensor and so Corollary \ref{taypar} implies for the
 infinite order Taylor series for its pull back $\exp_p^*I$ to the tangent
 space under the exponential map
 $$
  (\,\exp^*_pI\,)_X
  \;\;\raise-6pt\hbox{$\stackrel\sim{\scriptstyle X\,\to\,0}$}\;\;
  \Phi(\,X\,)\,\circ\,I_p\,\circ\,\Phi^{-1}(\,X\,)
 $$
 whereas the holomorphicity of $\knc_p$ is equivalent to the statement
 $(\knc^*_pI)_X\,=\,I_p$ for all $X\,\in\,T_pM$. Hence the differential
 $K_{*,\,X}:\,T_X(T_pM)\longrightarrow T_{KX}(T_pM)$, which becomes $DK(X)$
 under the identification of both domain and target with $T_pM$, intertwines
 $I_p$ with $(\exp^*_pI)_{KX}$:
 $$
  \Phi(\,KX\,)\,\circ\,I_p\,\circ\,\Phi(\,KX\,)^{-1}\,\circ\,DK(\,X\,)
  \;\;=\;\;
  DK(\,X\,)\,\circ\,I_p
 $$
 Dropping the explicit mention of the point $p\,\in\,M$ on $I\,=\,I_p$
 we may write this equation
 \begin{equation}\label{inter}
  [\;I,\;\Phi^{-1}(\,KX\,)\,\circ\,DK(\,X\,)\;]
  \;\;=\;\;
  0
  \;\;=\;\;
  [\;I,\;DK^{-1}(\,X\,)\,\circ\,\Phi(\,X\,)\;]
 \end{equation}
 where we replaced $X$ by $K^{-1}X$ in the second equality. Note that this
 result is compatible with the interpretation of the K\"ahler backward
 parallel transport $\Psi^{-1}(X)\,:=\,\Phi^{-1}(KX)\,\circ\,DK(X)$ as
 the differential of the K\"ahler normal coordinates $\knc_p$.

 For a moment let us forget that $K$ is a diffeomorphism and think of it as
 a vector field with Taylor series $K\,\in\,\overline\S\,T^*_pM\otimes T_pM$.
 By construction the symbolic differential operator $\D$ equals $\D K\,=\,
 \pr_{\Sigma^1}(\,DK\,)$ on vector fields, whereas $\D^*K\,=\,0$. Applying
 the anticommutator $\Delta\,=\,\{\,\D,\,\D^*\,\}$ to the vector field $K$
 we thus obtain the decisive formula
 \begin{equation}\label{orec}
  \Delta K
  \;\;=\;\;
  \D^*\,\pr_{\Sigma^1}\Big(\;DK\;\Big)
  \;\;=\;\;
  \D^*\,\pr_{\Sigma^1}\Big(\;(\,\id\,-\,[\,\Phi^{-1}\circ K\,]\,)\,DK\;\Big)
 \end{equation}
 where $[\,\Phi^{-1}\circ K\,]\,\circ\,DK$ commutes with $I$ and thus
 vanishes under $\pr_{\Sigma^1}$ according to (\ref{inter}). The resulting
 formula is actually a recursion formula, which allows us to calculate the
 term $K_n$ in the Taylor series of $K$ from the terms $K_2,\,\ldots,\,
 K_{n-2}$ up to addition by an element of $\ker\,\Delta$, because
 $\id-\Phi^{-1}(KX)\,=\,-\frac16\,R_{X,\,\cdot}X\,+\,O(X^3)$ is at
 least quadratic in $X$.

 In the ensuing induction on the degree $n\,\geq\,2$ is more convenient to
 begin with the inductive step. By induction hypothesis we may thus assume
 that for some $n\,\geq\,3$ all the terms $K_2,\,\ldots,\,K_{n-1}$ are
 universal composition polynomials in the complex structure $I$, the
 curvature tensor and its iterated covariant derivatives $R,\,\nabla R,
 \,\ldots$. Under this assumption the homogeneous term of degree $n$ in
 $X$ on the right hand side of the recursion formula (\ref{orec}) is a
 universal composition polynomial in these generators and so then is the
 corresponding term on the left hand side. The formal Laplace operator
 $2\Delta\,=\,(N\otimes\id+\Der_I\otimes I)$ on the other hand has only
 finitely many non--zero eigenvalues for fixed degree $N\,=\,n$ so that
 its partial inverse (\ref{pi}) can be written as a polynomial in
 $\Der_I\otimes I$. In consequence the partial inverse maps composition
 polynomials in $I,\,R,\,\nabla R,\,\ldots$ to composition polynomials
 and so
 $$
  K_n
  \;\;=\;\;
  H_n\;+\;\textrm{composition polynomial in\ }
  I,\,R,\,\nabla R,\,\nabla^2R,\,\ldots
 $$
 for some $H_n\,\in\,\H^n$. Incidentally the same conclusion is valid in
 the special case $n\,=\,2$ forming the base of our induction, because the
 right hand side of the recursion formula (\ref{orec}) is at least cubic
 in $X$ as $\D^*$ raises the degree by $1$ and $\id-\Phi^{-1}(KX)\,=\,O(X^2)$.

 In order to verify the base of induction and complete the induction step we
 need to show $H_n\,=\,0$ for all $n\,\geq\,2$. For this purpose we calculate
 the term homogeneous of degree $n$ in $X$ in the K\"ahler backward parallel
 transport observing $DK(X)\,X\,=\,(NK)\,X$:
 \begin{eqnarray*}
  [\,\Psi^{-1}(\,X\,)\,X\,]_n
  &=&
  [\,\Phi^{-1}(\,KX\,)\,DK(\,X\,)\,X\,]_n
  \\
  &=&
  n\,H_nX\;+\;\textrm{composition polynomial in\ }
  I,\,R,\,\nabla R,\,\nabla^2R,\,\ldots
 \end{eqnarray*}
 Inserting this expression into equation (\ref{pott}) we obtain the
 expansion
 \begin{equation}\label{decis}
  [\;(\,N^2\,+\,\Der_I^2\,)\,\theta\;]_{n+1}(\,X\,)
  \;\;=\;\;
  8\,n\,g(\;H_nX,\;X\;)\;+\;\ldots
 \end{equation}
 for the homogeneous term of degree $n+1$ in $X$ in the K\"ahler normal
 potential, where the ellipsis denotes a finite sum of homogeneous
 polynomials of degree $n+1$ in $X$ of the form $g(\,\mathbb{A},\,
 \mathbb{B}\,)$ with composition polynomials $\mathbb{A}$ and $\mathbb{B}$
 in the complex structure $I$, the curvature tensor $R$ and its iterated
 covariant derivatives. For every polynomial of this form we find
 $$
  \Der_I\;[\;g(\;\mathbb{A},\;\mathbb{B}\;)\;]
  \;\;=\;\;
  g(\;\delta\,\mathbb{A},\;\mathbb{B}\;)\;+\;
  g(\;\mathbb{A},\;\delta\mathbb{B}\;)
 $$
 because the two additional terms $g(I\mathbb{A},\mathbb{B})\,+\,g(\mathbb{A},
 I\mathbb{B})$ on the right hand side obtained upon expanding $\delta$
 cancel out by the skew symmetry of $I$. Of course $\mathbb{A}$ and
 $\mathbb{B}$ can not be composition polynomials in $I$ only, since
 we would not get a polynomial of degree $n+1\,\geq\,3$, hence at least
 one of $\mathbb{A}$ or $\mathbb{B}$ is at least linear in $R,\,\nabla R,
 \,\nabla^2R,\,\ldots$. Diagram (\ref{diag}) thus tells us the possible
 combinations of degree and weight of the polynomial $g(\mathbb{A},
 \mathbb{B})$, from which we deduce:
 $$
  g(\;\mathbb{A},\;\mathbb{B}\;)
  \;\;\equiv\;\;0
  \qquad\textrm{mod}\;\;\S^{\geq(2,2)}T^*_pM
 $$
 In consequence equation (\ref{decis}) is a congruence modulo
 $\S^{\geq(2,2)}T^*_pM$, which we may write
 $$
  [\;(\,N^2\,+\,\Der_I^2\,)\,\theta\;]_{n+1}
  \;\;\equiv\;\;8\,n\,(\,\cl\,H_n\,)
 $$
 in light of the definition (\ref{cl}) of the closure map. The normalization
 constraint imposed on the K\"ahler normal potential implies on the other
 hand the congruence $[\,(N^2+\Der_I^2)\,\theta\,]_{n+1}\,\equiv\,0$ modulo
 $\S^{\geq(2,2)}T^*_pM$ for all $n+1\,\neq\,2$ so that $\cl\,H_n\,=\,0$ and
 a fortiori $H_n\,=\,0$ due to the exactness of the sequence (\ref{ses}).
 \qed

 \pfill
 Although the preceeding proof of Theorem \ref{univ} takes the existence
 of K\"ahler normal coordinates stipulated in Theorem \ref{knc} for granted,
 it is easily rearranged to prove existence on the fly alongside the induction.
 In this setup we start with arbitrary anchored holomorphic coordinates
 $\knc^{(1)}_p:\,T_pM\longrightarrow M$ thought of as a first approximation
 to the K\"ahler normal coordinates we want to construct and use the flow
 of the real part $H_2$ of a quadratic holomorphic vector field to modify
 $\knc^{(1)}_p$ to the better approximation $\knc^{(2)}_p$ characterized
 by $H_2\,=\,0$. In turn we use the real part $H_3$ of a cubic holomorphic
 vector field to find an even better approximation $\knc^{(3)}_p$
 characterized by $H_2\,=\,0\,=\,H_3$ and so on. The advantage of this
 modified proof is that it applies verbatim to a wider class of complex
 affine manifolds:

 \begin{Definition}[Balanced Complex Affine Manifolds]
 \hfill\label{scam}\break
  A complex affine manifold is a smooth manifold $M$ endowed with a torsion
  free connection $\nabla$ on its tangent bundle and an almost complex
  structure $I$, which is parallel $\nabla I\,=\,0$ and thus integrable
  by the Theorem of Newlander--Nierenberg. A complex affine manifold $M$
  is balanced provided the curvature tensor $R$ of the connection $\nabla$
  is a $(1,1)$--form in the sense:
  $$
   R_{I\,\cdot\,,\,I\,\cdot\,}\;\;=\;\;R_{\,\cdot\,,\,\cdot\,}
   \qquad\Longleftrightarrow\qquad
   \delta R\;\;=\;\;0
  $$
 \end{Definition}

 \noindent
 In every point $p\,\in\,M$ of a balanced complex affine manifold $M$
 there thus exist unique anchored holomorphic coordinates $\knc_p:\,
 T_pM\longrightarrow M$ characterized by the congruence
 $$
  K
  \;\;:=\;\;
  \exp_p^{-1}\,\circ\,\knc_p
  \;\;\equiv\;\;
  \id
  \qquad\textrm{mod}\;\;\H^\perp
 $$
 imposed directly on the infinite order Taylor series of the difference
 element $K$. In consequence this Taylor series is given by the very same
 universal composition power series in $I,\,R,\,\nabla R,\,\ldots$ we
 found in the K\"ahler case. Although originally formulated in terms of
 the potential and so apparently depending on the metric structure of a
 K\"ahler manifold, the concept of K\"ahler normal coordinates turns out
 to arise from the underlying complex affine structure! The following
 corollary about totally geodesic submanifolds, whose proof is left to
 the reader, is a nice illustration of this dependence on the affine
 structure:

 \begin{Corollary}[Totally Geodesic Complex Submanifolds]
 \hfill\label{sub}\break
  For every totally geodesic complex submanifold $N\,\subset\,M$ of a
  K\"ahler manifold $M$ the restriction of the K\"ahler normal coordinates
  $\knc_p:\,T_pM\longrightarrow M$ in a point $p\,\in\,N$ to the vector
  subspace $T_pN\,\subset\,T_pM$ are the K\"ahler normal coordinates for
  the K\"ahler manifold $N$.
 \end{Corollary}

 \noindent
 Instead of using the difference element $K$ we may consider using its
 inverse in the proof of Theorem \ref{univ}. Lacking an analogue of equation
 (\ref{pott}) it is more difficult to relate the normalization constraint
 imposed on $K^{-1}\,=\,\knc_p^{-1}\,\circ\,\exp_p$ to the normalization
 constraint imposed on the potential $\theta$, however the choice of $K^{-1}$
 allows us to use the Gau\ss\ Lemma to simplify the analogue of the recursion
 formula (\ref{orec}). Mimicking the argument we find
 $$
  (\,N\,\otimes\,\id\,+\,\Der_I\,\otimes\,I\,)\;K^{-1}
  \;\;=\;\;
  2\,\{\,\D,\,\D^*\,\}\;K^{-1}
  \;\;=\;\;
  2\,\D^*\,\pr_{\Sigma^1}\Big(\;DK^{-1}\,\circ\,(\,\id\;-\;\Phi\,)\;\Big)
 $$
 since $DK^{-1}\,\circ\,\Phi$ still commutes with $I$ according to
 (\ref{inter}). Inserting the definitions of the projection $\pr_{\Sigma^1}$
 and the boundary operator $\D^*$ we find that the right hand side simplifies
 due to the Gau\ss\ Lemma (\ref{gss}) written in the form $(\id-\Phi(X))\,X
 \,=\,0$, moreover we may replace the cumbersome notation involving $DK^{-1}$
 by the derivation extension of the endomorphism $F\,:=\,(\id-\Phi(X))\,\circ
 \,I$ using the identity $(\,\Der_FK^{-1}\,)(X)\,=\,DK^{-1}(X)\,FX$:

 \begin{Remark}[Explicit Version of Recursion Formula]
 \hfill\label{recs}\break
  The Taylor series of the inverse difference element $K^{-1}$ satisfies the
  recursion formula:
  \begin{eqnarray*}
   \Big[\;(\,N\,\otimes\,\id\,+\,\Der_I\,\otimes\,I\,)\,K^{-1}\;\Big](\;X\;)
   &=&
   I\,\Big[\,\Der_{(\id-\Phi(X))\,I}K^{-1}\,\Big](\,X\,)
   \\
   &=&
   +\;\sum_{k\,\geq\,0}\frac{k+1}{(k+3)!}\,
   (\nabla^k_{X,\,\ldots,\,X}R)^{}_{X,\,IX}IX
   \;\;+\;\;O(\,R^2\,)
  \end{eqnarray*}
  where $O(R^2)$ signifies composition polynomials at least quadratic in
  $R,\,\nabla R,\,\nabla^2R,\,\ldots$.
 \end{Remark}

 \noindent
 The relatively simple recursion formula for $K^{-1}$ allows us to calculate
 the lowest order terms of the Taylor series for all relevant objects on a
 K\"ahler manifold like the difference element
 \begin{eqnarray*}
  K^{-1}X
  &=&
  X\;+\;\frac1{12}\,R_{X,\,IX}IX\;+\;\frac1{96}\,
  \Big(\,3\,(\nabla_XR)_{X,IX}IX\;-\;(\nabla_{IX}R)_{IX,X}X\;\Big)
  \\
  &&
  \hphantom{X\;}+\;\frac1{960}\,\Big(\,7\,(\nabla^2_{X,X}R)_{X,IX}IX\,+\,
  2\,(\nabla^2_{IX,X}R)_{X,IX}X\,+\,2\,(\nabla^2_{X,IX}R)_{X,IX}X
  \\
  &&
  \qquad\;-\,(\nabla^2_{IX,IX}R)_{X,IX}IX\,\Big)\,-\,\frac1{120}\,
  R_{X,R_{X,IX}X}IX\,-\,\frac1{96}\,R_{X,R_{X,IX}IX}X\,+\,O(\,X^6\,)
 \end{eqnarray*}
 which provides the expansion of $K$ by formal inversion. Equation (\ref{pott})
 can be used to expand
 \begin{eqnarray*}
  \theta(\,X\,)
  &=&
  g(\,X,\,X\,)\;-\;\frac18\,g(\,R_{X,IX}IX,X\,)
  \;-\;\frac1{24}\,g(\,(\nabla_XR)_{X,IX}IX,X\,)
  \\[3pt]
  &&
  \qquad-\;\frac1{576}\,\Big(\,5\,g(\,(\nabla^2_{X,X}R)_{X,IX}IX,X\,)
  \,-\,g(\,(\nabla^2_{IX,IX}R)_{X,IX}IX,X\,)\,\Big)\qquad\qquad
  \\[3pt]
  &&
  \qquad+\;\frac1{48}\;g(\,R_{X,IX}X,R_{X,IX}X\,)\;+\;O(\,X^7\,)
 \end{eqnarray*}
 the K\"ahler normal potential and Remark \ref{dfk} to expand the Riemannian
 distance function:
 \begin{eqnarray*}
  \dist^2_g(\,p,\,\knc^{}_pX\,)
  &=&
  g(\,X,X\,)
  \;-\;\frac16\,g(\,R_{X,IX}IX,X\,)
  \;-\;\frac1{16}\,g(\,(\nabla_XR)_{X,IX}IX,X\,)
  \\[3pt]
  &&
  \quad\;-\;\frac1{480}\,\Big(\,7\,g(\,(\nabla^2_{X,X}R)_{X,IX}IX,X\,)
  \,-\,g(\,(\nabla^2_{IX,IX}R)_{X,IX}IX,X\,)\,\Big)
  \\[3pt]
  &&
  \quad\;+\;\frac{23}{720}\;g(\,R_{X,IX}X,R_{X,IX}X\,)\;+\;O(\,X^7\,)
 \end{eqnarray*}
 With somewhat more effort one obtains the congruences (\ref{tcong}),
 (\ref{dcong}) given in the introduction.

 \pfill
 Before closing this section we want to discuss the analogue of the
 exponentially extended vector fields $Z\,\rightsquigarrow\,Z^{\exp}$
 defined in relation with the forward and backward parallel transport
 in Section \ref{affine}. Replacing the affine exponential map $\exp_p:\,
 T_pM\longrightarrow M$ used in this construction with K\"ahler normal
 coordinates $\knc_p:\,T_pM\longrightarrow M$ we can define for every
 tangent vector $Z\,\in\,T_pM$ the holomorphically extended vector field
 $Z^\knc$ on $T_pM$ by choosing a representative curve $\gamma:\,
 ]-\,\varepsilon,\,+\varepsilon\,[\longrightarrow M$ for $Z$ with
 associated parallel transport $\PT^\nabla_\gamma(t)$:
 $$
  Z^\knc
  \;\;:=\;\;
  \left.\frac d{dt}\right|_0
  \Big(\;T_pM\;\stackrel{\PT^\nabla_\gamma(t)}\longrightarrow\;
  T_{\gamma(t)}M\;\stackrel{\knc_{\gamma(t)}}\longrightarrow\;
  M\;\stackrel{(\knc_p)^{-1}}\longrightarrow\;T_pM\;\Big)
  \;\;\in\;\;
  \Gamma(\;T\,(T_pM)\;)
 $$
 By construction $Z^\knc$ is actually the real part of a holomorphic vector
 field on $T_pM$, because the parallel transport $\PT^\nabla_\gamma(t):\,
 T_pM\longrightarrow T_{\gamma(t)}M$ intertwines the complex structures
 $I_p$ and $I_{\gamma(t)}$ and thus can be thought of as a biholomorphism
 between $T_pM$ and $T_{\gamma(t)}M$ considered as complex manifolds, in
 consequence $Z^\knc$ is the derivative of a family of biholomorphisms.
 Expanding the composition defining this family of biholomorphisms
 in order to bring the properties of the difference element $K$ to bear
 we end up with the alternative formulation
 $$
  Z^\knc
  \;\;:=\;\;
  \left.\frac d{dt}\right|_0K^{-1}_p\;\circ\;
  \Big(\,\exp_p^{-1}\circ\exp_{\gamma(t)}\circ\PT^\nabla_\gamma(t)\,\Big)
  \;\circ\;\Big(\,\PT^\nabla_\gamma(\,t\,)^{-1}\circ K_{\gamma(t)}\circ
  \PT^\nabla_\gamma(t)\,\Big)
 $$
 in which the stipulated derivative in $t\,=\,0$ leads to the rather
 simple formula:
 \begin{equation}\label{kncx}
  Z^\knc(\;X\;)
  \;\;=\;\;
  K^{-1}_{*,\,KX}\;\Big(\;Z^{\exp}(\,KX\,)\;+\;(\,\nabla_ZK\,)\,X\;\Big)
 \end{equation}
 The vector field $Z^\knc$ thus reflects the covariant derivative of the
 difference element $K$ in the direction of $Z\,\in\,T_pM$ modified by adding
 the exponentially extended vector field $Z^{\exp}$. On the other hand
 $Z^\knc$ reflects the covariant derivative of the pull back $\theta_p
 \,:=\,\theta_p^\loc\circ\knc_p$:
 \begin{eqnarray*}
  (\,Z^\knc\theta_p\,)(\;X\;)
  &=&
  \left.\frac d{dt}\right|_0\theta_p\Big(\,\knc^{-1}_p
  (\,\knc_{\gamma(t)}\,\PT^\nabla_\gamma(t)X\,)\,\Big)
  \\
  &=&
  \left.\frac d{dt}\right|_0\Big(\,\theta^\loc_p\,-\,\theta^\loc_{\gamma(t)}
  \,\Big)\,\Big(\,\knc_{\gamma(t)}\,\PT^\nabla_\gamma(t)\,X\,\Big)
  \;+\;\left.\frac d{dt}\right|_0
  \theta_{\gamma(t)}\Big(\,\PT^\nabla_\gamma(t)\,X\,\Big)
  \\[3pt]
  &=&
  \textrm{real part of holomorphic function}\;+\;
  (\,\nabla_Z\theta\,)_p(\;X\;)
 \end{eqnarray*}
 Recall at this point that the difference of the two local potentials
 $\theta^\loc_p$ and $\theta^\loc_{\gamma(t)}$ is the real part of a
 holomorphic function and so then is the derivative $\left.\frac d{dt}
 \right|_0(\theta^\loc_p-\theta^\loc_{\gamma(t)})$. A little consideration
 reveals that the only possible candidate for this real part is a multiple
 of $Z^\sharp$ and hence:
 \begin{equation}\label{kx}
  Z^\knc\,\theta_p
  \;\;=\;\;
  2\,Z^\sharp\;+\;(\,\nabla_Z\theta\,)_p
 \end{equation}
 Decomposing this equation into its bihomogeneous components (\ref{sdec})
 it can be solved simultaneously for the vector field $Z^\knc$ on $T_pM$
 and the covariant derivative $(\nabla_Z\theta)_p$ of the K\"ahler normal
 potential. For this purpose let us decompose the normal potential into three
 parts
 $$
  \theta_p
  \;\;=\;\;
  2\,g_p\;+\;\theta^\free_p
  \;\;=\;\;
  2\,g_p\;+\;(\;\theta^\crit_p\;+\;\theta^\rest_p\;)
 $$
 where $\theta^\rest_p\,\in\,\overline\S^{\geq(3,3)}T^*_pM$ and the critical
 part $\theta^\crit_p$ comprises all bihomogeneous components of $\theta_p$,
 which are exactly quadratic in the holomorphic or antiholomorphic coordinates.
 The irritating notation for the quadratic polynomial $(2g_p)(X)\,=\,g_p(X,X)$
 is certainly a drawback of our convention relating symmetric forms with
 polynomials.

 \begin{Lemma}[Spencer Connection on K\"ahler Potential]
 \hfill\label{spcon}\break
  The directional derivative $Z\,\ins\,\theta^\crit$ of the critical part
  $\theta^\crit$ of the K\"ahler normal potential $\theta$ in the direction
  of a tangent vector $Z\,\in\,T_pM$ in a point $p\,\in\,M$ decomposes into
  the sum
  $$
   Z\,\ins\,\theta^\crit
   \;\;=\;\;
   \pr_{[\,1\,]}(\;Z\;\ins\;\theta^\crit\;)
   \;+\;\pr_{\geq(2,2)}(\;Z\;\ins\;\theta^\crit\;)
  $$
  where $\pr_{\geq(2,2)}(\,Z\,\ins\,\theta^\crit\,)$ is a power series at
  least quadratic in both the holomorphic and antiholomorphic coordinates
  and $\pr_{[\,1\,]}(\,Z\,\ins\,\theta^\crit\,)\,\in\,\overline\S_{[\,1\,]}
  T^*_pM$ is the closure of the real part of a holomorphic vector field
  essentially equal to the holomorphically extended vector field:
  $$
   Z^\knc
   \;\;=\;\;
   Z\;-\;\frac12\,\cl^{-1}
   \Big(\;\pr_{[\,1\,]}(\,Z\,\ins\,\theta^\crit\,)\;\Big)
  $$
  Moreover the covariant derivative of the K\"ahler normal potential is
  given by a projection
  $$
   \nabla_Z\theta
   \;\;=\;\;
   \pr_{\geq(2,2)}\,\Big(\;Z\;\ins\;\theta^\free\;\Big)
   \;-\;\frac12\,\cl^{-1}\,
   \Big(\;\pr_{[\,1\,]}(\;Z\;\ins\;\theta^\crit\;)\;\Big)\;\theta^\free
  $$
  of its formal derivative $Z\,\ins\,\theta^\free$ and a term depending
  bilinearly on $\theta^\crit$ and $\theta^\free$.
 \end{Lemma}

 \noindent
 In order to justify Lemma \ref{spcon} let us solve equation (\ref{kx})
 with respect the holomorphically extended vector field $Z^\knc$ associated
 to a tangent vector $Z\,\in\,T_p$. Equation (\ref{kncx}) together with the
 explicit expansions of $Z^{\exp}$ and $K$ immediately imply that $Z^\knc(X)
 \,=\,Z\,+\,O(X^2)$ has no linear term and thus can be reconstructed without
 ambiguity from its closure, which implicitly appears on the left hand side
 of the equation considered. In fact the decomposition of power series into
 their bihomogeneous components (\ref{sdec}) is parallel and hence the right
 hand side of equation (\ref{kx}) agrees with $\nabla_Z\theta\,=\,\nabla_Z
 \theta^\free\,\in\,\Gamma(\,\overline\S^{\geq(2,2)}T^*M\,)$ up to the linear
 term $2Z^\sharp$. Decomposing the left hand side $Z^\knc\theta$ analoguously
 into a linear term, a term in $\Gamma(\,\overline\S^{\geq(2,2)}T^*M\,)$ and
 a necessarily vanishing remainder term we find
 \begin{eqnarray*}
  (\,Z^\knc\theta\,)(\,X\,)
  &=&
  2\,g(\,Z,X\,)\;+\;\Big(\;2\,g(\,Z^\rest(X),\,X\,)
  \;+\;\pr_{[\,1\,]}(\,Z\,\ins\,\theta^\crit\,)(\,X\,)\;\Big)
  \\
  &&
  \hphantom{2\,g(\,Z,X\,)}\;+\;
  \Big(\;\pr_\rest(\,Z\,\ins\,\theta^\crit\,)(\,X\,)\;+\;(\,Z\,\ins\,
  \theta^\rest\,)(\,X\,)\;+\;(\,Z^\rest\,\theta^\free\,)\;\Big)
 \end{eqnarray*}
 where $Z^\rest\,=\,-\frac12\,\cl^{-1}\,\pr_{[\,1\,]}(\,Z\,\ins\,
 \theta^\crit\,)$ denotes all non--constant components of the holomorphically
 extended vector field $Z^\knc\,=\,Z\,+\,Z^\rest$. Taking on the other hand a
 closer look at Definition \ref{hhc} we observe that for all $k\,\geq\,4$ the
 bihomogeneous component of the higher holomorphic sectional curvature tensor
 $S_k\,\in\,\S^kT^*_pM$ of bidegree $(k-2,2)$ equals:
 $$
  S_{k-2,\,2}(\;X\;)
  \;\;=\;\;
  \frac1{(k-4)!}\;g(\;(\nabla_{\pr^{1,0}X,\,\pr^{1,0}X,\,\ldots,\,
  \pr^{1,0}X}R)_{X,\,IX}IX,\;X\;)
 $$ 
 In turn the congruence (\ref{tcong}) implies the following explicit
 formula for the critical part
 \begin{eqnarray*}
  \theta^\crit(\,X\,)
  &=&
  -\;\frac18\,g(\;R_{X,\,IX}IX,\;X\;)\;-\;\sum_{k\,>\,4}
  \frac1{4\,(k-2)\,(k-3)}\;\Big(\;S_{k-2,\,2}(\,X\,)
  \;+\;S_{2,\,k-2}(\,X\,)\;\Big)
  \\
  &=&
  -\;\frac18\,g(\,R_{X,\,IX}IX,\,X\,)\;-\;\sum_{k\,>\,4}
  \frac1{2\,(k-2)!}\;\Re\;g\Big(\,(\nabla_{\pr^{1,0}X,\,\ldots,\,
  \pr^{1,0}X}R)_{X,\,IX}IX,\,X\,\Big)
 \end{eqnarray*}
 of the K\"ahler normal potential, from which the closure $\pr_{[\,1\,]}
 (\,Z\,\ins\,\theta^\crit\,)$ of the holomorphically extended vector field
 $-2\,Z^\knc$ associated to a tangent vector $Z\,\in\,T_pM$ is easily
 calculated:

 \begin{Corollary}[Holomorphically Extended Vector Fields]
 \hfill\label{hex}\break
  The Taylor series of the holomorphically extended vector field $Z^\knc$ on
  $T_pM$ associated to $Z\,\in\,T_pM$ depends linearly on the curvature tensor
  $R_p$ and its iterated covariant derivatives:
  $$
   Z^\knc(\;X\;)
   \;\;=\;\;
   Z\;+\;\sum_{k\,\geq\,4}\frac1{2^{k-3}\,(k-2)!}\,\Re\;\Big(
   \;(\nabla^{k-4}_{X-iIX,\,\ldots,\,X-iIX}R)_{Z+iIZ,\;IX}(IX+iX)\;\Big)
  $$
 \end{Corollary}
\section{Hermitean Symmetric Spaces}
\label{spaces}
 In differential geometry locally symmetric spaces form an important
 laboratory to test new concepts and ideas for plausibility, because
 the characteristic vanishing $\nabla R\,=\,0$ of the covariant derivative
 of the curvature tensor $R$ on a symmetric space sets up a Lie theoretic
 framework for doing calculations, which are unfeasible or even impossible
 to do on general affine manifolds. In this section we use this framework
 to derive formulas for the difference elements, the normal potential and
 the holomorphically extended vector fields on hermitean symmetric spaces,
 although all results except Corollary \ref{kps} are valid for all symmetric
 among the balanced complex affine spaces. A standard reference for the
 Lie theoretic framework, albeit for Riemannian symmetric spaces only,
 is \cite{s}.

 \pfill
 The starting point of our discussion of hermitean symmetric spaces is the
 affine Killing equation for a vector field on an affine manifold $M$
 endowed with a torsion free connection $\nabla$ on its tangent bundle.
 In terms of the Lie derivative of connections this equation reads
 $$
  0
  \;\;\stackrel?=\;\;
  (\mathfrak{Lie}_X\nabla)_YZ
  \;\;:=\;\;
  [\,X,\,\nabla_YZ\,]\,-\,\nabla_{[X,Y]}Z\,-\,\nabla_X[\,Y,\,Z\,]
  \;\;=\;\;
  R_{X,\,Y}Z\,+\,\nabla^2_{Y,\,Z}X
 $$
 its solutions $X\,\in\,\Gamma(\,TM\,)$ are called affine Killing fields. The
 extended affine Killing field $X^\ext\,:=\,\X\oplus X$ with $\X\,:=\,\nabla X$
 associated to a solution satisfies the extended equation
 \begin{equation}\label{ake}
  \nabla_YX\;\;=\;\;\X\,Y\qquad\qquad
  \nabla_Y\X\;\;=\;\;-R_{X,\,Y}
 \end{equation}
 the extended affine Killing field is thus a parallel section
 $\nabla^\Killing_Y(\,\X\oplus X\,)\,=\,0$ of the vector bundle
 $\End\,TM\,\oplus\,TM$ with respect to the Killing connection defined by:
 \begin{equation}\label{kc}
  \nabla^{\Killing}_Y(\,\X\,\oplus\,X\,)
  \;\;:=\;\;
  (\;\nabla_Y\X\;+\;R_{X,\,Y}\;)\,\oplus\,(\;\nabla_YX\,-\,\X\,Y\;)
 \end{equation}
 A direct consequence of this construction is that the Lie subalgebra of
 affine Killing vector fields $\mathfrak{aff}(\,M,\,\nabla\,)\,\subset\,
 \Gamma(\,TM\,)$ is a vector space of dimension $\dim\,\mathfrak{aff}
 (\,M,\,\nabla\,)\,\leq\,m^2\,+\,m$ on every connected manifold $M$ of
 dimension $m$ with equality only on affine spaces. Straightforward, but
 slightly tedious calculations result in an explicit formula for the
 curvature
 $$
  R^{\Killing}_{Y,\,Z}(\,\X\,\oplus\,X\,)
  \;\;=\;\;
  \Big(\;(\nabla_XR)_{Y,\,Z}\,-\,(\X\star R)_{Y,\,Z}\;\Big)\;\oplus\;0
  \;\;\stackrel!=\;\;
  (\mathfrak{Lie}_XR)_{Y,\,Z}\;\oplus\;0
 $$
 of the Killing connection, the relation with the Lie derivative of the
 curvature tensor $R$ in the second equality requires $\X\,=\,\nabla X$.
 Due to torsion freeness the Lie bracket of two affine Killing fields
 $X,\,Y\,\in\,\mathfrak{aff}(\,M,\,\nabla\,)$ can be calculated via
 $[\,X,\,Y\,]\,=\,\mathfrak{Y}X\,-\,\X Y$ from their associated
 extended affine Killing fields $X^\ext,\,Y^\ext\,\in\,\Gamma(\,\End\,TM
 \oplus\,TM\,)$. Taking the covariant derivative and using (\ref{ake}) we
 find that the $C^\infty M$--bilinear ``algebraic'' bracket
 \begin{equation}\label{lie}
  \Big[\;\X\,\oplus\,X,\;\Y\,\oplus\,Y\;\Big]_\alg
  \;\;:=\;\;
  (\;[\,\Y,\,\X\,]\;-\;R_{Y,\,X}\;)\;\oplus\;(\;\Y\,X\;-\;\X\,Y\;)
 \end{equation}
 on the vector bundle $\End\,TM\,\oplus\,TM$ reduces to $[\,X^\ext,\,Y^\ext
 \,]_\alg\,=\,[\,X,\,Y\,]^\ext$ for extended affine Killing fields. In
 general this algebraic bracket does not satisfy the Jacobi identity and
 thus does not define a fiberwise Lie algebra structure on the vector bundle
 $\End\,TM\oplus TM$, the degree of failure of the Jacobi identity is however
 measured by the cyclic sum:
 \begin{eqnarray*}
  \lefteqn{[\;\X\,\oplus\,X,\;[\;\Y\,\oplus\,Y,\;
  \mathfrak{Z}\,\oplus\,Z\;]_\alg\,]_\alg\;+\;\textrm{cyclic permutations}}
  \qquad\qquad
  &&
  \\
  &=&
  \Big(\;(\,\X\,\star\,R\,)_{Y,\,Z}\;+\;(\,\Y\,\star\,R\,)_{Z,\,X}
  \;+\;(\,\mathfrak{Z}\,\star\,R\,)_{X,\,Y}\;\Big)\;\oplus\;0
 \end{eqnarray*}
 Last but not least we remark that the covariant derivative of
 $[\,,\,]_\alg$ under $\nabla^\Killing$ equals:
 $$
  (\nabla^\Killing_Z[\;,\,]_\alg)
  (\,\X\oplus X,\,\Y\oplus Y\,)
  \;\;=\;\;
  \Big(\,(\nabla_ZR)_{X,\,Y}+(\X\star R)_{Y,\,Z}-(\mathfrak{Y}\star R)_{X,\,Z}
  \,\Big)\,\oplus\,0
 $$
 Let us now specify the preceeding equations to a hermitean locally symmetric
 space, this is a K\"ahler manifold $M$ with covariantly parallel curvature
 tensor $\nabla R\,=\,0$. In this case the joint stabilizer of the parallel
 Riemannian metric $g$, complex structure $I$ and curvature tensor $R$
 defines a subbundle parallel with respect to the Levi--Civita connection
 $\nabla$
 \begin{equation}\label{z2}
  \g M
  \;\;:=\;\;
  \stab(\,g,\,I,\,R\,)\;\oplus\;TM
  \;\;\subset\;\;
  \End\,TM\;\oplus\;TM
 \end{equation}
 which is moreover parallel under the Killing connection $\nabla^\Killing$
 defined in equation (\ref{kc}), because the curvature tensor $R$ takes
 values in the joint stabilizer $\stab(g,I,R)$ of the parallel
 sections $g,\,I$ and $R$. On the subbundle $\mathfrak{g}M$ the algebraic
 bracket $[\,,\,]_\alg$ satisfies the Jacobi identity and is parallel with
 respect to the Killing connection, thus $\mathfrak{g}M$ becomes a bundle
 of $\mathbb{Z}_2$--graded Lie algebras endowed with the flat algebra
 connection $\nabla^\Killing$, moreover the odd subbundle is exactly the
 tangent bundle $TM$. In passing we remark that the Lie algebra structure
 (\ref{lie}) on $\g_pM$ encodes the Jacobi operators of the Riemannian
 metric with
 \begin{equation}\label{ad}
  (\,\ad^2X\,)\,A
  \;\;\widehat=\;\;
  [\;0\,\oplus\,X,\;[\;0\,\oplus\,X,\;0\,\oplus\,A\;]\;]
  \;\;=\;\;
  0\,\oplus\,R_{X,\,A}X
  \;\;\widehat=\;\;
  R_{X,\,A}X
 \end{equation}
 for all $X,\,A\,\in\,T_pM$. Under the ansatz $\Phi^{-1}(\,X\,)\,=\,
 \varphi^{-1}(\,\ad\,X\,)$ suggested by this identity the parallel transport
 equation of Lemma \ref{pte} becomes the ordinary differential equation
 $$
  x\frac d{dx}\,(\,x\frac d{dx}\,+\,1\,)\,\varphi^{-1}(\,x\,)
  \;\;=\;\;
  x^2\,\varphi^{-1}(\,x\,)
  \qquad\qquad
  \varphi^{-1}(\,0\,)\;\;=\;\;1
 $$
 for the unknown power series $\varphi^{-1}(\,x\,)\,\in\,\Q[\![\,x\,]\!]$
 with unique solution $\varphi^{-1}(\,x\,)\,=\,\frac{\sinh\,x}x$ or:
 \begin{equation}\label{pts}
  \Phi^{-1}(\;X\;)\;\;=\;\;\frac{\sinh\;\ad\,X}{\ad\,X}
  \qquad\qquad
  \Phi(\;X\;)\;\;=\;\;\frac{\ad\,X}{\sinh\;\ad\,X}
 \end{equation}
 For hermitean locally symmetric spaces we will use in addition the following
 specific identity:

 \begin{Lemma}[Fundamental Commutator Identity]
 \hfill\label{fcc}\break
  The Jacobi operator $(\ad^2X):\,T_pM\longrightarrow T_pM,\,A\longmapsto
  R_{X,\,A}X,$ associated to a tangent vector $X\,\in\,T_pM$ commutes on
  every hermitean symmetric space $M$ with the Jacobi operator $\ad^2IX$
  associated to the tangent vector $IX\,\in\,T_pM$. For general K\"ahler
  manifolds the commutator of these Jacobi operators is a linear
  combination of iterated covariant derivatives of $R$:
  \begin{eqnarray*}
   3\;[\;R_{X,\,\cdot}X,\;R_{IX,\,\cdot}IX\;]\;A
   &=&
   {\textstyle\frac12}\,
   (\nabla^2_{X,\,IX}R-\nabla^2_{IX,\,X}R)_{X,\,IX}A
   \;+\;(\nabla^2_{X,\,IX}R-\nabla^2_{IX,\,X}R)_{A,\,X}IX
   \\
   &&
   \,+\;(\nabla^2_{A,\,X}R-\nabla^2_{X,\,A}R)_{X,\,IX}IX
   \;+\;(\nabla^2_{A,\,IX}R-\nabla^2_{IX,\,A}R)_{X,\,IX}X
  \end{eqnarray*}
 \end{Lemma}

 \proof
 Although a Lie theoretic argument is significantly shorter we prefer to
 prove the general formula for arbitrary K\"ahler manifolds. It is rather
 difficult to believe though that such a complicated argument could be made
 up without knowledge of the Lie theoretic background. Using the first Bianchi
 identity three times we obtain the identity
 \begin{eqnarray*}
  \lefteqn{+\;R_{R_{A,\,U}V,\,W}Z\;-\;R_{R_{A,\,V}U,\,W}Z
   \;-\;R_{R_{A,\,W}U,\,V}Z\;+\;R_{R_{A,\,W}V,\,U}Z}\qquad
  &&
  \\
  \lefteqn{-\;R_{R_{A,\,Z}U,\,V}W\;+\;R_{R_{A,\,Z}V,\,U}W
   \;+\;R_{R_{A,\,Z}W,\,U}V\;-\;R_{R_{A,\,Z}W,\,V}U}\qquad
  &&
  \\
  &=&
  -\;R_{R_{U,\,V}A,\,W}Z\;+\;(\,R_{A,\,W}R\,)_{U,\,V}Z
  \;-\;[\,R_{A,\,W},\,R_{U,\,V}\,]\,Z
  \\
  &&
  +\;(\,R_{A,\,Z}R\,)_{U,\,V}W\;-\;[\,R_{A,\,Z},\,R_{U,\,V}\,]\,W
  \;-\;R_{U,\,V}R_{A,\,Z}W
  \\
  &=&
  +\;(\,R_{U,\,V}R\,)_{A,\,W}Z\;+\;R_{A,\,R_{U,\,V}W}Z
  \;+\;(\,R_{A,\,W}R\,)_{U,\,V}Z
  \\
  &&
  +\;(\,R_{A,\,Z}R\,)_{U,\,V}W\;-\;R_{A,\,Z}R_{U,\,V}W
  \\
  &=&
  -\;R_{R_{U,\,V}W,\,Z}A\;+\;(\,R_{U,\,V}R\,)_{A,\,W}Z
  \;+\;(\,R_{A,\,W}R\,)_{U,\,V}Z\;+\;(\,R_{A,\,Z}R\,)_{U,\,V}W
 \end{eqnarray*}
 for arbitrary tangent vectors $A$ and $U,\,V,\,W,\,Z$. Expressing
 the action of the curvature on curvature $(R_{U,V}R)_{A,W}Z\,=\,
 (\,\nabla^2_{U,V}R\,-\,\nabla^2_{V,U}R\,)_{A,W}Z$ by skew symmetrized
 iterated covariant derivatives and specifying $U\,=\,X\,=\,W$ and
 $V\,=\,IX\,=\,Z$ we conclude:
 \begin{eqnarray*}
  3\,[\;R_{X,\,\cdot}X,\;R_{IX,\,\cdot}IX\;]\,A
  &=&
  -\;R_{R_{X,\,IX}X,\,IX}A
  \;+\;(\nabla^2_{X,\,IX}R-\nabla^2_{IX,\,X}R)_{A,\,X}IX
  \\
  &&
  +\;(\nabla^2_{A,\,X}R-\nabla^2_{X,\,A}R)_{X,\,IX}IX
  \;+\;(\nabla^2_{A,\,IX}R-\nabla^2_{IX,\,A}R)_{X,\,IX}X
 \end{eqnarray*}
 In light of Remark \ref{scc} the latter identity implies the
 formula in question.
 \qed

 \pfill
 The most important consequence of Lemma \ref{fcc} is that we may evaluate
 a doubly even power series in two variables in the commuting Jacobi operators
 $\ad^2X$ and $\ad^2IX$ associated to a vector $X\,\in\,T_pM$ tangent to a
 hermitean symmetric space $M$. In particular the technical Lemma \ref{dde}
 is formulated in terms of the endomorphisms $F(\,\ad\,X,\,\ad\,IX\,)\,\in\,
 \End\,T_pM$ associated to tangent vectors $X\,\in\,T_pM$ and a power series
 $F\,\in\,\Q[\![\,x,\,\overline x\,]\!]$, which is doubly even in the sense
 $F(-x,\overline x)\,=\,F(x,\overline x)\,=\,F(x,-\overline x)$. More
 specifically we are interested in doubly even power series $F^\ext\,
 \in\,\Q[\![\,x,\,\overline x\,]\!]$ arising from even power series
 $F\,\in\,\Q[\![\,x\,]\!]$ via:
 \begin{equation}\label{ext}
  F^\ext(\;x,\;\overline x\;)
  \;\;:=\;\;
  \frac1{2\,x}\Big(\;(\,x\,+\,\overline x\,)\,F(\,x\,+\,\overline x\,)\;+\;
  (\,x\,-\,\overline x\,)\,F(\,x\,-\,\overline x\,)\;\Big)  
 \end{equation}
 With $F$ being even the expression $(x+\overline x)\,F(x+\overline x)
 \,+\,(x-\overline x)\,F(x-\overline x)$ vanishes at $x\,=\,0$ so that
 $F^\ext\,\in\,\Q[\![\,x,\,\overline x\,]\!]$ is indeed a well--defined
 power series with doubly even expansion
 \begin{equation}\label{fexp}
  F^\ext(\,x,\,\overline x\,)
  \;\;:=\;\;
  \sum_{k\,\geq\,0}\sum_{\mu\,=\,0}^k
  {2k\,+\,1\choose 2\mu}\;F_{2k}\;(\,x^2\,)^{k\,-\,\mu}(\,\overline x^2\,)^\mu
 \end{equation}
 in terms of the coefficients $(\,F_{2k}\,)_{k\geq0}$ of $F$. Usually it
 is simpler to calculate $F^\ext$ directly from its definition (\ref{ext}),
 for the even power series $\frac{\tanh\frac x2}{\frac x2}\,\in\,
 \Q[\![\,x\,]\!]$ for example we find:
 \begin{eqnarray*}
  \Big(\;\frac{\tanh\,\frac x2}{\frac x2}\;\Big)^\ext(\;x,\,\overline x\;)
  &=&
  \frac1x\;\Big(\;\frac{e^{x+\overline x}\,-\,1}{e^{x+\overline x}\,+\,1}
  \;+\;\frac{e^{x-\overline x}\,-\,1}{e^{x-\overline x}\,+\,1}\;\Big)
  \\
  &=&
  \frac2x\;\frac{(\,e^{2x}\,-\,1\,)\,e^{\overline x}}
  {(\,e^{x+\overline x}\,+\,1\,)\,(\,e^x\,+\,e^{\overline x}\,)}
  \;\;=\;\;
  4\;\frac{\sinh\,x}x\;\frac{e^x\;e^{\overline x}}
  {(\,e^{x+\overline x}\,+\,1\,)\,(\,e^x\,+\,e^{\overline x}\,)}
 \end{eqnarray*}
 In passing we remark that the power series $F^\ext\,\in\,\Q[\![\,x,\,
 \overline x\,]\!]$ satisfies the congruences
 \begin{equation}\label{cong}
  F^\ext(x,\overline x)
  \;\;\equiv\;\;
  F(x)
  \quad\textrm{mod}\;\;(\,\overline x^2\,)
  \qquad
  F^\ext(x,\overline x)
  \;\;\equiv\;\;
  \Big(\,\overline x\,\frac d{d\overline x}\,+\,1\,\Big)\,F(\overline x)
  \quad\textrm{mod}\;\;(\,x^2\,)
 \end{equation}
 modulo the ideals generated respectively by $\overline x^2$ and $x^2$. Both
 congruences can be read off from expansion (\ref{fexp}) or are readily
 derived from definition (\ref{ext}) using L'H\^opital's Rule.

 \begin{Lemma}[Technical Lemma]
 \hfill\label{dde}\break
  Consider an even formal power series $F\,\in\,\Q[\![\,x\,]\!]$ in one
  variable and its associated doubly even extension $F^\ext\,\in\,\Q
  [\![\,x,\,\overline x\,]\!]$ defined in equation (\ref{ext}). For every
  hermitean locally symmetric space $M$ the evaluation of $F^\ext$ at the
  commuting Jacobi operators $\ad^2X$ and $\ad^2IX$ associated to
  $X\,\in\,T_pM$ results in an endomorphism $F^\ext(\,\ad\,X,\,\ad\,IX\,)
  \,\in\,\End\,T_pM$ of the tangent space in $p$, which makes prominent
  appearance in the directional derivative:
  $$
   D\Big[\,F(\ad\,IX)X\,\Big]\,A
   \;\;:=\;\;
   \left.\frac d{dt}\right|_0F\,\Big(\,\ad\,I(X+tA)\,\Big)\,(X+tA)
   \;\;=\;\;
   F^\ext(\,\ad\,X,\,\ad\,IX\,)\,A
  $$
  Moreover the odd endomorphism $\ad\,[\,F(\ad\,IX)X\,]$ of the
  $\mathbb{Z}_2$--graded Lie algebra $\g_pM$ defined in equation
  (\ref{z2}) as $\stab_p(g,I,R)\,\oplus\,T_pM$ can be written
  in terms of $F^\ext(\,\ad\,X,\,\ad\,IX\,)$:
  \begin{eqnarray*}
   \left.\ad\,\Big(\;F(\,\ad\,IX\,)\,X\;\Big)\right|_{T_pM}\qquad\,
   &=&
   (\,\ad\,X\,)\;\circ\;F^\ext(\;\ad\,X,\;\ad\,IX\;)
   \\
   \left.\ad\,\Big(\;F(\,\ad\,IX\,)\,X\;\Big)\right|_{\stab_p(g,I,R)}
   &=&
   F^\ext(\;\ad\,X,\;\ad\,IX\;)\;\circ\;(\,\ad\,X\,)
  \end{eqnarray*}
  In consequence the square of this endomorphism restricted to the tangent
  space $T_pM$ reads:
  $$
   \ad^2\,\Big(\;F(\,\ad\,IX\,)\,X\;\Big)
   \;\;=\;\;
   (\,\ad^2\,X\,)\;F^\ext(\;\ad\,X,\;\ad\,IX\;)^2
  $$
 \end{Lemma}

 \proof
 Recall that the even subalgebra $\stab_p(g,I,R)$ of the
 $\mathbb{Z}_2$--graded Lie algebra $\g_pM$ consists of endomorphisms
 of $T_pM$ commuting with $I$. In consequence the symmetries of the
 curvature tensor $R$ of K\"ahler type of $M$ imply $[\,IX,\,I\hat X\,]
 \,=\,R_{IX,\,I\hat X}\,=\,[\,X,\,\hat X\,]$ and so
 $$
  [\,IX,\,[\,I\hat X,\,\Y\,]\,]
  \;\;=\;\;
  [\,IX,\,I\,[\,\hat X,\,\Y\,]\,]
  \;\;=\;\;
  [\,X,\,[\,\hat X,\,\Y\,]\,]
 $$
 for all $X,\,\hat X\,\in\,T_pM$ and all $\Y\,\in\,\stab_p(g,I,R)$.
 Generalizing this identity we conclude
 \begin{equation}\label{flip}
  \ldots[IX_{\mu+1},[IX_\mu,\ldots,[X_1,Y]\ldots]]\ldots
  \;\;=\;\;
  \ldots[X_{\mu+1},[X_\mu,\ldots,[X_1,Y]\ldots]]\ldots
 \end{equation}
 for all $X_1,\,\ldots,\,X_\mu,\,X_{\mu+1},\,\ldots$ in $T_pM$ provided
 $\mu$ is even and $Y\,\in\,T_pM$ or alternatively $\mu$ is odd and $Y\,\in\,
 \stab_p(g,I,R)$. A straightforward induction on $r$ or simply the binomial
 formula in the universal envelopping algebra generalizes the Jacobi identity
 for $\g_pM$ to:
 \begin{equation}\label{gig}
  \ad\,\Big(\;(\,\ad^r Y\,)\,Z\;\Big)
  \;\;=\;\;
  \sum_{\mu\,=\,0}^r(-1)^\mu\;{r\choose\mu}\;(\,\ad^{r-\mu}Y\,)
  \;(\,\ad\,Z\,)\;(\,\ad^\mu Y\,)
 \end{equation}
 Using this identity together with $\sum_{r\,=\,\mu}^{2k-1}{r\choose
 \mu}\,=\,{2k\choose\mu+1}$ we calculate for all $k\,\in\,\N_0$:
 \begin{eqnarray*}
  \lefteqn{D\,\Big(\,(\ad^{2k}IX)\,X\,\Big)\,A}
  \qquad
  &&
  \\
  &=&
  \left.\frac d{dt}\right|_0\Big(\,\ad^{2k}\,I(X+tA)\,\Big)\,(X+tA)
  \\
  &=&
  (\,\ad^{2k}IX\,)\,A\;-\;
  \sum_{r\,=\,0}^{2k-1}(\,\ad^{2k-r-1}IX\,)\;[\;(\,\ad^rIX\,)\,X,\;IA\;]
  \\
  &=&
  (\,\ad^{2k}IX\,)\,A\;-\;
  \sum_{\mu\,=\,0}^{2k-1}(-1)^\mu\,{2k\choose\mu+1}\,
  (\,\ad^{2k-\mu-1}IX\,)\,(\,\ad\,X\,)\,(\,\ad^\mu IX\,)\;IA
 \end{eqnarray*}
 In order to evaluate this sum we treat the summands differently depending
 on the parity of the index $\mu$: Identity (\ref{flip}) allows us to remove
 all $I$'s to the right of $\ad\,X$ without changing the result provided the
 index $\mu$ is odd, in the opposite case with even index $\mu$ we employ
 identity (\ref{flip}) to replace $(\ad\,X)(\ad\,IX)$ by $-(\ad\,IX)(\ad\,X)$
 and then remove all $I$'s even further to the right without doing any harm.
 The net result of all these modifications reads
 \begin{eqnarray*}
  \lefteqn{D\,\Big(\,(\ad^{2k}IX)\,X\,\Big)\,A}
  \qquad
  &&
  \\
  &=&
  (\,\ad^{2k}IX\,)\,A\;+\;
  \sum_{{\scriptstyle\mu\,=\,0}\atop{\scriptstyle\mu\;\mathrm{odd}}}^{2k}
  {2k\choose\mu+1}\,(\,\ad^2IX\,)^{k-\frac{\mu+1}2}\;
  (\,\ad^2X\,)^{\frac{\mu+1}2}\;A
  \\
  &&
  \hphantom{(\,\ad^{2k}IX\,)\,A}\;+\;
  \sum_{{\scriptstyle\mu\,=\,0}\atop{\scriptstyle\mu\;\mathrm{even}}}^{2k}
  {2k\choose\mu+1}\,(\,\ad^2IX\,)^{k-\frac\mu2}\;(\,\ad^2X\,)^{\frac\mu2}\;A
  \\
  &=&
  (\,\ad^{2k}IX\,)\,A\;+\;
  \sum_{\nu\,=\,0}^k\Big(\;{2k\choose 2\nu+1}\;+\;{2k\choose 2\nu}
  \;-\;\delta_{\nu\,=\,0}\;\Big)\,(\,\ad^2IX\,)^{k-\nu}\,(\,\ad^2X\,)^\nu\,A
  \\
  &=&
  \sum_{\nu\,=\,0}^k
  {2k+1\choose 2\nu+1}\,(\,\ad^2IX\,)^{k-\nu}\,(\,\ad^2X\,)^\nu\,A
  \;\;=\;\;
  (\;x^{2k}\;)^\ext(\;\ad\,X,\;\ad\,IX\,)\;A
 \end{eqnarray*}
 where $\nu\,=\,\frac\mu2$ or $\nu\,=\,\frac{\mu+1}2$ depending on the parity
 of $\mu$ and the Kronecker delta $\delta_{\nu\,=\,0}$ is needed to cancel
 the non existent summand with $\mu\,=\,-1$ after switching to $\nu\,=\,0$.
 With the directional derivative being linear we deduce the validity of
 the first statement of the lemma for all even power series $F\,\in\,
 \Q[\![\,x\,]\!]$ from its validity on the basis $(\,x^{2k}\,)_{k\,\in\,\N_0}$.

 \pfill
 The argument for the other two statements follows this line of reasoning
 very closely, the starting point for both statements is to use identity
 (\ref{gig}) to expand $\ad\,[\,(\ad^{2k}IX)\,X\,]$ into:
 $$
  \ad\,[\,(\,\ad^{2k}IX\,)\,X\,]
  \;\;=\;\;
  \sum_{\mu\,=\,0}^{2k}(-1)^\mu\;{2k\choose\mu}
  (\,\ad^{2k-\mu}IX\,)\;(\,\ad\,X\,)\;(\,\ad^\mu IX\,)
 $$
 In this situation making a case distinction on the parity of the index
 $\mu$ is not sufficient to employ identity (\ref{flip}) and so we have
 to restrict this endomorphism of $\g_pM$ either to the tangent space
 $T_pM$ or to $\stab_p(g,I,R)$. Restriction to $T_pM$ allows us to remove
 all $I$ to the left of $\ad\,X$ without changing the result provided $\mu$
 is even, for odd $\mu$ we have to replace $(\ad\,IX)(\ad\,X)$ by $-(\ad\,X)
 (\ad\,IX)$ first using (\ref{flip}) before removing all $I$ further to the
 left. In analogy we remove all $I$ to the right of $\ad\,X$ on $\stab_p
 (g,I,R)$ provided $\mu$ is even, otherwise we replace $(\ad\,X)(\ad\,IX)
 \,=\,-(\ad\,IX)(\ad\,X)$ first and proceed as before.
 \qed

 \begin{Theorem}[Difference Elements of Hermitean Symmetric Spaces]
 \hfill\label{main}\break
  For every hermitean locally symmetric space $M$ the difference element
  $K\,:=\,\exp_p^{-1}\,\circ\,\knc_p$ measuring the deviation between the
  Riemannian and K\"ahlerian normal coordinates reads:
  $$
   KX
   \;\;=\;\;
   \frac{\mathrm{artanh}\,(\,\frac12\,\ad\,IX\,)}{\frac12\,\ad\,IX}\;X
   \qquad\qquad
   K^{-1}X
   \;\;=\;\;
   \frac{\tanh(\,\frac12\,\ad\,IX\,)}{\frac12\,\ad\,IX}\;X
  $$
 \end{Theorem}

 \proof
 Before verifying the recursion formula (\ref{recs}) for $K^{-1}$ let
 us first use the technical Lemma \ref{dde} to show that $K$ and $K^{-1}$
 are actually composition inverses of each other. In order to simplify the
 exposition of this argument we consider the general case of two even power
 series $F,\,\hat F\,\in\,\Q[\![\,x\,]\!]$ parametrizing power series
 on $T_pM$ with values in $T_pM$ via
 \begin{equation}\label{swap}
  F(\,X\,)
  \;\;:=\;\;
  F(\;\ad\,IX\;)\;X
  \;\;\stackrel!=\;\;
  -\,I\,F(\;\ad\,X\;)\;IX
 \end{equation}
 the second equality rises from equation (\ref{flip}) in the form $[IX,[IX,
 Y]]\,=\,-I\,[X,[X,IY]]$ valid for all $X,\,Y\,\in\,T_pM$. The technical
 Lemma \ref{dde} implies for the composition:
 \begin{eqnarray*}
  F(\,\hat F(\,X\,)\,)
  &=&
  -\,I\,F\left(\,\sqrt{\ad^2[\,\hat F(\ad\,IX)X\,]}\,\right)\,
  I\,\hat F(\ad\,IX)X
  \\
  &=&
  -\,I\,F\left(\,\sqrt{(\ad^2X)\,\hat F^\ext(\ad\,X,\ad\,IX)^2}\,\right)\;
  \hat F(\ad\,X)\,IX
  \\[3pt]
  &=&
  -\,I\,F(\,\textrm{``}(\ad\,X)\,\hat F(\ad\,X)\textrm{''}\,)\;
  \hat F(\ad\,X)\,IX
  \;\;=\;\;
  -\,I\,(\,F\,\circ\,\hat F\,)(\,\ad\,X\,)\,IX
 \end{eqnarray*}
 where $(F\circ\hat F)(x)\,:=\,F(x\hat F(x))\,\hat F(x)$. For the series
 corresponding to $K$ and $K^{-1}$ we find
 $$
  \frac{\tanh\,\frac x2}{\frac x2}\;\circ\;
  \frac{\mathrm{artanh}\;\frac x2}{\frac x2}
  \;\;=\;\;
  \frac{\tanh\,\frac12\,(\,2\,\mathrm{artanh}\;\frac x2\,)}
  {\frac12\,(\,2\,\mathrm{artanh}\;\frac x2\,)}
  \;\frac{\mathrm{artanh}\;\frac x2}{\frac x2}
  \;\;=\;\;
  1
 $$
 and conclude $K^{-1}(KX)\,=\,X$, mutatis mutandis we obtain $K(K^{-1}X)\,=
 \,X$ as well. Recall now that the difference element $K^{-1}$ on a locally
 symmetric space is a composition polynomial in the curvature tensor alone
 as $\nabla R\,=\,0$ and all iterated covariant derivatives vanish, thus it
 has weight $0$ in the sense $\delta K^{-1}\,=\,0$. Rewriting the definition
 (\ref{wo}) of the weight operator
 $$
  \delta
  \;\;:=\;\;
  \id\,\otimes\,I\;-\;\Der_I\,\otimes\,\id
  \;\;\stackrel!=\;\;
  (\,\id\,\otimes\,I\,)\,(\,\Der_I\,\otimes\,I\;+\;\id\,)
 $$
 we conclude $(\Der_I\otimes I)\,K^{-1}\,=\,-K^{-1}$ so that the recursion
 formula of Remark \ref{recs} becomes:
 \begin{equation}\label{crit}
  \Big[\;(\,N\,-\,1\,)\,K^{-1}\;\Big](\,X\,)
  \;\;=\;\;
  I\,DK^{-1}(\,X\,)\,(\,\id\,-\,\Phi(X)\,)\,IX
 \end{equation}
 In order to verify that the stipulated power series $K^{-1}X\,=\,
 \frac{\tanh\,\frac12\,\ad\,IX}{\frac12\,\ad\,IX}\,X$ satisfies this
 formal differential equation characterizing the difference element
 let us consider the power series
 $$
  F(\,x\,)
  \;\;:=\;\;
  x\,\frac d{dx}\,\Big(\,\frac{\tanh\,\frac x2}{\frac x2}\,\Big)
  \;\;=\;\;
  \frac 1{\cosh^2\frac x2}\,-\,\frac{\tanh\,\frac x2}{\frac x2}
  \;\;=\;\;
  -\,\frac{\tanh\,\frac x2}{\frac x2}\;(\;1\,-\,\frac x{\sinh x}\;)
 $$
 where we have used the addition theorem $\frac1{\cosh^2\frac x2}
 \,=\,\frac{\tanh\frac x2}{(\cosh\frac x2)(\sinh\frac x2)}
 \,=\,\frac{\tanh\frac x2}{\frac12\sinh x}$ for the hyperbolic
 sine in the second equality. Conveniently using the swap identity
 (\ref{swap}) we conclude
 \begin{eqnarray*}
  [\,(\,N\,-\,1\,)\,K^{-1}\,](\,X\,)
  &=&
  F(\,\ad\,IX\,)\,X
  \;\;=\;\;
  -\,I\,F(\,\ad\,X\,)\,IX
  \\
  &=&
  I\,\frac{\tanh\,\frac12\ad\,X}{\frac 12\ad\,X}\,\Big(\,\id\,-\,
  \frac{\ad\,X}{\sinh\,\ad\,X}\,\Big)\,IX
 \end{eqnarray*}
 where the shift from $N-1$ to $x\frac d{dx}$ reflects the fact that $x
 \frac d{dx}$ counts all arguments but one of the power series $K^{-1}$. On
 the other hand the technical Lemma \ref{dde} and equation (\ref{pts}) for
 the forward parallel transport $\Phi(X)\,=\,\frac{\ad\,X}{\sinh\,\ad\,X}$
 on symmetric spaces allow us to write the right hand side of the formal
 differential equation (\ref{crit}) for $K^{-1}X\,=\,\frac{\tanh\,\frac12\ad
 \,IX}{\frac12\ad\,IX}X$ in the form
 \begin{eqnarray*}
  I\,DK^{-1}(\,X\,)\,(\,\id\,-\,\Phi(X)\,)\,IX
  &=&
  I\,\Big(\,\frac{\tanh\,\frac x2}{\frac x2}\;\Big)^\ext(\,\ad\,X,\,\ad\,IX\,)
  \Big(\,\id\,-\,\frac{\ad\,X}{\sinh\,\ad\,X}\,\Big)\,IX
  \\
  &=&
  I\,\frac{\tanh\,\frac12\ad\,X}{\frac 12\ad\,X}\,\Big(\,\id\,-\,
  \frac{\ad\,X}{\sinh\,\ad\,X}\,\Big)\,IX  
 \end{eqnarray*}
 where we may reduce $(\frac{\tanh\frac x2}{\frac x2})^\ext$ modulo the
 ideal $(\,\overline x^2\,)$ according to the congruences (\ref{cong}),
 because $\ad^2IX$ commutes with the power series $\frac{\ad\,X}{\sinh
 \ad\,X}$ in $\ad^2X$ and kills its argument $IX$. Comparing the results
 for the left and right hand sides we conclude that the stipulated power
 series $K^{-1}$ satisfies in fact the formal differential equation
 (\ref{crit}) uniquely characterizing the difference element $K^{-1}
 \,=\,\knc^{-1}_p\,\circ\,\exp_p$ as claimed.
 \qed

 \begin{Corollary}[K\"ahler Potential on Symmetric Spaces]
 \hfill\label{kps}\break
  On every hermitean locally symmetric space $M$ the normal potential takes
  the form:
  $$
   \theta_p(\;X\;)
   \;\;=\;\;
   g_p\,\Big(\;X,\;\frac{\log(\,\id\,-\,\frac14\,\ad^2IX\,)}
   {-\,\frac14\,\ad^2IX}\,X\;\Big)
  $$
 \end{Corollary}

 \proof
 Essentially the proof consists of a calculation of the expression
 $\Psi(X)X$ on hermitean locally symmetric spaces using again the key
 technical Lemma \ref{dde}. In a prelude to the proof we remark that the
 power series identity $2\,\mathrm{artanh}\,x\,=\,\log\frac{1+x}{1-x}$
 is a direct consequence of the differential equation $\mathrm{artanh}'x
 \,=\,\frac1{1-x^2}$. The second congruence in (\ref{cong}) thus becomes
 $$
  \Big(\,\frac{\mathrm{artanh}\,\frac x2}{\frac x2}\,\Big)^\ext
  (\;x,\;\overline x\;)
  \;\;\equiv\;\;
  \Big(\;\overline x\,\frac d{d\overline x}\;+\;1\,\Big)\;
  \frac{\mathrm{artanh}\,\frac{\overline x}2}{\frac{\overline x}2}
  \;\;=\;\;
  \frac1{1\,-\,(\frac{\overline x}2)^2}
 $$
 modulo the ideal generated by $x^2$. In consequence the formula (\ref{pts})
 for the backward parallel transport on symmetric spaces and the technical
 Lemma \ref{dde} allow us to expand the definition $\Psi^{-1}(X)\,=\,
 \Phi^{-1}(KX)\,K_{*,X}$ of the K\"ahler backward parallel transport
 to the effect that
 \begin{eqnarray*}
  \Psi^{-1}(\,X\,)\,X
  &=&
  \frac{\sinh\,\ad\,KX}{\ad\,KX}\;
  \Big(\,\frac{\mathrm{artanh}\,\frac x2}{\frac x2}\,\Big)^\ext
  (\,\ad\,X,\,\ad\,IX\,)\,X
  \\
  &=&
  \Big(\,\sum_{k\,\geq\,0}\frac1{(2k+1)!}\,\ad^{2k}KX\,\Big)
  \Big(\,1\,-\,{\textstyle\frac14}\,\ad^2IX\,\Big)^{-1}X
  \;\;=\;\;
  \Big(\,1\,-\,{\textstyle\frac14}\,\ad^2IX\,\Big)^{-1}X
 \end{eqnarray*}
 because $\ad^2X$ kills the eventual argument $X$ in the first line and
 so then does
 $$
  \ad^2KX
  \;\;=\;\;
  (\,\ad^2X\,)\;\circ\;
  \Big(\,\frac{\mathrm{artanh}\,\frac x2}{\frac x2}\,\Big)^\ext
  (\,\ad\,X,\,\ad\,IX\,)
 $$
 in the second. On the other hand we know that the K\"ahler normal
 potential $\theta$ on a hermitean locally symmetric space is necessarily
 of weight $0$ in the sense $\Der_I\theta\,=\,0$, because it is essentially
 a composition polynomial in the curvature tensor alone. With the Jacobi
 operator $\ad^2IX\,=\,R_{IX,\,\cdot}IX$ being a symmetric endomorphism
 we thus simplify equation (\ref{pott}) to:
 $$
  [\;N^2\,\theta\;](\,X\,)
  \;\;=\;\;
  4\,g(\;\Psi^{-1}(X)\,X,\;\Psi^{-1}(X)\,X\;)
  \;\;=\;\;
  4\,g(\;X,\;(\,1\,-\,{\textstyle\frac14}\,\ad^2IX\,)^{-2}\,X\;)
 $$
 The corollary now follows from the identity $(x\frac d{dx})^2\log(1-x^2)
 \,=\,-4x^2\,(1-x^2)^{-2}$, no additional shift is needed in this argument
 between the number operator $N$ and $x\frac d{dx}$.
 \qed

 \pfill
 Before closing this section we want to verify equation (\ref{kx})
 describing the Spencer connection for the K\"ahler normal potential
 on hermitean locally symmetric spaces. In light of the explicit formulas for
 the forward and backward parallel transport (\ref{pts}) it seems justified
 to make the ansatz $\Theta(X)\,=\,\theta(\ad\,X)$ for the power series
 $\Theta$ describing the exponentially extended vector fields on symmetric
 spaces with an unknown even power series $\theta\,\in\,\Q[\![\,x\,]\!]$ as
 parameter. Under this ansatz the equation (\ref{eext}) for $\Theta$ becomes
 the ordinary differential equation
 \begin{eqnarray*}
  x\frac d{dx}\,\Big(\,x\frac d{dx}\,-\,1\,\Big)\,\theta(\,x\,)
  &=&
  x\frac d{dx}\Big(\;(\,1\,-\,\frac x{\sinh x}\,)\,\frac x{\sinh x}\;\Big)
  \;+\;\Big(\;x\frac d{dx}\,\frac{\sinh x}x\;\Big)\,
  \Big(\;\frac x{\sinh x}\;\Big)^2
  \\
  &=&
  -\,x\frac d{dx}\Big(\;\frac x{\sinh x}\;\Big)^2
 \end{eqnarray*}
 with initial value $\theta(0)\,=\,1$, the simplication in the second line
 is due to the standard identity $(x\frac d{dx}f^{-1})\,f^2\,=\,-x\frac d{dx}f$
 valid for every invertible, commutative power series $f$. Hence
 $$
  \Big(\,x\frac d{dx}\,-\,1\,\Big)\,\theta(\,x\,)
  \;\;=\;\;
  -\,\Big(\,\frac x{\sinh x}\,\Big)^2\;+\;\textrm{integration constant}
  \qquad\qquad
  \theta(\,0\,)\;\;=\;\;1
 $$
 where the initial value $\theta(0)\,=\,1$ forces the integration constant
 to vanish. The unique even power series solving the latter ordinary
 differential equation equals $\theta(x)\,=\,\frac x{\tanh x}$ due to
 $\tanh'x\,=\,\frac1{\cosh^2x}$, hence exponentially extended vector
 fields on symmetric spaces read:
 \begin{equation}\label{expt}
  Z^{\exp}(\,X\,)
  \;\;=\;\;
  \Theta(\,X\,)\,Z
  \;\;=\;\;
  \frac{\ad\,X}{\tanh\;\ad\,X}\;Z
 \end{equation}
 Incidentally this formula is exactly the formula describing the unique Killing
 vector field with vanishing covariant derivative and value $Z\,\in\,T_pM$ in
 a point $p\,\in\,M$ in exponential coordinates \cite{s}. With the difference
 element $K$ being parallel $\nabla K\,=\,0$ on locally symmetric spaces
 formula (\ref{kncx}) tells us that the holomorphically extended vector fields
 are simply these transvection Killing vector fields written alternatively
 in K\"ahler normal coordinates:

 \begin{Corollary}[Holomorphically Extended Vector Fields]
 \hfill\label{symevf}\break
  In accordance with Corollary \ref{hex} the holomorphically extended vector
  field $Z^\knc$ asscoiated to a tangent vector $Z\,\in\,T_pM$ of a hermitean
  symmetric space equals the qudratic vector field
  $$
   Z^\knc(\;X\;)
   \;\;=\;\;
   Z\;-\;\frac14\,\Big(\;R_{Z,\,X}X\;-\;R_{Z,\,IX}IX\;\Big)
  $$
  on the tangent space $T_pM$. In particular all Killing vector fields on
  a hermitean symmetric space become at most quadratic vector fields
  when written in K\"ahler normal coordinates.
 \end{Corollary}

 \proof
 Recall first of all that the curvature tensor of a locally symmetric space
 is covariantly constant $\nabla R\,=\,0$ and so then is the difference
 element $K$. In consequence equation (\ref{kncx}) describing the Taylor
 series of holomorphically extended vector fields becomes
 \begin{equation}\label{wout}
  Z^\knc(\;X\;)
  \;\;=\;\;
  K^{-1}_{*,\,KX}\;Z^{\exp}(\,KX\,)
  \;\;=\;\;
  DK^{-1}(\;KX\;)\;\frac{\ad\,KX}{\tanh\;\ad\,KX}\;Z
 \end{equation}
 in light of the formula (\ref{expt}) for exponentially extended vector
 fields. Using the technical Lemma \ref{dde} and the swap identity (\ref{swap})
 we may rewrite $\ad^2(KX)$ and $\ad^2(IKX)$ in the form
 $$
  \begin{array}{ccccl}
   \ad^2(\;\,KX\,\;)
   &=&
   (\,\ad\,X\,)^2\;\Big(\,{\displaystyle\frac{\mathrm{artanh}\,\frac x2}
   {\frac x2}}\,\Big)^\ext(\;\ad\,X,\;\ad\,IX\;)^2
   &=&
   F^2(\;\ad\;X,\;\ad\,IX\;)
   \\[5pt]
   \ad^2(\,IKX\,)
   &=&
   -\,I\,\circ\,(\,\ad^2KX\,)\,\circ\,I
   &=&
   F^2(\;\ad\,IX,\;\ad\;X\;)
  \end{array}
 $$
 defining the power series $F\,\in\,\Q[\![\,x,\,\overline x\,]\!]$ in two
 variables arising in this way as:
 $$
  F(\;x,\;\overline x\;)
  \;\;:=\;\;
  x\,\Big(\,\frac{\mathrm{artanh}\,\frac x2}{\frac x2}\,\Big)^\ext
  (\;x,\;\overline x\;)
  \;\;=\;\;
  \frac12\;\Big(\;\log\,\frac{2+x+\overline x}{2-x-\overline x}\;+\;
  \log\,\frac{2+x-\overline x}{2-x+\overline x}\;\Big)
 $$
 The identity $2\,\mathrm{artanh}\,\frac x2\,=\,\log\frac{2+x}{2-x}$ can
 be used conveniently to arrive at the extension (\ref{ext}) needed in the
 second equality. In turn the description (\ref{wout}) of holomorphically
 extended vector fields expands with the help of the technical
 Lemma \ref{dde} into the power series
 \begin{equation}\label{hh}
  Z^\knc(\;X\;)
  \;\;=\;\;
  DK^{-1}(\;KX\;)\;\frac{\ad\,KX}{\tanh\;\ad\,KX}\;Z
  \;\;=\;\;
  H(\;\ad^2X,\;\ad^2IX\;)\,Z
 \end{equation}
 in $\ad^2X$ and $\ad^2IX$ alone, where the doubly even power series
 $H\,\in\,\Q[\![\,x,\,\overline x\,]\!]$ reads
 \begin{eqnarray*}
  H(\;x,\;\overline x\;)
  &:=&
  \Big(\,\frac{\tanh\,\frac x2}{\frac x2}\,\Big)^\ext
  (\;F(\,x,\,\overline x\,),\;F(\,\overline x,\,x\,)\;)\;
  \frac{F(\,x,\,\overline x\,)}{\tanh\;F(\,x,\,\overline x\,)}
  \\
  &=&
  4\,\cosh\,F(\,x,\,\overline x\,)\;
  \frac{e^{F(\,x,\,\overline x\,)\,+\,F(\,\overline x,\,x\,)}}
  {(\,e^{F(\,x,\,\overline x\,)\,+\,F(\,\overline x,\,x\,)}\,+\,1\,)
  \,(\,e^{F(\,x,\,\overline x\,)}\,+\,e^{F(\,\overline x,\,x\,)}\,)}
 \end{eqnarray*}
 the exemplary calculation of $(\frac{\tanh\,\frac x2}{\frac x2})^\ext$
 following definition (\ref{ext}) has been used here in the second line. In
 order to simplify the power series $F(x,\overline x)$ and $F(\overline x,x)$
 in this formula we change variables from $x,\,\overline x$ to the new
 variables $a\,:=\,\frac{x+\overline x}2$ and $b\,:=\,\frac{x-\overline x}2$
 and obtain:
 $$
  F(\,x,\,\overline x\,)
  \;\;=\;\;
  \frac12\,\Big(\log\frac{1+a}{1-a}\;+\;\log\frac{1+b}{1-b}\,\Big)
  \qquad\quad
  F(\,\overline x,\,x\,)
  \;\;=\;\;
  \frac12\,\Big(\log\frac{1+a}{1-a}\;+\;\log\frac{1-b}{1+b}\,\Big)
 $$
 Inserting these expressions into the power series $H$ simplifies it
 drastically, namely we find
 \begin{eqnarray*}
  H(\;x,\,\overline x\;)
  &=&
  2\;\frac{(\,\frac{1+a}{1-a}\,)\,
  \Big(\,\sqrt{\frac{1+a}{1-a}}\,\sqrt{\frac{1+b}{1-b}}
  \,+\,\sqrt{\frac{1-a}{1+a}}\,\sqrt{\frac{1-b}{1+b}}\,\Big)}
  {(\,\frac{1+a}{1-a}\,+\,1\,)\,
  \Big(\,\sqrt{\frac{1+a}{1-a}}\,\sqrt{\frac{1+b}{1-b}}
  \,+\,\sqrt{\frac{1+a}{1-a}}\,\sqrt{\frac{1-b}{1+b}}\,\Big)}
  \\
  &=&
  \hphantom{2\;}\frac{(1+a)\,\sqrt{\frac{1+b}{1-b}}
  \;+\;(1-a)\,\sqrt{\frac{1-b}{1+b}}}
  {\sqrt{\frac{1+b}{1-b}}\;+\;\sqrt{\frac{1-b}{1+b}}}
  \;\;=\;\;
  1\,+\,ab
  \;\;=\;\;
  1\,+\,\frac14\,(\,x^2\,-\,\overline x^2\,)
 \end{eqnarray*}
 by using the identity $\frac{\sqrt t}{\sqrt t+\sqrt{t^{-1}}}\,=\,
 \frac t{t+1}$ twice. Reinserting this result into our description
 of holomorphically extended vector fields (\ref{hh}) on hermitean
 symmetric spaces we get eventually
 $$
  Z^\knc(\;X\;)
  \;\;=\;\;
  H(\,\ad^2X,\,\ad^2IX\,)\,Z
  \;\;=\;\;
  Z\;+\;\frac14\,\Big(\;(\,\ad^2X\,)\,Z\;-\;(\,\ad^2IX\,)\,Z\;\Big)
 $$
 which converts via equation (\ref{ad}) or $(\ad^2X)Z\,=\,-R_{Z,X}X$ into the
 first statement of the Lemma. The second statement follows from the fact that
 the complementary Killing vector fields with vanishing value in $p$ are
 linear vector fields in K\"ahler normal coordinates.
 \qed

 \pfill
 Having established Corollary \ref{symevf} about holomorphically extended
 vector fields on hermitean locally symmetric spaces we can eventually verify
 the imporant equation (\ref{kx}). According to Corollary \ref{kps} the
 K\"ahler normal potential on hermitean locally symmetric spaces can
 be written in the form $\theta(X)\,=\,g(\,X,\,F(\ad\,IX)\,X\,)$ with
 the even power series:
 $$
  F(\,x\,)
  \;\;:=\;\;
  \frac{\log(\,1\,-\,(\frac x2)^2\,)}{-(\frac x2)^2}
  \;\;=\;\;
  \sum_{k\,\geq\,0}\frac1{k+1}\,(\,\frac x2\,)^{2k}
 $$
 Recall now that the Jacobi operators $\ad^2X\,=\,-R_{\cdot,\,X}X$ and
 $\ad^2IX$ are commuting symmetric endomorphisms and so then is $F^\ext
 (\,\ad\,X,\,\ad\,IX\,)$. Using the technical Lemma (\ref{dde}) we find
 \begin{eqnarray*}
  \left.\frac d{dt}\right|_0\,\theta(\;X\;+\;t\,Z\;)
  &=&
  g\Big(\;Z,\;
  \left[\;F(\,\ad\,IX\,)\;+\;F^\ext(\,\ad\,X,\,\ad\,IX\,)\;\right]\,X\;\Big)
  \\
  &=&
  2\,g(\;Z,\;(\,\id\,-\,{\textstyle\frac14}\,\ad^2IX\,)^{-1}X\;)
 \end{eqnarray*}
 for every fixed direction $Z\,\in\,T_pM$ in light of the congruence of
 power series
 $$
  F(\,\overline x\,)\;+\;F^\ext(\,x,\,\overline x\,)
  \;\;\equiv\;\;
  (\;\overline x\,\frac d{d\overline x}\;+\;2\;)\;F(\,\overline x\,)
  \;\;=\;\;
  2\,\Big(\;1\;-\;(\,\frac{\overline x}2\,)^2\;\Big)^{-1}
  \qquad
  \textrm{mod}\;\;(\,x^2\,)
 $$
 in the variables $x,\,\overline x$ modulo the ideal generated by $x^2$, the
 equality in this congruence is most easily deduced from the power series
 expansion of $F$. In consequence we obtain:
 \begin{eqnarray*}
  (\,Z^\knc\theta\,)(\;X\;)
  &=&
  2\,g(\;Z\;+\;{\textstyle\frac14}\,(\;\ad^2X\;-\;\ad^2IX\;)\,Z,
  \;(\,\id\,-\,{\textstyle\frac14}\,\ad^2IX\,)^{-1}X\;)
  \\[5pt]
  &=&
  2\,g(\;Z,\;(\,\id\,-\,{\textstyle\frac14}\,\ad^2IX\,)\;
  (\,\id\,-\,{\textstyle\frac14}\,\ad^2IX\,)^{-1}X\;)
  \;\;=\;\;
  2\,g(\;Z,\;X\;)
 \end{eqnarray*}
\appendix
\section{K\"ahler Normal Coordinates in Examples}
\label{examples}
 In order to illustrate the concept of K\"ahler normal coordinates we
 want to desribe these coordinates and the associated normal potentials
 for a couple of examples in this appendix, namely the four infinite series
 of compact hermitean symmetric spaces: The complex Grassmannians, which
 include the complex projective spaces as a special case, the real
 Grassmannians of oriented planes and the real and quaternionic twistor
 spaces. A classical reference for formulas pertaining to Riemannian
 symmetric spaces like the relation between the curvature tensor and
 the Lie algebra structure on the Lie algebra of Killing vector fields
 exploited in Section \ref{spaces} is certainly \cite{s}. From among
 these formulas we will use in particular the description of the Riemannian
 exponential map in terms of matrix exponentials.
 
 \pfill
 Probably the most prominent examples of K\"ahler manifolds besides flat
 $\C^n$ are the complex projective spaces, which we will discuss in the
 context of complex Grassmann manifolds in the first part of this appendix.
 The K\"ahler metric of choice on the complex Grassmannian $\Gr_kV$ of
 $k$--dimensional subspaces of a complex vector space $V$ is a generalization
 of the Fubini--Study metric on the set $\mathbb{P}\,V\,=\,\Gr_1V$ of
 lines defined in terms of a positive definite hermitean form $h:\,
 \overline V\times V\longrightarrow\C$ on $V$ and the associated
 identification of the vector space
 $$
  \Hom(\;P,\;P^\perp\;)\;\stackrel\cong\longrightarrow\;T_P\Gr_kV,
  \qquad X\;\longmapsto\;\left.\frac d{dt}\right|_0
  \im\,(\;\id\,+\,t\,X:\quad P\,\longrightarrow\,V\;)
 $$
 tangent to $\Gr_kV$ in a point $P\,\in\,\Gr_kV$ with $h$--orthogonal
 complement $P^\perp\,\subset\,V$. Under this identification the complex
 structure $I$ becomes the multiplication by $i$ in $\Hom(P,P^\perp)$ and
 $$
  h_\FS:\quad \Hom(\,P,\,P^\perp\,)\;\times\;\Hom(\,P,\,P^\perp\,)
  \;\longrightarrow\;\C,\qquad (\,X,\,Y\,)\;\longmapsto\;\tr_P(\,X^*Y\,)
 $$
 defines the hermitean Fubini--Study metric with Riemann metric $g_\FS
 \,:=\,\Re\,h_\FS$ and K\"ahler form $\Im\,h_\FS$, where $X^*\,\in\,\Hom
 (P^\perp,P)$ denotes the adjoint of $X$ with respect to $h$. In order to
 calculate the exponential map we apply the matrix exponential $\exp:\,
 \End\,V\longrightarrow\mathbf{GL}\,V$ to the Lie algebra element
 $X\,-\,X^*\,\in\,\mathfrak{su}\,V$ corresponding to the tangent
 vector $X$
 $$
  \exp\;\pmatrix{0 &-X^*\cr X & 0}
  \;\;=\;\;
  \pmatrix{\cos\sqrt{X^*X} & * \cr
  X\,\frac{\sin\sqrt{X^*X}}{\sqrt{X^*X}} & *}
 $$
 in which the ill--defined square root never materializes, because both
 $\cos\,x$ and $\frac{\sin\,x}x$ are even power series in $x$. The Riemannian
 exponential is now covered by the matrix exponential in the isometry group
 for symmetric spaces like $\Gr_kV$, see for example \cite{s}, in consequence
 the exponential map $\exp_P:\,T_P\Gr_kV\longrightarrow\Gr_kV$ for the
 Grassmannians can be written:
 $$
  \exp_P:\quad\Hom(\,P,\,P^\perp\,)\;\longrightarrow\;\Gr_kV,
  \quad X\;\longmapsto\;\im\,\Big(\,v\,\longmapsto\,\cos\sqrt{X^*X}\,v
  \,+\,X\frac{\sin\sqrt{X^*X}}{\sqrt{X^*X}}\,v\,\Big)
 $$
 Based on this description we simply make the following guess for K\"ahler
 normal coordinates
 \begin{equation}\label{kgr}
  \knc_P:\quad\Hom(\,P,\,P^\perp\,)\;\longrightarrow\;\Gr_kV,
  \quad X\;\longmapsto\;\im\,(\,v\,\longmapsto\,v\,+\,Xv\,)
 \end{equation}
 which we will verify in due course by calculating the normal potential. In
 any case the image of these K\"ahler normal coordinates is exactly the big
 cell consisting of all subspaces transversal to $P^\perp$, since these big
 cells form the coordinate charts in the standard holomorphic atlas for
 $\Gr_kV$ the proposed map $\knc_P:\,T_P\Gr_kV\longrightarrow\Gr_kV$ is
 certainly holomorphic. Whenever the operator norm of the tangent vector
 $X\,\in\,\Hom(P,P^\perp)$ satisfies $|\!|\,X\,|\!|\,<\,\frac\pi2$, then
 the linear map $\cos\sqrt{X^*X}\,\in\,\End\,P$ is invertible and the
 resulting equality
 $$
  \im\,\Big(\,v\,\longmapsto\,\cos\sqrt{X^*X}\,v
  \,+\,X\frac{\sin\sqrt{X^*X}}{\sqrt{X^*X}}\,v\,\Big)
  \;\;=\;\;
  \im\,\Big(\,v\,\longmapsto\,v\,+\,X\frac{\tan\sqrt{X^*X}}{\sqrt{X^*X}}
  \,v\,\Big)
 $$
 implies the following explicit formulas for the difference element $K$ and
 its inverse
 \begin{equation}\label{kk}
  K^{-1}X\;\;=\;\;X\,\frac{\tan\sqrt{X^*X}}{\sqrt{X^*X}}
  \qquad\Longleftrightarrow\qquad
  KX\;\;=\;\;X\,\frac{\arctan\sqrt{X^*X}}{\sqrt{X^*X}}
 \end{equation}
 because in this way $\knc_P(\,K^{-1}X\,)\,=\,\exp_PX$. To calculate the
 pull back of the Fubini--Study metric $g_\FS$ on $\Gr_kV$ to $T_pM$ via
 $\knc_P$ we embed $\Gr_kV$ isometrically into the real vector space of
 self adjoint endomorphisms of $V$ with respect to the hermitean form $h$
 \begin{equation}\label{emb1}
  \iota:\quad\Gr_kV\;\longrightarrow\;\End_{\mathrm{herm}}V,
  \qquad \hat P\;\longmapsto\;\pr_{\hat P}
 \end{equation}
 where $\pr_{\hat P}\,=\,\pr_{\hat P}^*$ is the orthogonal projection onto
 $\hat P$, in turn the composition $\iota\circ\knc_P$ reads
 $$
  \iota\,\Big(\;\im\,(\,v\,\longmapsto\,v\,+\,X\,v\,)\;\Big)
  \;\;=\;\;
  \pmatrix{\id\cr X}\,\pmatrix{\id\,+\,X^*X}^{-1}\,\pmatrix{\id&X^*}
  \;\;=\;\;
  \pmatrix{Q & Q\,X^*\cr X\,Q & X\,Q\,X^*}
 $$
 with $Q\,:=\,(\,\id+X^*X\,)^{-1}$. In passing we observe that the very same
 calculation implies that the differential of the embedding $\iota$ in the
 chosen, but arbitrary point $P\,\in\,\Gr_kV$ is given by
 $$
  \iota_{*,\,P}:\quad
  \Hom(\,P,\,P^\perp\,)\;\longrightarrow\;\End_{\mathrm{herm}}V,\qquad
  X\;\longmapsto\;\pmatrix{0& X^*\cr X&0}
 $$
 as $Q\,=\,\id\,+\,O(X^2)$, hence $\iota$ is an isometric embedding as claimed
 provided we choose the positive definite scalar product $G(F,\hat F)\,:=\,
 \frac12\,\tr_VF\hat F$ on $\End_{\mathrm{herm}}V$. Using the standard formula
 $\delta M\,=\,-M\,(\,\delta M^{-1}\,)\,M$ for the variation of inverses we
 may calculate the differential
 \begin{eqnarray*}
  \lefteqn{(\iota\circ\knc_P)_{*,\,X}A}\qquad
  &&
  \\
  &=&
  \pmatrix{0\cr A}Q\pmatrix{\id&X^*}
  \;-\;\pmatrix{\id\cr X}Q(A^*X+X^*A)Q\pmatrix{\id&X^*}
  \;+\;\pmatrix{\id\cr X}Q\pmatrix{0&A^*}
 \end{eqnarray*}
 of the composition $\iota\,\circ\,\knc_P:\,T_P\Gr_kV\longrightarrow
 \End_{\mathrm{herm}}V$ and find after a rather lengthy calculation
 better done separately for the 9 summands and simplifying $Q\,(\id+X^*X)
 \,=\,\id$:
 \begin{eqnarray*}
  (\,\knc_P^*g_\FS\,)_X(\;A,\,B\;)
  &\stackrel!=&
  (\,(\,\iota\circ\knc_P\,)^*G\,)_X(\;A,\,B\;)
  \\
  &=&
  {\textstyle\frac12}\,\tr_V(\;(\iota\circ\knc_P)_{*,\,X}A\;
  (\iota\circ\knc_P)_{*,\,X}B\;)
  \\
  &=&
  {\textstyle\frac12}\,\tr_P(\;QA^*B\,+\,QB^*A\,-\,QA^*XQX^*B
  \,-\,QX^*AQB^*X\,)
 \end{eqnarray*}
 Note that the peculiar arrangement of the factors in the third and
 fourth summand make the result real as the trace of sums of hermitean
 matrices. Having calculated the pull back of the metric we are in the
 position to verify that the anchored holomorphic coordinates $\knc_P:\,
 T_P\Gr_kV\longrightarrow\Gr_kV$ proposed in equation (\ref{kgr}) are actually
 the unique K\"ahler normal coordinates in the point $P\,\in\,\Gr_kV$. For
 this purpose we define $\theta_P\,\in\,C^\infty(\,T_P\Gr_kV\,)$ by
 \begin{equation}\label{npgr}
  \theta_P(\,X\,)
  \;\;=\;\;
  \tr_P\,\log(\;\id\;+\;X^*X\;)
  \;\;=\;\;
  \log\,{\textstyle\det_P}(\;\id\;+\;X^*X\;)
 \end{equation}
 and observe that for given tangent vectors $A,\,B\,\in\,T_P\Gr_kV\,\cong\,
 \Hom(P,P^\perp)$
 \begin{eqnarray*}
  \frac{\partial^2}{\partial A\,\partial B}\,\theta_P(\;X\;)
  &=&
  \frac{\partial}{\partial A}\,\tr_P(\;Q\,(\,B^*X\,+\,X^*B\,)\;)
  \\
  &=&
  \tr_P(\;Q\,(\,B^*A\,+\,A^*B\,)\;-\;Q\,(\,A^*X\,+\,X^*A\,)
  \,Q\,(\,B^*X\,+\,X^*B\,)\;)
 \end{eqnarray*}
 due to the standard logarithmic derivative $\delta\log(\det M)\,=\,
 \tr(M^{-1}\delta M)$ of the determinant and the definition $Q\,:=\,
 (\,\id\,+\,X^*X\,)^{-1}$ of $Q$. Comparing this result with the formula
 for the pull back $\knc_P^*g_\FS$ of the Fubini--Study metric from $\Gr_kV$
 to $T_P\Gr_kV$ we conclude
 $$
  (\,\knc_P^*g_\FS\,)_X(\;A,\,B\;)
  \;\;=\;\;
  \frac14\,\Big(\;\frac{\partial^2}{\partial A\,\partial B}\,
  \theta_P(\;X\;)\;+\;\frac{\partial^2}{\partial IA\,\partial IB}\,
  \theta_P(\;X\;)\;\Big)
 $$
 because the troublesome terms with either none or both of $A$ or $B$
 starred drop out in averaging over $A,\,B$ and $IA,\,IB$. In consequence the
 function $\theta_P$ is some potential function for the Riemannian metric
 $\knc_P^*g_\FS$ on $T_P\Gr_kV$ and since it evidently satisfies the
 normalization constraint imposed on the unique normal potential it
 is actually equal to this potential, in turn $\knc_P:\,T_P\Gr_kV
 \longrightarrow\Gr_kV$ are K\"ahler normal coordinates as claimed.

 \pfill
 Although the formula for the difference element of the complex
 Grassmannians obtained in this way is quite explicit, it is more useful to
 recast it in terms of the curvature of $\Gr_kV$. According to our discussion
 of the Lie algebra of Killing vector fields on symmetric spaces the curvature
 of the complex Grassmannians can be calculated from the Lie algebra of
 Killing vector fields on $\Gr_kV$ by means of equation (\ref{ad}), more
 precisely we find
 $$
  R_{U,\,V}W
  \;\;=\;\;
  -\,[\,[\,U-U^*,\,V-V^*\,],\,W-W^*\,]
  \;\;=\;\;
  U\,V^*W\,-\,V\,U^*W\,-\,W\,U^*V\,+\,W\,V^*U
 $$
 for all $U,\,V,\,W\,\in\,\Hom(P,P^\perp)$ with adjoints $U^*,\,V^*,\,W^*
 \,\in\,\Hom(P^\perp,P)$. The complex structure on the Grassmannians $\Gr_kV$
 is now defined by declaring the isomorphism $T_P\Gr_kV\,\cong\,\Hom(P,
 P^\perp)$ of real vector spaces to be complex linear $IX\,:=\,iX$ and
 so the formula for the curvature becomes for the arguments $U\,=\,IX\,=\,W$
 and $V\,=\,X\,(X^*X)^k$
 $$
  (\,\ad^2IX\,)\Big(\;X\,(\,X^*X\,)^k\;\Big)
  \;\;:=\;\;
  R_{iX,\,X(X^*X)^k}\,iX
  \;\;=\;\;
  -\,4\,X\,(\,X^*X\,)^{k+1}
 $$
 for all $k\,\geq\,0$. In consequence equation (\ref{kk}) can be written
 in the form:
 \begin{equation}\label{kex}
  KX
  \;\;=\;\;
  \frac{\arctan(\,-\,\frac14\,\ad^2IX\,)^{\frac12}}
  {(\,-\,\frac14\,\ad^2IX\,)^{\frac12}}\;X
  \;\;=\;\;
  \frac{\mathrm{artanh}(\,\frac12\;\ad\,IX\,)}{\frac12\,\ad\,IX}\;X
 \end{equation}
 This calculation was the motivation for the authors to look for a proof of
 Lemma \ref{main} describing the difference element for arbitrary hermitean
 symmetric spaces by exactly this formula. In the same vein the normal
 potential for the complex Grassmannians can be rewritten in terms of the
 Lie algebra structure in order to reflect Corollary \ref{kps}:
 \begin{equation}\label{grpot}
  (\,\knc_P^*\theta_P\,)(\;X\;)
  \;\;=\;\;
  \tr_P\log(\;\id\;+\;X^*X\;)
  \;\;\stackrel!=\;\;
  g_P\Big(\;X,\;\frac{\log(\,\id\,-\,\frac14\,\ad^2IX\,)}
  {-\,\frac14\,\ad^2IX}\,X\;\Big)
 \end{equation}

 \pfill
 The second family of examples of K\"ahler manifolds discussed in this
 appendix are the real Grassmannians of oriented planes in a real vector
 space $V$, the universal covering spaces of the standard real Grassmannians
 $\Gr_2V$ of planes in $V$, elements of $\Gr^\ort_2V$ are thus $2$--dimensional
 subspaces $P\,\subset\,V$ endowed with an orientation determining a sense of
 of counterclockwise rotation. Choosing a positive definite scalar product
 $g:\,V\times V\longrightarrow\R$ we may identify the vector space tangent
 to $\Gr^\ort_2V$ in an oriented plane $P\,\in\,\Gr^\ort_2V$ with the vector
 space
 $$
  \Hom(\;P,\;P^\perp\;)\;\stackrel\cong\longrightarrow\;T_P\Gr^\ort_2V,
  \qquad X\;\longmapsto\;\left.\frac d{dt}\right|_0
  \im\,(\;\id\,+\,t\,X:\quad P\,\longrightarrow\,V\;)
 $$
 of linear maps from $P$ to its orthogonal complement $P^\perp$ and define
 the Fubini--Study metric
 $$
  g_\FS:\quad \Hom(\,P,\,P^\perp\,)\;\times\;\Hom(\,P,\,P^\perp\,)
  \;\longrightarrow\;\R,\qquad (\,X,\,Y\,)\;\longmapsto\;\tr_P(\,X^*Y\,)
 $$
 in complete analogy to the complex Grassmannians, where $X^*$ is now the
 adjoint of $X$ with respect to the scalar product $g$. In difference to
 the complex Grassmannians however $\Hom(\,P,\,P^\perp\,)$ is not a priori
 a complex vector space so that it is impossible to define an almost complex
 structure $I$ on $\Gr_2^\ort V$ simply by multiplication with $i$.
 Nevertheless every oriented plane $P\,\subset\,V$ carries a unique isometry
 $J\,\in\,\mathbf{O}(\,P,g\,)$ satisfying $J^2\,=\,-\id_P$, namely the
 rotation by $+\,90^\circ$. In turn the tangent space $T_P\Gr^\ort_2V\,=\,
 \Hom(\,P,\,P^\perp\,)$ to the Grassmannian of oriented planes in $V$ becomes
 a complex vector space by precomposing
 $$
  I:\quad \Hom(\;P,\;P^\perp\;)\;\longrightarrow\;\Hom(\;P,\;P^\perp\;),
  \qquad X\;\longmapsto\;XJ
 $$
 with $J\,\in\,\mathbf{O}(\,P,g\,)$, moreover this complex structure on
 $\Hom(P,P^\perp)$ is orthogonal with respect to the Fubini--Study metric
 $g_\FS(IX,IY)\,=\,\tr_P(J^*X^*YJ)\,=\,g_\FS(X,Y)$ due to $J^*\,=\,J^{-1}
 \,=\,-J$. Evidently no similar construction exists for the real Grassmannians
 of oriented or unoriented subspaces of $V$ of dimensions other than $2$.

 Leaving the question of integrability of the almost complex structure $I$
 on $\Gr^\ort_2V$ aside for the moment we recall that the complex bilinear
 extension of the scalar product $g$ to the complexification $V_\C\,:=\,
 V\otimes_\R\C$ of the real vector space $V$ defines the complex quadric
 $$
  Q_g(\,V\,)
  \;\;:=\;\;
  \{\;\;[\,p\,]\,\in\,\mathbb{P}\,V_\C\;\;|\;\;\;
  v\textrm{\ isotropic vector with\ }g(\,p,\,p\,)\,=\,0\;\;\;\}
  \;\;\subset\;\;
  \mathbb{P}\,V_\C
 $$ 
 of isotropic lines in $V_\C$, which is a smooth complex submanifold of
 $\mathbb{P}\,V_\C$ and thus a K\"ahler manifold itself endowed with the
 restriction of the Fubini--Study metric associated to the positive definite
 hermitean form $h:\,\overline{V_\C}\times V_\C\longrightarrow\C,\,(v,w)
 \longmapsto g(\overline v,w),$ arising from $g$ and the real structure
 on $V_\C$.

 Evidently the real and imaginary part of every vector $p\,=\,
 (\Re\,p)\,+\,i(\Im\,p)\,\in\,V_\C$ representing an isotropic line
 $[\,p\,]\,\in\,Q_g(\,V\,)$ are orthogonal vectors of the same length
 in $V$ and vice versa. In this way the Grassmannian $\Gr^\ort_2V$ of
 oriented planes in $V$ embeds canonically into complex projective space
 and becomes the quadric $Q_g(\,V\,)\,\subset\,\mathbb{P}\,V_\C$ of isotropic
 lines. More precisely the canonical embedding sends an oriented plane
 $P\,\subset\,V$ to the isotropic line spanned by the vector $p\,:=\,e_1-ie_2$
 encoding an oriented orthonormal basis $e_1,\,e_2$ for $P$:
 $$
  \iota:\quad\Gr^\ort_2V\;\longrightarrow\;Q_g(\,V\,)\;\;\subset\;\;
  \mathbb{P}\,V_\C,\qquad P\;\longmapsto\;[\,p\,]
 $$
 With $P$ being a plane in $V$ all its oriented orthonormal bases are related
 by rotations, thus the representative vector $p\,\in\,V_\C$ is only defined up
 to multiplication by an element of $S^1$. Independent of this $S^1$--ambiguity
 in the choice of $p\,\in\,V_\C$ the following identities hold true
 \begin{equation}\label{tid}
  \begin{array}{ccl}
   g(\,\overline p,\,Fp\,)\;
   &=&
   \tr_PF
   \\[2pt]
   |\;g(\,p,\,Fp\,)\;|^2
   &=&
   \tr^2_PF\;-\;4\,\mathrm{det}_PF
   \;\;=\;\;
   2\,\tr_P(\,F^2\,)\;-\;\tr_P^2F
  \end{array}
 \end{equation}
 for every symmetric endomorphism $F\,\in\,\End\,P$ as the reader may easily
 verify using the matrix coefficients of $F$ in the orthonormal basis $e_1,
 \,e_2$. In terms of the identification of the vector space tangent to
 $\mathbb{P}\,V_\C$ in the point $\iota(P)\,=\,[\,p\,]$ with the vector
 space of complex linear maps from the line $\C p$ to its orthogonal
 complement $\{\,p\,\}^\perp$ with respect to the hermitean form $h$
 the differential of $\iota$ in an oriented plane $P\,\in\,\Gr^\ort_2V$
 can be written
 $$
  \iota_{*,\,P}:\quad\Hom(\;P,\;P^\perp\;)\;\longrightarrow\;
  \Hom(\;\C\,p,\,\{\,p\,\}^\perp\;),\qquad X\;\longmapsto\;
  \left.X\right|_{\C p}
 $$
 in fact $\iota_{*,\,P}:\,\left.\frac d{dt}\right|_0\im(\,\id\,+\,tX\,)
 \;\longmapsto\;\left.\frac d{dt}\right|_0[\,p\,+\,t\,Xp\,]$. In particular
 the embedding $\iota$ is actually an holomorphic embedding with
 $\iota_{*,\,P}(IX)\,=\,\left.XJ\right|_{\C p}\,=\,\left.iX\right|_{\C p}$
 due to $Jp\,=\,ip$ so that the almost complex structure $I$ on $\Gr^\ort_2V$
 is necessarily integrable. However $\iota$ is not an isometric embedding for
 the Fubini--Study metric $g^{\mathbb{P}\,V_\C}$ on the target, to be
 precise
 $$
  g^{\mathbb{P}\,V_\C}_{[\,p\,]}
  (\;\left.X\right|_{\C p},\;\left.Y\right|_{\C p}\;)
  \;\;:=\;\;
  \frac{h(\,Xp,\,Yp\,)}{h(\,p,\,p\,)}
  \;\;=\;\;
  \frac12\,g(\,X\overline{p},\,Yp\,)
  \;\;=\;\;
  \frac12\,\tr_P(\,X^*Y\,)
 $$
 equals half the Fubini--Study metric $\frac12\,g_\FS$ on $\Gr^\ort_2V$.
 Interestingly this observation is sufficient to construct the K\"ahler
 normal coordinates for the real Grassmannian $\Gr^\ort_2V$ by setting:
 $$
  \iota\,\circ\,\knc_p:\quad\Hom(\;P,\,P^\perp\;)\;\longrightarrow\;
  Q_g(\,V\,),\qquad X\;\longmapsto\;[\,p\,+\,X\,p\,-\,{\textstyle\frac14}
  \,g(Xp,Xp)\,\overline{p}\,]
 $$
 The $S^1$--ambiguity of the isotropic vector $p\,\in\,V_\C$
 representing the oriented plane $P\,\in\,\Gr^\ort_2V$ has no bearance
 on the complex line spanned by the vector $p+Xp-\frac14g(Xp,Xp)\overline{p}$,
 moreover the image vector is isotropic due to $Xp\,\in\,P^\perp\otimes_\R\C$
 and $g(\overline p,p)\,=\,2$. Last but not least we observe that the map
 $\iota\circ\knc_P:\,\Hom(P,P^\perp)\longrightarrow Q_g(\,V\,)$ is covered
 by the complex quadratic polynomial $\Hom(P,P^\perp)\longrightarrow V_\C
 \setminus\{\,0\,\},\,X\longmapsto p+Xp-\frac14g(Xp,Xp)\overline{p},$ and
 in turn is holomorphic, recall that after all we have
 $(IX)p\,:=\,XJp\,=\,i(Xp)$.

 The implicitly defined map $\knc_P:\,\Hom(P,P^\perp)\longrightarrow
 \Gr^\ort_2V$ arising from $\iota\circ\knc_P$ is thus well--defined and
 holomorphic, to verify that $\knc_P$ are the K\"ahler normal coordinates
 it thus suffices to find a local K\"ahler potential for $g_\FS$ satisfying
 the normalization condition of Definition \ref{knc}. For this purpose we
 pull back the K\"ahler normal potential of $\mathbb{P}\,V_\C$ in the point
 $\iota(\,P\,)\,=\,[\,p\,]$ back to $\Hom(P,P^\perp)$ and multiply by $2$
 to account for the homothety $(\iota\circ\knc_P)^*g^{\mathbb{P}\,V_\C}
 \,=\,\frac12\,g_\FS$. According to our discussion (\ref{grpot}) of the
 complex Grassmannians
 $$
  \theta^{\mathbb{P}\,V_\C}_{[\,p\,]}(\,[\,p\,+\,q\,]\,)
  \;\;=\;\;
  \log\Big(\;1\;+\;\frac{h(\,q,\,q\,)}{h(\,p,\,p\,)}\;\Big)
 $$
 is the K\"ahler normal potential for $g^{\mathbb{P}\,V_\C}$ whenever
 $q\,\in\,\{\,p\,\}^\perp$ holds true, in turn we conclude
 \begin{eqnarray*}
  2\,(\,\iota\circ\knc_P\,)^*\theta^{\mathbb{P}\,V_\C}(\;X\;)
  &=&
  2\;\theta^{\mathbb{P}\,V_\C}_{[\,p\,]}
  (\;[\,p\,+\,X\,p\,-\,{\textstyle\frac14}\,g(Xp,Xp)\,\overline{p}\,]\;)
  \\
  &=&
  2\,\log\Big(\;1\,+\,{\textstyle\frac12}\,g(X\overline p,Xp)\,+\,
  {\textstyle\frac1{16}}\,|\,g(Xp,Xp)\,|^2\;\Big)
  \\
  &=&
  2\,\log\Big(\;1\,+\,{\textstyle\frac12}\,\tr_P(\,X^*X\,)
  \,+\,{\textstyle\frac18}\,\tr_P(\,X^*X\,)^2\,-\,{\textstyle\frac1{16}}
  \,\tr^2_P(\,X^*X\,)\;\Big)
 \end{eqnarray*}
 by using $h(p,p)\,=\,2$ and the identities (\ref{tid}). The result is
 actually a power series invariant under $\Der_I$ due to $(IX)^*(IX)\,=\,
 -J^*(X^*X)J$ and thus satisfies the normalization constraint required by
 Definition \ref{knc}, in consequence $\knc_P:\,\Hom(P,P^\perp)\longrightarrow
 \Gr^\ort_2V$ are the unique K\"ahler normal coordinates of $\Gr^\ort_2V$
 in the point $P\,\in\,\Gr^\ort_2V$.
 
 \pfill
 In order to calculate the difference elements $K$ and $K^{-1}$ for the
 real Grassmannians $\Gr^\ort_2V$ we still have to compare the formula
 for $\knc_P$ established above with the analoguous formula for the Riemannian
 exponential map. The calculations can be done in complete analogy to the
 case of the complex Grassmannians, because the Riemannian exponential is
 still covered by a suitable version of the matrix exponential, and the
 final results reads:
 \begin{equation}\label{exp2}
  \iota(\,\exp_PX\,)
  \;\;=\;\;
  \Big[\;(\cos\sqrt{X^*X})\,p\,+\,
  X\,\frac{\sin\sqrt{X^*X}}{\sqrt{X^*X}}\,p\;\Big]
 \end{equation}
 Solving the equation $\exp_PY\,=\,\knc_PX$ with respect to the tangent
 vector $Y\,=\,KX$ for a given argument vector $X\,\in\,\Hom(P,P^\perp)$
 is thus equivalent to solving the equation
 $$
  \Big[\;(\cos\,\sqrt{Y^*Y})\,p\;+\;
  Y\,\frac{\sin\,\sqrt{Y^*Y}}{\sqrt{Y^*Y}}\,p\;\Big]
  \;\;=\;\;
  \Big[\;p\;+\;X\,p\;-\;\frac14\,g(Xp,Xp)\,\overline{p}\;\Big]
 $$
 for points in $\mathbb{P}\,V_\C$, which we may decouple into two
 independent equations
 \begin{equation}\label{deco}
  \tau\,\Big(\,(\,\cos\,\sqrt{Y^*Y}\,)\,p\,\Big)
  \;\;=\;\;
  p\;-\;\frac14\,g(Xp,Xp)\,\overline{p}
  \qquad\quad
  \tau\,\Big(\,Y\,\frac{\sin\,\sqrt{Y^*Y}}{\sqrt{Y^*Y}}\,p\,\Big)
  \;\;=\;\;
  X\,p
 \end{equation}
 by taking the orthgonal decomposition $V\,=\,P\,\oplus\,P^\perp$ into
 account and introducing a non--zero slack variable $\tau\,\in\,\C^*$.
 Applying the complex bilinear scalar product $g$ with $\overline{p}$
 to the first equation and using the identities (\ref{tid}) we find that
 $\tau\,\in\,\R^+$ is actually positive
 $$
  \tau\,(\,\tr_P\cos\,\sqrt{Y^*Y}\,)
  \;\;=\;\;
  g\Big(\;\overline{p},\;p\,-\,\frac14\,g(Xp,Xp)\,\overline{p}\;\Big)
  \;\;=\;\;
  2
 $$
 and thus:
 \begin{equation}\label{knc2}
  X
  \;\;=\;\;
  K^{-1}Y
  \;\;=\;\;
  \frac2{\tr_P(\,\cos\,\sqrt{Y^*Y}\,)}
  \,Y\,\frac{\sin\,\sqrt{Y^*Y}}{\sqrt{Y^*Y}}
 \end{equation}
 Unluckily we wanted to solve the equation $\exp_PY\,=\,\knc_PX$ for $Y$
 and not for $X$, hence we regress to the decoupled equations (\ref{deco})
 and read the identities (\ref{tid}) backwards to obtain
 \begin{eqnarray*}
  \tau^2\,(\,\tr_P\sin^2\sqrt{Y^*Y}\,)
  &=&
  g\Big(\;\tau\,(\,Y\,\frac{\sin\,\sqrt{Y^*Y}}{\sqrt{Y^*Y}}\,)\,\overline p,
  \;\tau\,(\,Y\,\frac{\sin\,\sqrt{Y^*Y}}{\sqrt{Y^*Y}}\,)\,p\;\Big)
  \\
  &=&
  g(\;X\,\overline p,\;X\,p\;)
  \\[8pt]
  \tau^2\,(\,\tr_P\cos^2\sqrt{Y^*Y}\,)
  &=&
  g(\;\overline p\,-\,{\textstyle\frac14}\,g(X\overline p,X\overline p)
  \,p,\;p\,-\,{\textstyle\frac14}\,g(Xp,Xp)\,\overline p\;)
  \\[3pt]
  &=&
  2\;+\;{\textstyle\frac18}\,|\,g(\,Xp,\,Xp\,)\,|^2
  \;\;=\;\;
  2\;+\;{\textstyle\frac14}\,\tr_P(\,X^*X\,)^2
  \;-\;{\textstyle\frac18}\,\tr_P^2(\,X^*X\,)
 \end{eqnarray*}
 With $\tr_P(\,\sin^2\sqrt{Y^*Y}\,+\,\cos^2\sqrt{Y^*Y}\,)\,=\,2$
 and $g(X\overline p,Xp)\,=\,\tr_P(\,X^*X\,)$ we can solve for $\tau$:
 $$
  \tau
  \;\;=\;\;
  \tau(\;X\;)
  \;\;=\;\;
  \Big(\;1\;+\;{\textstyle\frac12}\,\tr_P(X^*X)
  \;+\;{\textstyle\frac18}\,\tr_P(X^*X)^2
  \;-\;{\textstyle\frac1{16}}\,\tr^2_P(X^*X)\;\Big)^{\frac12}
 $$
 Incidentally we observe that this expression is exactly the argument
 of the logarithm in the formula $(\knc_P^*\theta^\FS_P\,)(\,X\,)\,=\,
 4\,\log\,\tau(X)$ for the K\"ahler normal potential of the real
 Grassmannian $\Gr^\ort_2V$. In any case we conclude by solving the
 second of the equations (\ref{deco}) that:
 \begin{equation}\label{knc3}
  Y
  \;\;=\;\;
  KX
  \;\;=\;\;
  \frac{X}{\tau(\,X\,)}\,\frac{\arcsin\sqrt{\frac{X^*X}{\tau^2(\,X\,)}}}
  {\sqrt{\frac{X^*X}{\tau^2(\,X\,)}}}
 \end{equation}
 Although this formula has little to no resemblance to the formula of
 Theorem \ref{main}, it can be shown by explicit power series expansion
 that both formulas actually amount to the same.

 \pfill
 In the last part of this appendix we want to discuss K\"ahler normal
 coordinates for a particularly interesting family of hermitean symmetric
 spaces: The twistor spaces of orthogonal complex structures on real
 vector spaces of even dimension and the closely related twistor spaces
 of quaternionic linear orthogonal complex structures on quaternionic
 vector spaces. Starting with the former we consider a real vector space
 $V$ of even dimension $2n$ endowed with a positive definite scalar product
 $g:\,V\times V\longrightarrow\R$ and define its twistor space as:
 $$
  \T(\,V,\,g\,)
  \;\;:=\;\;
  \{\;\;J\,\in\,\End\,V\;\;|\;\;J\textrm{\ orthogonal endomorphism with\ }
  J^2\,=\,-\id_V\;\;\}
 $$
 Since an endomorphism squaring to $-\id_V$ is orthogonal, if and only if
 it is skew symmetric with respect to $g$, the twistor space is actually
 a submanifold $\T(V,g)\,\subset\,\mathfrak{so}(V,g)$ of the Lie algebra
 of skew symmetric endomorphisms of $V$. In particular we may identify the
 tangent space $T_J\T(V,g)\,\subset\,\mathfrak{so}(V,g)$ of the twistor
 space in a point $J\,\in\,\T(V,g)$ with the vector space
 $$
  \Sigma^1_\skew(\;J\;)
  \;\;:=\;\;
  \{\;\;X\,\in\,\mathfrak{so}(\,V,\,g\,)\;\;|\;\;\;
  X\textrm{\ skew symmetric and\ }XJ\,+\,JX\,=\,0\;\;\;\}
 $$
 of skew symmetric endomorphisms $X$ of $V$ anticommuting with $J$ by
 differentiation:
 $$
  T_J\T(\;V,\,g\;)\;\stackrel\cong\longrightarrow\;
  \Sigma^1_\skew(\,J\,),\qquad\left.\frac{d}{dt}\right|_0
  J_t\;\longmapsto\;\dot J_0
 $$
 In fact $\dot J_0\,\in\,\Sigma^1_\skew(\,J\,)$ anticommutes with
 $J\,=\,J_0$, because $J_t^2\,=\,-\id_V$ for all $t$. The Riemannian
 metric of choice on $\T(V,g)\,\subset\,\mathfrak{so}(V,g)$ is simply
 the restriction of the standard scalar produce $g(X,Y)\,=\,-\frac12\,
 \tr_V(XY)$ on $\mathfrak{so}(V,g)$, moreover the almost complex structure
 on $\T(V,g)$ in a point $J\,\in\,\T(V,g)$ is the right multiplication with
 $J$ on the tangent space
 $$
  I_J:\quad\Sigma^1_\skew(\,J\,)\;\longrightarrow\;\Sigma^1_\skew(\,J\,),
  \qquad X\;\longmapsto\;X\,J
 $$
 because $XJ$ still is skew and anticommutes with $J$. Since $\T(V,g)\,
 \subset\,\mathfrak{so}(V,g)$ is actually a union of two adjoint orbits
 the Riemannian exponential is easily seen to be covered by the matrix
 exponential for the Lie group $\mathbf{SO}(\,V,\,g\,)$, the only subtlety
 here is that the Lie algebra element corresponding to the tangent vector
 $X\,\in\,\Sigma^1_\skew(\,J\,)$ is not $X$ itself, but $-\frac12\,XJ$, due
 to the identity $[\,-\frac12\,XJ,\,J\,]\,=\,X$. In this way we find the
 explicit formula
 $$
  \exp_J:\quad\Sigma^1_\skew(\,J\,)\;\longrightarrow\;\T(\;V,\,g\;),\qquad
  X\;\longmapsto\;e^{-\,\frac12\,XJ}\,J\,e^{+\,\frac12\,XJ}
 $$
 for the Riemannian exponential of the twistor space $\T(V,g)$, and since
 $X$ anticommutes with $J$ we may simplify this to read $\exp_JX\,=\,e^{-XJ}J
 \,=\,(\cosh\,X)J\,+\,(\sinh\,X)$ using the observation $(-XJ)^2\,=\,X^2$.
 In order to establish the integrability of the almost complex structure $I$
 and calculate the K\"ahler normal coordinates it turns out to be convenient
 to identify the twistor space $\T(V,g)$ with the Grassmannian of Lagrangian
 subspaces for the complex bilinear extension of the scalar product $g$ on
 $V$ to $V_\C\,=\,V\otimes_\R\C$, which is a complex submanifold of the
 Grassmannian $\Gr_nV_\C$ of $n$--dimensional subspaces of $V_\C$ defined by:
 $$
  \LGr_nV_\C
  \;\;=\;\;
  \{\;\;L\;\in\;\Gr_nV_\C\;\;|\;\;\;
  L\textrm{\ isotropic with respect to\ }g\;\;\;\}
 $$
 Explicitly the diffeomorphism between $\T(V,g)$ and $\LGr_nV_\C$ reads
 $$
  \iota:\quad\T(\;V,\,g\;)\;\stackrel\cong\longrightarrow\;\LGr_nV_\C,
  \qquad J\;\longmapsto\;V^{1,\,0}_J
 $$
 where $V^{1,0}_J\,\subset\,V_\C$ is the eigenspace for the eigenvalue $+i$
 of the complex linear extension of the orthogonal complex structure $J\,\in
 \,\T(V,g)$ to $V_\C$. In passing we remark that $\LGr_nV_\C$ can also be
 realized as the quadratic projective variety of pure spinors. 

 The argument demonstrating the surjectivity of $\iota$ uses the canonical
 real structure on the complexification $V_\C\,=\,V\otimes_\R\C$ in the
 form of an involution $L\,\longmapsto\,\overline L$ of $\LGr_nV_\C$
 satisfying $L\,\cap\,\overline L\,=\,\{\,0\,\}$, after all $g$ is positive
 definite on the real subspace $V\,\subset\,V_\C$. We may thus associate
 to every Lagrangian subspace $L\,\in\,\LGr_nV_\C$ the endomorphism
 $J\,:=\,i\pr_L\,-\,i\pr_{\overline L}$ of $V_\C\,=\,L\,\oplus\,\overline L$,
 which commutes by construction with the real structure and thus comes from
 a complex structure $J\,\in\,\End\,V$ of the underlying real vector space
 $V$. With its eigenspaces $L$ and $\overline L$ being isotropic subspaces
 $J$ is automatically skew symmetric and thus orthogonal with respect to
 $g$, moreover we find $V^{1,0}_J\,=\,L$ so that the orthogonal complex
 structure $J\,\in\,\T(V,g)$ is a preimage of the subspace $L$ we
 started with.

 Let us now turn to the calculation of the differential of $\iota$ in order
 to show that it is an holomorphic embedding $\T(V,g)\longrightarrow\Gr_nV_\C$.
 For a curve $t\longmapsto J_t$ of orthogonal complex structures on $V$
 representing a tangent vector $X\,=\,\left.\frac d{dt}\right|_0J_t\,\in\,
 \Sigma^1_\skew(\,J\,)$ in $J\,=\,J_0$ we find
 $$
  \iota_{*,\,J}\Big(\;\left.\frac d{dt}\right|_0J_t\;\Big)
  \;\;=\;\;
  \left.\frac d{dt}\right|_0\im\Big(\;\;V^{0,1}_J\,\longrightarrow\,
  V^{0,1}_{J_t},\;v\,\longmapsto\,\frac12\,(\,v\,-\,iJ_tv\,)\;\Big)
  \;\;\widehat=\;\;
  -\,\frac i2\,\left.\frac d{dt}\right|_0J_t
 $$
 and hence conclude that the differential of $\iota:\,\T(V,g)\longrightarrow
 \Gr_nV_\C$ reads:
 $$
  \iota_{*,\,J}:\quad\Sigma^1_\skew(\,J\,)\;\longrightarrow\;
  \Hom(\;V^{1,0}_J,\;V^{0,1}_J\;),\qquad X\;\longmapsto\;-\,
  \frac i2\left.X\right|_{V^{1,0}_J}
 $$
 In particular $\iota_{*,\,J}(\,IX\,)\,=\,-\frac i2\!\left.XJ
 \right|_{V^{1,0}_J}\,=\,i\,\iota_{*,\,J}(\,X\,)$ due to $J\,=\,+i$
 on $V^{1,0}_J$ making $\iota$ a holomorphic embedding as claimed,
 necessarily then the almost complex structure $I$ on the twistor
 space $\T(V,g)$ is integrable. Moreover the pull back of the
 Fubini--Study metric on $\Gr_nV_\C$ via $\iota$ equals $\frac14$
 times the chosen Riemannian metric on $\T(V,g) \,\subset\,
 \mathfrak{so}(V,g)$, because
 \begin{equation}\label{factor}
  (\,\iota^*g_\FS\,)_J(\;X,\;Y\;)
  \;\;=\;\;
  \Re\;\tr_{V^{1,0}_J}(\;(-{\textstyle\frac i2}X)^*(-{\textstyle\frac i2}Y)\;)
  \;\;=\;\;
  -\,{\textstyle\frac18}\,\tr_V(\,XY\,)
 \end{equation}
 recall here that $X^*\,=\,-X$ is skew symmetric and that
 $\overline{\tr_{V^{1,0}_J}F}\,=\,\tr_{V^{0,1}_J}F$ for every endomorphism
 $F\,\in\,\End\,V$ of the underlying real vector space $V$.

 For the next step we need to discuss another characterization of symmetric
 spaces as manifolds endowed with a binary operation $*$ generalizing in a
 sense the multiplication of a Lie group, \cite{br} is a very good reference
 for this point of view. The binary operation $*:\,\Gr_nV_\C\times\Gr_nV_\C
 \longrightarrow\Gr_nV_\C$ associated to the complex Grassmannian $\Gr_nV_\C$
 is most easily defined in terms of a suitable modification of the embedding
 (\ref{emb1}) used previously
 $$
  \mathrm{inv}:\quad\Gr_nV_\C\;\longrightarrow\;\End_{\mathrm{herm}}V_\C,
  \qquad\hat P\;\longmapsto\;2\,\pr_{\hat P}\;-\;\id
 $$
 which identifies $\Gr_nV_\C$ with the set of self adjoint involutions
 of $V_\C$ with zero trace, to wit
 $$
  \mathrm{inv}(\;P\,*\,\hat P\;)
  \;\;:=\;\;
  \mathrm{inv}(\,P\,)\,\circ\,\mathrm{inv}(\,\hat P\,)\,\circ\,
  \mathrm{inv}(\,P\,)
 $$
 where the right hand side is still a self adjoint involution with zero
 trace. For obvious reasons the composition $\mathrm{inv}\circ\iota:\,
 \T(V,g)\longrightarrow\End_{\mathrm{herm}}V_\C$ maps $J$ to $\mathrm{inv}
 (\,V^{1,0}_J\,)\,=\,-iJ$ so that $\iota$ is actually a homomorphism of
 symmetric spaces in the sense:
 $$
  \iota(\,J\,)\;*\;\iota(\,\hat J\,)
  \;\;=\;\;
  \mathrm{inv}^{-1}(\;(-iJ)\,(-i\hat J)\,(-iJ)\;)
  \;\;=\;\;
  \iota(\,J\,*\,\hat J\,)
  \;\;:=\;\;
  \iota(\,J\,\hat J^{-1}\,J\,)
 $$
 Like every other homomorphism of symmetric spaces $\iota$ is thus a totally
 geodesic embedding, in turn Corollary \ref{sub} tells us that the K\"ahler
 normal coordinates for the twistor space $\T(V,g)$ are simply the restriction
 of the K\"ahler normal coordinates of the complex Grassmannian:
 $$
  \iota\,\circ\,\knc_J:\quad\Sigma^1_\skew(\,J\,)\;\longrightarrow\;\LGr_nV_\C,
  \qquad X\;\longmapsto\;\knc_{V^{1,0}_J}(\,\iota_{*,\,J}X\,)
 $$
 In fact it is easily verified that the image subspace
 $$
  \knc_{V^{1,0}_J}(\,-\,{\textstyle\frac i2}X\,)
  \;\;=\;\;
  \im\;\Big(\;\;\id-{\textstyle\frac i2}X\;:\;\;
  V^{1,0}_J\;\longrightarrow\;V_\C\;\;\Big)
  \;\;\stackrel!\in\;\;
  \LGr_nV_\C
 $$
 is actually a Lagrangian subspace of $V_\C$. It remains to find the orthogonal
 complex structure corresponding to this image subspace. For this purpose we
 recall that the square $X^2$ of the skew symmetric endomorphism $X\,\in\,
 \Sigma^1_\skew(\,J\,)$ is diagonalizable with non--positive eigenvalues
 so that the endomorphisms $\id\pm\frac12XJ$ are always invertible due to:
 $$
  4\,(\,\id\,+\,{\textstyle\frac12}\,XJ\,)\,
  (\,\id\,-\,{\textstyle\frac12}\,XJ\,)
  \;\;=\;\;
  4\,-\,(XJ)^2
  \;\;=\;\;
  4\,-\,X^2
 $$
 Hence we may write the image of $\id-\frac i2X:\,V^{1,0}_J\longrightarrow
 V_\C$ as the image of the endomorphism
 \begin{eqnarray*}
  (\,\id\,-\,{\textstyle\frac i2}\,X\,)\;\circ\;\pr_{V^{1,0}_J}
  &=&
  {\textstyle\frac12}\,(\,\id\,-\,{\textstyle\frac i2}\,X\,)
  \,(\,\id\,-\,i\,J\,)
  \\
  &=&
  {\textstyle\frac12}\,(\,\id\,-\,i\,[\,(J+{\textstyle\frac12}X)\,
  (\id-{\textstyle\frac12}XJ)^{-1}\,]\,)\,\circ\,(\,\id\,-\,
  {\textstyle\frac12}\,X\,J\,)
  \\[2pt]
  &=&
  \pr_{V^{1,0}_{\knc_JX}}\;\circ\;(\,\id\,-\,{\textstyle\frac12}\,X\,J\,)
 \end{eqnarray*}
 of $V_\C$ where the orthogonal complex structure $\knc_JX\,\in\,\T(V,g)$
 is given explicitly by:
 \begin{equation}
  \knc_JX
  \;\;=\;\;
  4\,\frac{(J+\frac12X)\,(\id+\frac 12XJ)}{4\;-\;X^2}
  \;\;=\;\;
  \frac{4\,+\,X^2}{4\,-\,X^2}\,J\;+\;\frac{4\,X}{4\,-\,X^2}
 \end{equation}
 The reader may find it amusing to verify the slightly surprising statement
 $\knc_JX\,\in\,\T(V,g)$ directly. Comparing this result with the formula
 for the exponential map of the twistor space
 $$
  \exp_J(\,KX\,)
  \;\;=\;\;
  (\,\cosh\,KX\,)\,J\;+\;(\,\sinh\,KX\,)
  \;\;\stackrel!=\;\;
  \frac{4\,+\,X^2}{4\,-\,X^2}\,J\;+\;\frac{4\,X}{4\,-\,X^2}
  \;\;=\;\;
  \knc_JX
 $$
 we find $\exp(\,KX\,)\,=\,\cosh\,KX\,+\,\sinh\,KX
 \,=\,\frac{(2+X)\,(2+X)}{(2+X)\,(2-X)}$ and thus conclude:
 \begin{equation}
  KX
  \;\;=\;\;
  \log\;\frac{2\,+\,X}{2\,-\,X}
  \;\;=\;\;
  2\,\mathrm{artanh}(\,{\textstyle\frac12}X\,)
  \qquad\Longleftrightarrow\qquad
  K^{-1}X
  \;\;=\;\;
  2\,\tanh(\,{\textstyle\frac12}X\,)
 \end{equation}
 Last but not least the normal potential of the twistor space equals the
 restriction of the K\"ahler potential of the complex Grassmannian $\Gr_nV_\C$
 to the K\"ahler submanifold $\LGr_nV_\C\,\cong\,\T(V,g)$
 $$
  \theta_J(\;X\;)
  \;\;=\;\;
  4\;\tr_{V^{1,0}_J}\log\Big(\;\id\;+\;
  (-{\textstyle\frac i2}X)^*(-{\textstyle\frac i2}X)\;\Big)
  \;\;=\;\;
  2\;\tr_V\log\Big(\;\id\;-\;{\textstyle\frac14}\,X^2\;\Big)
 $$
 where the factor $4$ accounts for the homothety $\iota^*g_\FS\,=\,
 \frac14\,g$. Needless to say this formula for the potential is compatible
 with the description $g_J(X,X)\,=\,-\frac12\,\tr_V(\,X^2\,)$ of the
 Riemannian metric induced on the union $\T(V,g)\,\subset\,\mathfrak{so}(V,g)$
 of two adjoint orbits.

 \pfill
 The discussion of the last family of hermitean symmetric spaces considered
 in this appendix, the quaternionic twistor spaces, can be kept very short,
 because we may interprete a quaternionic vector space endowed with a positive
 definite quaternionic hermitean form $h$ as a real vector space $V$ of
 dimension $4n$ endowed with a scalar multiplication $\mathbb{H}\times V
 \longrightarrow V$ by quaternions and a compatible positive definite scalar
 product $g\,=\,\Re\,h$ in the sense $g(qv,w)\,=\,g(v,\overline q w)$ for all
 $q\,\in\,\mathbb{H}$ and all $v,\,w\,\in\,V$. Under this reinterpretation
 the quaternionic twistor space of all orthogonal complex structures on
 $V$ commuting with $\mathbb{H}$
 $$
  \T^{\mathbb{H}}(\,V,\,g\,)
  \;\;:=\;\;
  \{\;\;J\,\in\,\End_{\mathbb{H}}V\;\;|\;\;J\textrm{\ orthogonal endomorphism
  with\ }J^2\,=\,-\id_V\;\;\}
 $$
 becomes a symmetric subspace $\T^{\mathbb{H}}(V,g)\,\subset\,\T(V,g)$
 of the real twistor space associated to $V$, because $J*\hat J$ commutes
 with the scalar multiplication by $\mathbb{H}$, whenever so do $J,\,\hat J
 \,\in\,\T^{\mathbb{H}}(V,g)$. In consequence the quaternionic twistor space
 is a totally geodesic K\"ahler submanifold and all the formulas pertaining
 to $\T(V,g)$ apply verbatim to $\T^{\mathbb{H}}(V,g)$.
\end{document}